\documentclass[11pt]{amsart}
\usepackage{graphicx}
\usepackage{amsmath,amsbsy}
\usepackage{amssymb, enumerate,color}
\usepackage{comment}
\usepackage{mathrsfs}      
\usepackage{helvet}         
\usepackage{courier}        
\usepackage{type1cm}        
\usepackage{esint}
\usepackage{subfig}

\includecomment{versiona}

\usepackage{tikz}
\usetikzlibrary{shapes.geometric}

\definecolor{ForestGreen} {cmyk}{0.91,0,0.88,0.12}
\definecolor{LemonChiffon}{rgb}{1.,0.98,0.8}
\definecolor{LightCoral}  {rgb}{0.94,0.5,0.5}
\definecolor{MediumYellow}{rgb}{1.,1.,0.50}
\definecolor{PaleGreen}   {rgb}{0.6,0.98,0.6}
\definecolor{PineGreen}   {cmyk}{0.92,0,0.59,0.25}
\definecolor{SeaGreen}    {cmyk}{1,0,0.70,0}
\definecolor{Salmon}      {cmyk}{0,0.53,0.38,0}

\newcommand{\rosso}[1]{{\color{black} #1}} 

\newcommand{\rv}[1]{{\color{black} #1}}
\newcommand{\rvf}[1]{{\color{black} #1}}

\newcommand{\rev}[1]{{\color{black} #1}}

\long\def\COMMENT#1{}
\long\def\NOTE#1{}
\long\def\TESI#1{} 

\long\def\QUESTION#1{}

\long\def\FABIO#1{}

\long\def\COMMENT#1{}
\long\def\NOTE#1{}

\usepackage{geometry,calc}
\newtheorem{thm}{Theorem}[section]
\newtheorem{cor}[thm]{Corollary}
\newtheorem{lem}[thm]{Lemma}
\newtheorem{prop}[thm]{Proposition}
\theoremstyle{definition}

\newtheorem{pb}[thm]{Problem}

\theoremstyle{remark}
\newtheorem{rem}[thm]{Remark}
\numberwithin{equation}{section}

\let\div\undefined\DeclareMathOperator{\div}{div} 
\DeclareMathOperator{\grad}{\nabla} 
\DeclareMathOperator*{\Grad}{\boldsymbol\nabla}
\DeclareMathAlphabet\mathbfcal{OMS}{cmsy}{b}{n}

\newcommand\elementvem{K} 
\newcommand{\edge}{e}

\newcommand{\Poly}[1]{\mathbb{P}_{#1}}
\newcommand{\Polydue}[1]{\pmb{\mathbb{P}}_{#1}(\elementvem)}
\newcommand{\PolydueM}[1]{\pmb{\mathbb{P}}_{#1}}
\newcommand{\Pinabla}[1]{\mathsf{\Pi}^\nabla_{#1}}
\newcommand{\PinablaG}[1]{\Pinabla{#1}}
\newcommand{\Pizero}[1]{\mathsf{\Pi}^0_{#1}}
\newcommand{\PizeroG}[1]{\Pizero{#1}}
\newcommand{\Pinablakpu}{\rvf{\Pinabla{k+1}}}

\newcommand{\Qnabla}[1]{\mathsf{Q}^\nabla_{#1}}
\newcommand\Qnablakpu{\rvf{\Qnabla{k+1}}}
\newcommand{\Qzero}[1]{\mathsf{Q}^0_{#1}}

\newcommand\Skpu{\rvf{S_{k+1}^\elementvem}}

\newcommand{\bK}{\partial \elementvem}

\newcommand\bolla{\mathfrak{b}}
\newcommand\bollatest{\mathfrak{d}}

\newcommand{\x}{\mathbf{x}}

\newcommand{\ww}{\mathbf{w}}
\newcommand{\m}{\mathbf{m}}

\renewcommand{\u}{\mathbf{u}}
\renewcommand{\v}{\mathbf{v}}
\newcommand{\utilde}{\widetilde{\u}_h}
\newcommand{\vtilde}{\widetilde{\v}_h}

\newcommand{\f}{\mathbf{f}}
\newcommand{\q}{\mathbf{q}}
\newcommand{\ppsi}{\boldsymbol{\psi}}

\newcommand{\V}{\mathbf{V}}
\newcommand{\VV}{\mathbb{V}}

\newcommand{\harmonic}[1]{\mathbb{H}_{#1}}

\newcommand{\bubble}[1]{\mathcal{B}_{#1}(\elementvem)}

\newcommand{\n}{\mathbf{n}}

\newcommand{\bH}{\mathbf{H}}

\newcommand{\HunozeroK}{\mathrm{H}^1_0(\elementvem)}
\newcommand{\dueHunozeroK}{\mathbf{H}^1_0(\elementvem)}
\newcommand{\HunoK}{\mathrm{H}^1(\elementvem)}

\newcommand{\HunozeroO}{\mathrm{H}^1_0(\Omega)}

\newcommand{\dueHunozeroO}{\bH^1_0(\Omega)}
\newcommand\Ldo{\mathrm{L}^2_0(\Omega)} 

\newcommand\ds{\mathrm{ds}} 
\newcommand\dx{\mathrm{d}\x} 

\newcommand\monomials{\mathcal{P}}
\newcommand\monomialsdue{\mathbfcal{P}}
\newcommand\identity{\mathrm{Id}}

\newcommand{\vem}[1]{W_{#1}}
\newcommand{\enlarged}[1]{Z_{#1}}
\newcommand{\enhanced}[1]{\widetilde{W}_{#1}}

\newcommand{\enhancedunoS}{\widetilde{\mathbf{V}}_k(\elementvem)}
\newcommand{\enhancedunoG}{\widetilde{\mathbf{V}}_k}

\newcommand{\bubbleSp}{\rvf{{\mathcal{B}}_{k+1}^{k-2}(\elementvem)}}
\newcommand{\bubbleS}{\rvf{\pmb{\mathcal{B}}_{k+1}^{k-2}(\elementvem)}}

\DeclareFontFamily{U}{matha}{\hyphenchar\font45}
\DeclareFontShape{U}{matha}{m}{n}{
	<-6> matha5 <6-7> matha6 <7-8> matha7
	<8-9> matha8 <9-10> matha9
	<10-12> matha10 <12-> matha12
}{}
\DeclareSymbolFont{matha}{U}{matha}{m}{n}
\DeclareFontFamily{U}{mathx}{\hyphenchar\font45}
\DeclareFontShape{U}{mathx}{m}{n}{
	<-6> mathx5 <6-7> mathx6 <7-8> mathx7
	<8-9> mathx8 <9-10> mathx9
	<10-12> mathx10 <12-> mathx12
}{}
\DeclareSymbolFont{mathx}{U}{mathx}{m}{n}
\DeclareMathDelimiter{\vvvert} {0}{matha}{"7E}{mathx}{"17}%

\newcommand\interVem{\mathrm{I}_\elementvem}
\newcommand\interEn{\widetilde{\mathrm{I}}_\elementvem}

\newcommand\mesh{\mathcal{T}}

\newcommand{\RE}{\mathbb{R}}

\renewcommand{\kappa}{k}

\newcommand{\Amatr}{\mathsf{A}}
\newcommand{\Bmatr}{\mathsf{B}}
\newcommand{\Cmatr}{\mathsf{C}}
\newcommand{\Fmatr}{\mathsf{F}}
\newcommand{\vmatr}{\mathsf{u}}
\newcommand{\pmatr}{\mathsf{p}}

\newcommand{\err}[2]{\mathtt{err^{#1}}(#2)}

\newcommand\mdeg{\sigma}
\newcommand{\Pclem}[1]{{#1_{I}^\mathrm{cl}}}

\newcommand{\PclemE}[1]{\rev{{\widetilde #1_{I}^\mathrm{cl}}}}
\newcommand{\hatv}{v}
\newcommand{\hatvcl}{\hatv^\mathrm{cl}_I}

\begin{document}
	
	\title[]{The MINI mixed virtual element\\for the Stokes equation}%
	
	\author{Silvia Bertoluzza}%
	\address{IMATI ``E. Magenes'', CNR, Pavia (Italy)}
	\email{silvia.bertoluzza@imati.cnr.it}%
	
	\author{Fabio Credali}%
	\address{CEMSE Division, King Abdullah Univeristy of Science and Technology, Thuwal, Saudi Arabia}
	\email{fabio.credali@kaust.edu.sa}%
	
	\author{Daniele Prada}%
	\address{IMATI ``E. Magenes'', CNR, Pavia (Italy)}
	\email{daniele.prada@imati.cnr.it}%
	
	\date{}
	
	\begin{abstract}
		We present and discuss a generalization of the popular MINI mixed finite element for the 2D Stokes equation by means of conforming virtual elements on polygonal meshes. 
		We prove optimal error estimates for both velocity and pressure. Theoretical results are confirmed by several numerical tests performed with different choices of polynomial accuracy and meshes. 
		
		\
		
		\textbf{MSC Class:} 65N12, 65N15, 65N30
	\end{abstract}
	
	\maketitle
	
	\section{Introduction}
	
	A wide range of physical phenomena in science and engineering are modeled by the Navier--Stokes equations, which are used to mathematically describe fluid flows. When the Reynolds number assumes very small values, the nonlinear convective term of the Navier--Stokes equation can be dropped out since the flow is dominated by diffusion. The Stokes problem is then considered in this limit situation~\cite{temam2001navier}.
	
	The necessity of simulating these equations originated a large interest in developing effective numerical methods and several mixed finite element methods have been designed for triangular and quadrilateral meshes~\cite{mixedFEM}. Among them, one of the most popular choices is the MINI element~\cite{arnold1984stable}, which was born as modification of the inf-sup unstable $\Poly{1}-\Poly{1}$ element. The idea behind the MINI element is to overcome the instability by enhancing the velocity space with the addition of local cubic bubble functions.
	
	In recent years, the scientific community increased its interest for the design of numerical methods dealing with generic polytopal meshes. We mention, for instance, composite finite elements~\cite{hackbusch1997composite,hackbusch1997composite2}, polygonal FEM~\cite{kuznetsov2003new,sukumar2004conforming,mu2012weak}, hybrid high-order methods~\cite{di2014arbitrary,di2015hybrid,cockburn2016bridging}, discontinuous Galerkin methods~\cite{cangiani2014hp,cockburn2016bridging} and, finally, virtual element methods~\cite{antonietti2022virtual,reviewVEM}.
	
	Since their introduction~\cite{beirao2013basic,ahmad2013equivalent}, virtual element methods have captured the attention of the scientific community thanks to their robustness and versatility, giving rise to a vast and dynamic literature investigating several aspects of the method: theoretical analysis~\cite{chen2018some,brenner2017some,beirao2017stability,cangiani2017posteriori,mascotto2018ill}, efficient implementation and solvers~\cite{beirao2014hitchhiker,bertoluzza2017bddc,antonietti2018multigrid,dassi2020parallel2,dassi2022robust,credali}, extensions~\cite{brezzi2014basic,beir2016basic,de2016nonconforming,gomez2024space}. Among the applications, we just mention elasticity~\cite{da2013virtual,gain2014virtual,dassi2021hybridization}, eigenvalue problems~\cite{mora2015virtual,gardini2018virtual,boffi2020gar,alzaben2024stabilization}, contact and deformation problems~\cite{wriggers2017low,cihan2022virtual}, fracture networks~\cite{berrone2022efficient,berrone2023virtual}, and Darcy flows~\cite{liu2019virtual,DASSI20221}.
	
	In this landscape, several virtual element formulations have been designed for the Stokes problem. We mention the stream formulation~\cite{antonietti2014stream}, divergence-free methods~\cite{da2017divergence,dassi2020bricks}, the theoretical study on the Stokes complex~\cite{beirao2019stokes}, \textit{a posteriori} estimates~\cite{wang2020posteriori}, conforming approximations~\cite{MANZINI-MAZZIA}, $p-$ and $hp-$ formulations~\cite{chernov2021p}. In all these discretization techniques, the virtual element machinery is employed just to construct the velocity space, while the pressure one is defined by means of pure polynomials. On the other hand, fully-virtual \textit{equal-order} formulations have been proposed in~\cite{guo2020new,li2024stabilized} as generalizations of the stabilized $\Poly{k}-\Poly{k}$ finite element. The stability of the equal-order elements is indeed ensured by introducing a pressure stabilization.
	
	In this landscape, we introduce our current work. We present a virtual element generalization of the MINI mixed finite element. More precisely, we consider equal-order virtual element spaces \rosso{of degree $k$} for velocity and pressure. The velocity space is then enriched by locally-defined ``virtual bubbles'' $H^1-$orthogonal \rosso{to polynomials of degree $k-2$.} A stabilization term is \rosso{then} employed for dealing with the non-polynomial contribution of the pressure. The proposed discretization technique is easy to implement and is proved to be stable and optimal since the norm of the polynomial contribution bounds the norm of the entire function.
	
	This paper is organized as follows. We recall the continuous formulation of the Stokes problem and the definition of the MINI mixed finite element in Sections~\ref{sec:cont_prob} and~\ref{sec:mini}, respectively. The main features of conforming virtual elements are recalled in Section~\ref{sec:conforming_vem}, while, in Section~\ref{sec:vem_bubbles}, we introduce the definition of virtual bubbles and discuss their properties. The construction of the virtual MINI element is described in Section~\ref{sec:mini_vem}, while its well-posedness and the error analysis are presented in Section~\ref{sec:well_posedness} and Section~\ref{sec:error_analysis}, respectively. Several numerical tests confirming our theoretical findings are discussed in Section~\ref{sec:num_tests}. In Section~\ref{sec:static_cond}, we describe how to obtain an equal-order virtual element method for the Stokes problem by static condensation of the MINI bubble functions. We finally draw some conclusions in Section~\ref{sec:conclusions}.
	
	\section{Notation}\label{sec:notation}
	Given an open bounded domain $D\subset\RE^2$, we denote by $C^0(D)$ the space of continuous functions over $D$. We denote by $\mathrm{L}^2(D)$ the space of square integrable functions endowed with the inner product $(u,v)_D=\int_D u\cdot v\,\dx$. The symbol $\mathrm{L}^2_0(D)$ refers to the subspace of $\mathrm{L}^2(D)$ containing null mean functions. Sobolev spaces with integrability exponent  equal to 2 \rosso{and differentiability order} $r$, are denoted by $\mathrm{H}^r(D)$, with the associated norm being $\|\cdot\|_{r,D}$. Moreover, $\mathrm{H}^1_0(D)$ denotes the subspace of $\mathrm{H}^1(D)$ of functions having zero trace on $\partial D$. Vector spaces will be denoted by bold letters. The letter $C$ will be used to denote a generic positive constant only depending on the shape of the domain \rosso{(see geometric assumptions (G1)--(G2) later on) and on the polynomial degree $k$, but independent of the mesh size. Moreover, given a mesh $\mesh_h$, we introduce the broken Sobolev space 
		\begin{equation*}
			\rv{\mathrm{H}^1(\mesh_h)} = \{v\in \mathrm{L}^2(\Omega):\,v_{|\elementvem}\in \mathrm{H}^1(\elementvem)\}
		\end{equation*}
		endowed with the broken  semi-norm and norm
		\begin{equation*}
			|v|_{1,h}^2 = \sum_{\elementvem\in\mesh_h} \|\grad v\|_{0,\elementvem}^2,
			\qquad
			\|v\|_{1,h}^2  = \|v\|_{0,\Omega}^2+|v|_{1,h}^2,
			\qquad \forall v\in \mathrm{H}^1(\mesh_h).
	\end{equation*}}
	
	$\Poly{k}(D)$ denotes the space of polynomials of degree $\le k$, with the convention that ${\Poly{-1}=\{0\}}$. Finally, we denote by  
	\begin{equation*}
		\rv{\monomials_{k}(D) = \{ m_\ell:\ell=1,\dots,\dim(\Poly{k}(D))\}}\qquad\text{with}\qquad \rosso{\monomials_{k}(D)\subset\monomials_{k+1}(D)}
	\end{equation*}
	a given basis for $\Poly{k}(D)$. Also in this case, we adopt the notation $\monomials_{-1}=\{0\}$.
	
	\section{The continuous problem}\label{sec:cont_prob}
	
	We consider a polygonal domain $\Omega\subset\RE^2$ and the incompressible Stokes problem: given an external force $\f$, we seek for the velocity $\u$ and the pressure $p$ governed by the following equations
	\begin{equation}\label{eq:stokes_strong}
		\begin{aligned}
			-\Delta\u + \grad p &= \f && \text{in }\Omega, \\
			\div\u &=0&& \text{in }\Omega, \\
			\u&=\mathbf{0}&& \text{on }\partial\Omega.
		\end{aligned}
	\end{equation}
	
	We introduce the functional spaces
	\begin{equation}\label{eq:cont_spaces}
		\V = \dueHunozeroO = [\HunozeroO]^2 \qquad\text{and}\qquad Q=\Ldo
	\end{equation}
	and the bilinear forms
	\begin{equation}
		\begin{aligned}
			&a:\V\times\V\longrightarrow\RE \qquad &&a(\u,\v) = \int_\Omega \Grad\u:\Grad\v\,\dx, \\
			&b:\V\times Q \longrightarrow\RE \qquad &&b(\v,q)=\int_\Omega q\,\div\v\,\dx,
		\end{aligned}
	\end{equation}
	where $\mathbf{A}:\mathbf{B}=\sum_{i,j=1}^{2}\mathbf{A}_{ij}\mathbf{B}_{ij}$ for all tensors $\mathbf{A}$ and $\mathbf{B}$.
	By standard manipulations of \eqref{eq:stokes_strong}, we can write the problem in variational form.
	
	\begin{pb}\label{pro:stokes}
		Find $\u\in\V$ and $p\in Q$ such that
		\begin{equation*}
			\begin{aligned}
				a(\u,\v) - b(\v,p) &= (\f,\v)_\Omega &&\forall\v\in\V,\\
				b(\u,q) &= 0 &&\forall q\in Q.
			\end{aligned}
		\end{equation*}
	\end{pb}
	
	This is a saddle point problem in which the pressure $p$ plays the role of Lagrange multiplier associated with the incompressibility condition \cite{mixedFEM}. Problem~\ref{pro:stokes} is stable and well-posed  since the bilinear form $a$ is continuous and coercive with constants $\mu,\gamma>0$, indeed
	\begin{equation}
		a(\u,\v)\le\mu\|\u\|_{1,\Omega}\|\v\|_{1,\Omega},\qquad a(\v,\v)\ge\gamma\|\v\|_{1,\Omega}^2\qquad\forall\u,\v\in\V,
	\end{equation}
	while $b$ is continuous 
	\begin{equation}
		b(\v,q) \le\zeta\|\v\|_{1,\Omega}\|q\|_{0,\Omega}\qquad\forall\v\in\V,\,\forall q\in Q
	\end{equation}
	and satisfies the following inf-sup condition for a positive constant~$\beta$
	\begin{equation}\label{continfsup}
		\inf_{q\in Q,q\neq0} \sup_{\v\in\V,\v\neq\mathbf{0}} \frac{b(\v,q)}{\|\v\|_{1,\Omega}\|q\|_{0,\Omega}} \ge \beta.
	\end{equation}
	
	Notice that, by introducing the product space $\VV=\V\times Q$ and the bilinear form
	\begin{equation}\label{eq:double_bil}
		B[(\u,p),(\v,q)] = a(\u,\v) - b(\v,p) +b(\u,q) \qquad \forall(\u,p),(\v,q)\in\VV,
	\end{equation}
	Problem~\ref{pro:stokes} can be written in the following equivalent form.
	\begin{pb}\label{pro:stokes2}
		Find $(\u,p)\in\VV$ such that
		\begin{equation*}
			B[(\u,p),(\v,q)] = (\f,\v)_\Omega \qquad \forall(\v,q)\in\VV.
		\end{equation*}
	\end{pb}
	\rosso{The problem is well-posed thanks to continuity and compatibility properties of $B$, i.e.
		\begin{equation}
			B[(\u,p),(\v,q)] \le C \, \|(\u,p)\|_{\VV} \|(\v,q)\|_{\VV},
			\qquad
			\sup_{(\v,q)\in\VV} \frac{B[(\u,p),(\v,q)]}{\|(\v,q)\|_{\VV}} \ge \omega\,\|(\u,p)\|_{\VV},
		\end{equation}
		with $\|(\cdot,\cdot)\|_{\VV}$ being the product norm defined as
		\begin{equation}
			\|(\v,q)\|_{\VV}^2 = \|\v\|_{1,\Omega}^2 + \|p\|_{0,\Omega}^2.
		\end{equation}
	}
	
	Finally, the following estimate holds true \cite{mixedFEM}.
	
	\begin{thm}
		Given $\f\in\bH^{-1}(\Omega)$, there exists a unique pair $(\u,p)\in\VV$ solving Problem~\ref{pro:stokes}/\ref{pro:stokes2} and satisfying
		\[
		\|\u\|_{1,\Omega} + \|p\|_{0,\Omega} \le C\,\|\f\|_{-1,\Omega}.
		\]
	\end{thm}
	
	\section{The MINI mixed finite element: recall}\label{sec:mini}
	
	The MINI finite element~\cite{arnold1984stable} was designed specifically for the Stokes equation as a modified version of the unstable pair of continuous piecewise linear velocity and pressure $\PolydueM{1}-\Poly{1}$, which do not satisfy the inf-sup conditions~\cite[Chap. 8, Sec. 3.1]{mixedFEM}. The idea behind MINI is to enrich the velocity space $\PolydueM{1}$ by adding a cubic bubble in each element of the domain discretization.
	
\begin{figure}
		\centering
		
		\includegraphics[width=0.7\linewidth]{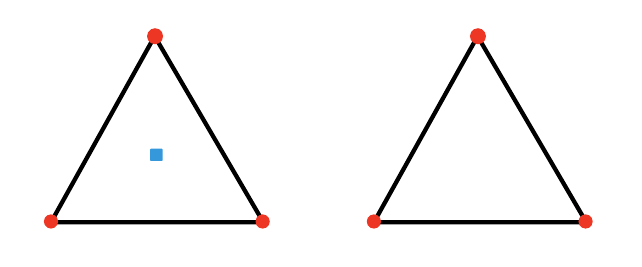}

		\caption{Graphical representation of the MINI mixed finite element. Degrees of freedom of the linear part are represented by red circles, whereas the bubble degree of freedom is indicated by a blue square.}
		\label{fig:mini_fem}
	\end{figure}

	More precisely, let us consider a conforming finite element triangulation $\mathfrak{T}_h$ of $\Omega$. The MINI velocity space is defined as
	\begin{equation}\label{eq:mini1}
		\begin{aligned}
			\V_h^{\textrm{mini}} = \{ \v_h\in\mathbf{C}^0(\Omega):\, &\v_{h|T} = \v_{h,1} + \v_{h,b}\\
			&\v_{h,1}\in\PolydueM{1}(T),\,\v_{h,b}\in\PolydueM{3}(T)\cap\mathrm{\mathbf{H}}_0^1(T)\quad\forall T\in\mathfrak{T}_h\},
		\end{aligned}
	\end{equation}
	while the pressure space is made up of continuous piecewise linear polynomials
	\begin{equation}\label{eq:mini2}
		Q_h^{\textrm{mini}} = \{ q_h\in C^0(\Omega):\,q_{h|T}\in\Poly{1}(T)\quad\forall T\in\mathfrak{T}_h\}.
	\end{equation}
	In Figure~\ref{fig:mini_fem}, we depict the degrees of freedom of the MINI element in a triangle $T\in\mathfrak{T}_h$. Given $\v_h\in\V_h^{\textrm{mini}}$, the linear contribution $\v_{h,1}$ is identified by its values at the vertices of $T$, whereas the bubble contribution $\v_{h,b}$ is identified by its value at the barycenter. The pressure $q_h\in Q_h^{\textrm{mini}}$ is described by the values at the vertices.
	
	From the standard theory \cite[Chap. 8, Sec. 4.2]{mixedFEM}, it is well-known that the MINI finite element satisfies the inf-sup conditions and yields linear convergence of the error. Indeed, if $(\u_h,p_h)$ is the discrete MINI solution to the Stokes equation, the following error estimate holds true provided that the exact solution $(\u,p)\in \bH^2(\Omega)\times \mathrm{H}^1(\Omega)$
	\begin{equation}
		\| \u - \u_h \|_{1,\Omega} + \| p - p_h \|_{0,\Omega} \le C\,h\,(\|\u\|_{2,\Omega}+\|p\|_{1,\Omega}).
	\end{equation}
	
	We remark that the bubble enrichment $\u_{h,b}$ does not really contribute to the approximation of the velocity, but it is just a mathematical tool for stabilizing the discrete formulation~\cite{verfurth1989posteriori,bank1990posteriori,bank1991posteriori,kim2000modified,russo1995posteriori}.
	
	\section{Conforming virtual elements: definition and properties}\label{sec:conforming_vem}
	
	In this section we recall the definition of plain~\cite{beirao2013basic} and enhanced~\cite{ahmad2013equivalent} virtual element spaces, also discussing some of their approximation properties~\cite{chen2018some}.
	
	\subsection{Domain discretization} \label{sec:geometry}
	
	Let $\{\mesh_h\}_h$ be a family of decompositions of $\Omega$ into polygonal elements $\elementvem$. We denote by $h_\elementvem$ the diameter of $\elementvem$ and by $e\in\partial\elementvem$ its edges. As usual in virtual elements theory, we assume that the following geometric requirements are satisfied by every element $\elementvem$ of $\mesh_h$ \cite{beirao2013basic}:
	\begin{enumerate}[(G1)]
		\item $\elementvem$ is star-shaped with respect to a ball of radius $\ge \rho_1 h_\elementvem$,
		\item the distance between every pair of vertices is $\ge \rho_2 h_\elementvem$. 
	\end{enumerate}
	Several variants of (G1) and (G2) have been proposed in literature, as well as additional requirements. An overview of the different geometrical assumptions considered in the VEM literature can be found in~\cite{sorgente2021vem}.
	
	\subsection{Projectors}\label{sec:projectors}
	
	When constructing virtual element discretizations, a key role is played by projectors onto polynomial spaces.
	
	Given $v\in\HunoK$, we define the elliptic projector $\Pinabla{k}:\HunoK\rightarrow\Poly{k}(\elementvem)$ as the solution of the following problem
	\begin{equation}\label{eq:def_pinabla_sca}
		\begin{aligned}			
			\int_{\elementvem} \nabla \Pinabla{k} v\cdot \nabla q \,\dx &= \int_{\elementvem} \nabla v \cdot \nabla q\,\dx \quad \forall q\in\Poly{k}(\elementvem),\\
			\int_{\bK} \Pinabla{k} v \,\ds &= \int_{\bK} v\,\ds.
		\end{aligned}
	\end{equation}
	
	We also introduce the orthogonal projection onto $\Poly{k}(\elementvem)$ as the operator $\Pizero{k}:\mathrm{L}^2(\elementvem)\rightarrow\Poly{k}(\elementvem)$ defined by
	\begin{equation}
		\int_\elementvem \Pizero{k}v\,q\,\dx = \int_\elementvem v\,q\,\dx \quad \forall q\in\Poly{k}(\elementvem).
	\end{equation}
	The following projection estimates hold true~\cite{beirao2016virtual}.
	\begin{lem}\label{lem:proj_err}
		For all $v\in \mathrm{H}^s(\Omega)$ and $\elementvem\in\mesh_h$, there holds
		\begin{equation*}
			\begin{aligned}
				&\| v-\Pizero{k}v\|_{r,\elementvem} \le Ch_\elementvem^{s-r}|v|_{s,\elementvem}&&\quad r,s\in\mathbb{N},\,0\le r\le s\le k+1,\\
				&\| v-\Pinabla{k}v\|_{r,\elementvem} \le Ch_\elementvem^{s-r}|v|_{s,\elementvem}&&\quad r,s\in\mathbb{N},\,0\le r\le s\le k+1,\,s\ge1.\\
			\end{aligned}
		\end{equation*}
	\end{lem}
	
	When dealing with vector fields, we will employ the same notation for componentwise projection:
	\begin{equation}
		\PinablaG{k}\v=(\Pinabla{k}v_x,\Pinabla{k}v_y),\qquad\PizeroG{k}\ww=(\Pizero{k}w_x,\Pizero{k}w_y),
	\end{equation}
	for $\v=(v_x,v_y)\in\mathrm{\mathbf{H}}^1(\elementvem)$ and $\ww=(w_x,w_y)\in\mathrm{\mathbf{L}}^2(\elementvem)$.
	
	\subsection{Discrete spaces}\label{sec:vem}
	
	The \textit{plain} local virtual element space \cite{beirao2013basic} of degree $k$ is given by
	\begin{equation}\label{eq:classsical_vem}
		\vem{k}(\elementvem) = \{ v\in\HunoK:\, v_{|\edge}\in\Poly{k}(\edge)\,\forall\edge\in\partial\elementvem,\,\Delta v\in\Poly{k-2}(\elementvem)\},
	\end{equation}
	that is each function $v\in\vem{k}(\elementvem)$ is, on each edge of $\elementvem$, a polynomial of degree $k$ and its Laplacian is a polynomial of degree $k-2$ in the interior. To define the \textit{enhanced} VEM space $\enhanced{k}(\elementvem)$ (see \cite{ahmad2013equivalent}), we let
	\begin{equation}
		\enlarged{k}(\elementvem) = \{ v\in\HunoK:\, v\in\Poly{k}(\edge)\,\forall\edge\in\partial\elementvem,\,\Delta v\in\Poly{k}(\elementvem)\},
	\end{equation}
	and we set
	\begin{equation}
		\enhanced{k}(\elementvem) = \{ v\in\enlarged{k}(\elementvem):\, (v-\Pinabla{k}v,q_k)_\elementvem=0\quad\forall q_k\in\monomials_{k}(\elementvem)\setminus\monomials_{k-2}(\elementvem)\}.
	\end{equation}
	
	We have $\dim(\vem{k}(\elementvem))=\dim(\enhanced{k}(\elementvem))$. We can describe both spaces by the same set of unisolvent degrees of freedom. We consider:
	\begin{itemize}
		\item the values of $v$ at the vertices of $\elementvem$;
		\item for $k\ge2$, the values of $v$ in $k-1$ points on each edge $\edge\in\partial\elementvem$;
		\item for $k\ge2$, the internal moments
		\begin{equation}\label{eq:moments_def}
			\frac1{|\elementvem|} \int_\elementvem v\,m\,\dx\qquad\forall m\in\monomials_{k-2}(\elementvem).
		\end{equation}
	\end{itemize}
	
	The global spaces are then obtained by glueing all the local spaces by continuity
	\begin{equation*}
		\begin{aligned}
			\vem{k} &= \{v\in\HunozeroO:\,v_{|\elementvem}\in\vem{k}(\elementvem)\quad\forall\elementvem\in\mesh_h\},\\
			\enhanced{k} &= \{v\in\HunozeroO:\,v_{|\elementvem}\in\enhanced{k}(\elementvem)\quad\forall\elementvem\in\mesh_h\}.\\
		\end{aligned}
	\end{equation*}
	
	\begin{rem}
		Notice that, for $v\in\vem{k}(\elementvem)$, the knowledge of the degrees of freedom allows only the computation of the elliptic projection~$\Pinabla{k}v$. On the other hand, if $v\in\enhanced{k}(\elementvem)$, we can also compute the orthogonal projection~$\Pizero{k}v$ since all the moments up to order~$k$ are known.
	\end{rem}
	
	Given a function $w\in\HunoK\cap C^0(\elementvem)$, we introduce the interpolant $\interVem w\in\vem{k}(\elementvem)$ defined by
	\begin{equation}
		\interVem w = w_I \quad\text{on }\partial\elementvem,\qquad \int_\elementvem \interVem w\,q\,\dx = \int_\elementvem w\,q\,\dx\quad\forall q\in\Poly{k-2}(\elementvem),
	\end{equation}
	\rosso{where $w_I$ denotes the standard nodal interpolant of $w$ on $\partial\elementvem$.} This operator satisfies the following interpolation error estimate~\cite{chen2018some}.
	\begin{prop}
		For \rv{$v\in \mathrm{H}^s(\elementvem)$, $2 \leq s \leq k+1$}, the following optimal order error estimate holds true
		\begin{equation}
			\| v-\interVem v \|_{0,\elementvem} + h_\elementvem | v-\interVem v|_{1,\elementvem} \leq C\,\rv{h_\elementvem^{s}} \|v\|_{\rv{s,\elementvem}.}
		\end{equation}
	\end{prop}
	
	In the same way, we can define the analogous interpolation operator $\interEn w$ for the enhanced space $\enhanced{k}(\elementvem)$ as
	\begin{equation}
		\interEn w = w_I \quad\text{on }\partial\elementvem,\qquad \int_\elementvem \interEn w\,q\,\dx = \int_\elementvem w\,q\,\dx\quad\forall q\in\Poly{k-2}(\elementvem),
	\end{equation}
	also satisfying the analogous estimate~\cite{chen2018some}.
	\begin{prop}\label{prop:interp_enh}
		For \rv{$v\in \mathrm{H}^{s}(\elementvem)$, $2\leq s \leq k+1$,} the following optimal order error estimate holds true
		\begin{equation}
			\| v-\interEn v \|_{0,\elementvem} + h_\elementvem | v-\interEn v|_{1,\elementvem} \leq C\,\rv{h_\elementvem^{s}} \|v\|_{\rv{s,\elementvem}.}
		\end{equation}
	\end{prop}
	Finally, it is also possible to define a Cl\'ement type interpolant as stated by the following proposition.
	\begin{prop}\label{clement} For $v \in H^1(\Omega)$ there exists an element $\PclemE v \in \enhanced{k}$ such that for all $\elementvem \in \mesh_h$ we have
		\begin{gather}
			\label{clementapprox}
			\| v - \PclemE v \|_{0,\elementvem} + h_\elementvem | v-\PclemE v|_{1,\elementvem} \leq \, C\, h_\elementvem | v |_{1,\rev{\widetilde{\elementvem}}},
		\end{gather}
		\rev{with $\widetilde{\elementvem}$ denoting the vertex patch of $\elementvem$.}
		Moreover it is possible to choose $\PclemE v$ in such a way that
		\[\label{clementmoments} \int_\elementvem (v - \PclemE v)\,q\,\dx = 0 \qquad\forall q\in\Poly{k-2}(\elementvem).\]
	\end{prop}
	\begin{proof}
		
		\rev{We first prove the proposition for a Cl\'ement type interpolant $\hatvcl$ defined in the plain virtual element space $\vem{k}$ and then we extend the result to $\PclemE v\in\enhanced{k}$.
			
			By \cite[Proof of Prop. 4.2]{mora2015virtual}, there exists an element $\hatv_I\in\vem{k}(\elementvem)$ such that 
			\begin{equation}\label{eq:clem_2}
				\| v - \hatv_I \|_{0,\elementvem} + h_\elementvem | v-\hatv_I|_{1,\elementvem} \leq \, C\, h_\elementvem | v |_{1,\rev{\widetilde{\elementvem}}}.
			\end{equation}
			Let $\hatvcl\in\vem{k}(\elementvem)$ be defined in such a way that the vertices and edge degrees of freedom coincide with those of $\hatv_I$, while the interior degrees of freedom are defined by
			\begin{equation}\label{eq:moments_def}
				\frac1{|\elementvem|} \int_\elementvem \hatvcl\,m\,\dx = 
				\frac1{|\elementvem|} \int_\elementvem v\,m\,\dx
				\qquad\forall m\in\monomials_{k-2}(\elementvem).
			\end{equation}
			We observe that, for $k=1$, $\hatvcl=\hatv_I$, which satisfies~\eqref{eq:clem_2}.
			We show that $\hatvcl$ satisfies \eqref{eq:clem_2} also for $k\ge2$. As $\hatv_I -\hatvcl = 0$ of $\partial\elementvem$, we have that
			\[
			\begin{aligned}
				| \hatv_I - \hatvcl |^2_{1,\elementvem} &= 
				- \int_\elementvem \Delta(\hatv_I - \hatvcl) (\hatv_I - \hatvcl)\,\dx = - \int_\elementvem \Delta(\hatv_I - \hatvcl) (\hatv_I - v) \,\dx\\
				& \leq \| \Delta(\hatv_I - \hatvcl) \|_{0,\elementvem} \| \hatv_I - v \|_{0,K} \leq C h_\elementvem^{-1} | \hatv_I - \hatvcl |_{1,\elementvem} h_\elementvem | v |_{1,\rev{\widetilde{\elementvem}}}.
			\end{aligned}
			\]
			By dividing both sides by $| \hatv_I - \hatvcl |_{1,\elementvem}$ we obtain the bound on the $H^1$ seminorm. The bound on the $L^2$ norm follows by Poincar\'e  inequality as $v - \hatvcl$ is average free.
			
			We now define $\PclemE v\in\enhanced{k}(\elementvem)$ as the element that shares the same degrees of freedom as $\hatvcl$, although defined in a different space. By applying the same reasoning as in~\cite[Proof of Theorem 5.4]{chen2018some}, we find that $\PclemE v$ satisfies the desired estimate.}
	\end{proof}
	
	\section{The space of virtual bubbles}\label{sec:vem_bubbles} 
	
	In this section, we define the space of virtual bubbles and discuss its main properties, which will be exploited for defining the MINI mixed virtual element.
	
	Given $\kappa\in\mathbb{Z}$, $\kappa\ge2$ and given a generic polygon $\elementvem$ satisfying (G1)--(G2), we define the local space of virtual bubbles as
	\begin{equation}
		\bubble{\kappa} = \{ \bolla \in\HunozeroK:\,\Delta \bolla\in\Poly{\kappa-2}(\elementvem)\}.
	\end{equation}
	Notice that $\bubble{1}=\{0\}$.
	The choice of degrees of freedom for $\bubble{\kappa}$ is a direct consequence of what described in \cite{beirao2013basic}: we do not need boundary degrees of freedom, hence we just consider the internal moments
	\begin{equation}
		\frac{1}{|\elementvem|}\int_\elementvem \bolla\,m\,\dx \qquad \forall m\in\monomials_{\kappa-2}(\elementvem).
	\end{equation}
	
	The first property we present relates the norm of $\bolla$ with the norm of its elliptic projection~$\Pinabla{\kappa}\bolla$.
	\begin{prop}\label{prop:norm_bub}
		Given $\bolla\in\bubble{\kappa}$, there exists a constant $\gamma_\sharp$ such that the following inequality holds true
		\begin{equation}
			|\bolla|_{1,\elementvem}^2 \le \gamma_\sharp |\Pinabla{\kappa}\bolla|_{1,\elementvem}^2.
		\end{equation}
	\end{prop}

	\begin{proof}
		\rv{	We start by recalling that, for all $v \in  \vem{k}(\elementvem) \cap \ker \Pinabla{k}$ it holds that 
			\begin{equation}\label{red2b}
				| v |_{1,K} \leq C | v |_{1/2,K}.
			\end{equation}
			Indeed, let $\widetilde p \in \Poly{k}(\elementvem)$ be such that $\Delta \widetilde p = \Delta v$ and $| \widetilde p |_{1,\elementvem} \leq C | v |_{1,\elementvem}$ (such a polynomial exists, thanks to
			\cite[Lemma 3.5]{beirao2017stability} and \cite[Lemma 10]{cangiani2017posteriori}). \rev{Integrating by parts and using} $\Pinabla{k} v = 0$ can write
			\begin{equation*}
				\begin{aligned}
					| v |_{1,\elementvem}^2 &=  \int_\rev{\partial\elementvem} \nabla (v - \widetilde p)\cdot\nu v \,\dx\leq \rev{\| \nabla(v - \widetilde p)\cdot\nu \|_{-1/2,\partial\elementvem} | v |_{1/2,\partial\elementvem}} \\&\leq C | v - \widetilde p |_{1,\elementvem}  | v |_{1/2,\partial\elementvem} \leq C | v |_{1,\elementvem}  | v |_{1/2,\partial\elementvem},
				\end{aligned}  
			\end{equation*}
			where we used that $v - \widetilde p$ is harmonic.}
		By applying the triangle inequality and \eqref{red2b}, we can write
		\begin{equation*}
			\begin{aligned}
				|\bolla|_{1,\elementvem}^2 &= \int_\elementvem |\grad\bolla|^2\,\dx
				\le\int_\elementvem |\grad\Pinabla{\kappa}\bolla|^2\,\dx + \int_\elementvem | \bolla-\Pinabla{\kappa}\bolla|^2\,\dx\\
				&\le  \int_\elementvem |\grad\Pinabla{\kappa}\bolla|^2\,\dx + C\, |\bolla-\Pinabla{\kappa}\bolla|_{1/2,\partial\elementvem}^2,
			\end{aligned}
		\end{equation*}
		hence, by taking into account that $\bolla=0$ on $\partial\elementvem$, we obtain
		\begin{equation}
			|\bolla|_{1,\elementvem}^2
			\le |\Pinabla{\kappa}\bolla|_{1,\elementvem}^2 +  C\, |\Pinabla{\kappa}\bolla|_{1/2,\partial\elementvem}^2 \le C\, |\Pinabla{\kappa}\bolla|_{1,\elementvem}^2.
		\end{equation}
	\end{proof}

	The second property of bubbles is their elliptic orthogonality to the space $\harmonic{\kappa}(\elementvem)$ of harmonic polynomials of degree $\kappa$, i.e.
	\begin{equation}
		\harmonic{\kappa}(\elementvem) = \{ q\in\Poly{\kappa}(\elementvem):\,\Delta q=0\}.
	\end{equation}
	Indeed, given $\bolla\in\bubble{\kappa}$, we have
	\begin{equation}\label{eq:horto}
		\int_\elementvem \grad \bolla\cdot\grad q\,\dx = -\int_\elementvem \bolla\,\Delta q\,\dx+\int_{\partial\elementvem}(\grad q\cdot\n_\elementvem)\bolla\,\ds = 0\quad\forall q\in\harmonic{\kappa}(\elementvem).
	\end{equation}
	Thus, the bubble space $\bubble{\kappa}$ is orthogonal to $\harmonic{\kappa}(\elementvem)$ with respect to the $H^1$ semi-scalar product.
	By combining~\eqref{eq:horto} with the definition of $\Pinabla{\kappa}\bolla$,
	\begin{equation*}
		\int_\elementvem \grad \bolla\cdot\grad q\,\dx = \int_\elementvem \grad\Pinabla{\kappa}\bolla\cdot\grad q\,\dx\quad\forall q\in\harmonic{\kappa}(\elementvem),
	\end{equation*}
	we easily prove the following orthogonality results for $\Pinabla{\kappa}\bolla$.
	\begin{prop}
		Given $\bolla\in\bubble{\kappa}$, it holds
		$$
		\int_\elementvem \grad\Pinabla{\kappa}\bolla\cdot\grad q\,\dx=0\quad\forall q\in\harmonic{\kappa}(\elementvem).
		$$
	\end{prop}
	
	\begin{cor}\label{cor:ort}
		The projected bubble space $\Pinabla{\kappa}(\bubble{\kappa})$ satisfies
		\begin{equation}
			\Pinabla{\kappa}(\bubble{\kappa})
			\subseteq  \Poly{\kappa}(\elementvem)\cap\harmonic{\kappa}^\perp(\elementvem),
		\end{equation}
		where $\harmonic{\kappa}^\perp(\elementvem) \subset \Poly{\kappa}(\elementvem)$ is the $H^1$-orthogonal subspace to $\harmonic{\kappa}(\elementvem)$.
	\end{cor}
	
	In the following we will need the space $\rvf{\bubbleSp}$ defined as
	\begin{equation}
		\rvf{\bubbleSp} = \{\bolla\in\bubble{\rvf{k+1}}:\,\bolla\perp\Poly{k-2}(\elementvem)\},
	\end{equation}
	\rv{where the subscript denotes the bubble degree, while the superscript indicates orthogonality with respect to polynomials of degree $k-2$.}
	\rv{By construction, $\dim(\bubbleSp)=k$ and} $\bolla\in\rvf{\bubbleSp}$ is uniquely identified by the moments
	\begin{equation*}
		\frac{1}{|\elementvem|}\int_\elementvem \bolla\,m\,\dx \qquad \forall m\in\monomials_{\rvf{\kappa-1}}(\elementvem)\setminus\monomials_{\kappa-2}(\elementvem),
	\end{equation*}
	while
	\begin{equation*}
		\frac{1}{|\elementvem|}\int_\elementvem \bolla\,m\,\dx =0 \qquad \forall m\in\monomials_{\kappa-2}(\elementvem).
	\end{equation*}
	It is immediate to see that, if $k=1$, then $\rvf{\mathcal{B}_{2}^{-1}(\elementvem)=\bubble{2}}$.
	
	\section{The new \rosso{$\PolydueM{k+b}-\Poly{k}$} Stokes \rosso{virtual} element}\label{sec:mini_vem}
	
	\rosso{In this section we present the new Stokes element. While the main goal of our work is to extend to polytopal meshes the MINI finite element $\PolydueM{1+b}-\Poly{1}$ recalled in~{\eqref{eq:mini1}--\eqref{eq:mini2}}, as the construction and analysis easily generalize to $\PolydueM{k+b}-\Poly{k}$, we directly present the latter. The MINI--VEM results from choosing $k=1$. By abuse of notation, from now on we will call MINI--VEM the new element of degree $k$.}
	
	Given $k\in\mathbb{Z}$, $k\ge1$ and a generic polygon $\elementvem$ satisfying the assumptions (G1)--(G2), the local velocity space $\V_k(\elementvem)$ is defined as
	\begin{equation}\label{eq:vel_mini}
		\V_k(\elementvem) = \enhancedunoS\oplus\bubbleS,
	\end{equation}
	where 
	\begin{equation}
		\enhancedunoS = [\enhanced{k}(\elementvem)]^2 \qquad\text{and}\qquad \bubbleS = [\bubbleSp]^2.
	\end{equation}
	\rosso{We  observe that~\eqref{eq:vel_mini} is a direct sum, that is $\enhancedunoS\cap\bubbleS=\{\mathbf{0}\}$.} In our discussion, we are going to denote trial and test functions of $\V_k(\elementvem)$ as
	\begin{equation}\label{eq:decomposition}
		\u_h = \utilde + \bolla_h, \quad \v_h = \vtilde + \bollatest_h
	\end{equation}
	with $\utilde,\vtilde\in\enhancedunoS$ and $\bolla_h,\bollatest_h\in\bubbleS$. \rosso{We prove that such splitting is stable.}
	\begin{prop}\label{prop:splitting}
		The splitting $\v_h=\vtilde+\bollatest_h$ is stable, that is, there exists a positive constant $\eta$, independent of $h$, such that
		\begin{equation}
			\eta (|\vtilde|_{1,\elementvem}+|\bollatest_h|_{1,\elementvem}) \le |\v_h|_{1,\elementvem} \le |\vtilde|_{1,\elementvem}+|\bollatest_h|_{1,\elementvem}.
		\end{equation}
	\end{prop}
	\begin{proof}
		
		The upper bound is easily obtained by the triangle inequality. Regarding the lower bound, we first notice that
		\begin{equation}
			|\vtilde|_{1,\elementvem}+|\bollatest_h|_{1,\elementvem}
			\le |\vtilde+\bollatest_h|_{1,\elementvem}+2|\bollatest_h|_{1,\elementvem}
			\le C\,(|\v_h|_{1,\elementvem}+|\bollatest_h|_{1,\elementvem}),
		\end{equation}
		therefore it remains to prove that
		\begin{equation}\label{eq:claim}
			|\bollatest_h|_{1,\elementvem}\le C\,|\v_h|_{1,\elementvem}.
		\end{equation}
		Proposition~\ref{prop:norm_bub}, together with the triangle inequality and the stability of the projection operator, yields
		\begin{equation}\label{eq:step1}
			\begin{aligned}
				|\bollatest_h|_{1,\elementvem}^2
				&\le \gamma_\sharp\,|\Pinablakpu\bollatest_h|_{1,\elementvem}^2
				\le \gamma_\sharp\,(|\Pinablakpu(\bollatest_h+\vtilde)|_{1,\elementvem}^2+|\Pinablakpu\vtilde|_{1,\elementvem}^2)\\
				&\le \gamma_\sharp\,(|\bollatest_h+\vtilde|_{1,\elementvem}^2+|\Pinablakpu\vtilde|_{1,\elementvem}^2).
			\end{aligned}
		\end{equation}
		In order to bound the projector norm, we exploit that $\Pinabla{k}\bollatest_h=\mathbf{0}$ and the projection stability, thus
		\begin{equation}\label{eq:step2}
			\begin{aligned}
				|\Pinablakpu\vtilde|_{1,\elementvem}^2
				&\le C\,(|\Pinablakpu\vtilde-\Pinabla{k}\vtilde|_{1,\elementvem}^2+|\Pinabla{k}(\vtilde+\bollatest_h)|_{1,\elementvem}^2)\\
				&\le C\,(|\Pinablakpu(\vtilde-\Pinabla{k}\vtilde)|_{1,\elementvem}^2+|\Pinabla{k}(\vtilde+\bollatest_h)|_{1,\elementvem}^2)\\
				&\le C\,(|\vtilde-\Pinabla{k}\vtilde|_{1,\elementvem}^2+|\vtilde+\bollatest_h|_{1,\elementvem}^2).
			\end{aligned}
		\end{equation}
		By applying \eqref{red2b} and a scaling of boundary norms, we can write
		\begin{equation}
			|\vtilde-\Pinabla{k}\vtilde|_{1,\elementvem}^2 \le C\, |\vtilde-\Pinabla{k}\vtilde|_{1/2,\partial\elementvem}^2 \le C\rv{h_\elementvem^{-1}}\,\|\vtilde-\Pinabla{k}\vtilde\|_{0,\partial\elementvem}^2.
		\end{equation}
		Using that $\bollatest_h=\mathbf{0}$ on $\partial\elementvem$, again $\Pinabla{k}\bollatest_h=\mathbf{0}$, and another scaling argument, we obtain
		\begin{equation*}
			\begin{aligned}
				|\vtilde-\Pinabla{k}\vtilde|_{1,\elementvem}^2
				&\le C\rv{h_\elementvem^{-1}}\,\|\vtilde+\bollatest_h-\Pinabla{k}(\vtilde+\bollatest_h)\|_{0,\partial\elementvem}^2\\
				&\le C\rv{h_\elementvem^{-2}}\,\|\vtilde+\bollatest_h-\Pinabla{k}(\vtilde+\bollatest_h)\|_{0,\elementvem}^2+|\vtilde+\bollatest_h-\Pinabla{k}(\vtilde+\bollatest_h)|_{1,\elementvem}^2,
			\end{aligned}
		\end{equation*}
		so that
		\begin{equation}\label{eq:step3}
			|\vtilde-\Pinabla{k}\vtilde|_{1,\elementvem}^2\le C\,|\vtilde+\bollatest_h|_{1,\elementvem}^2.
		\end{equation}
		Finally, combining~\eqref{eq:step1}, \eqref{eq:step2}, \eqref{eq:step3} we obtain~\eqref{eq:claim} and the result is proved.
	\end{proof}
	
	The choice of degrees of freedom for this velocity space is a direct consequence of its construction. First, we notice that
	\[
	\dim(\V_k(\elementvem)) = \dim(\enhancedunoS)+\dim(\bubbleS).
	\]
	Therefore, the degrees of freedom for $\V_k(\elementvem)$ are inherited from $\enhancedunoS$ and $\bubbleS$ respectively. In particular, for $\v_h=\vtilde + \bollatest_h \in\V_k(\elementvem)$, we choose
	\begin{itemize}
		\item \textbf{dofs(v)\#1,} the values of $\vtilde$ at the vertices of $\elementvem$;
		\item \textbf{dofs(v)\#2,} for $k\ge2$, the values of $\vtilde$ at $k-1$ internal points on each edge $\edge\in\partial\elementvem$;
		\item \textbf{dofs(v)\#3,} for $k\ge2$, the internal moments \rosso{of $\vtilde$}
		\begin{equation*}
			\frac1{|\elementvem|}\int_\elementvem \vtilde:\m\,\dx \qquad \forall\m\in\monomialsdue_{k-2}(\elementvem).
		\end{equation*}
		\item \textbf{dofs(v)\#4,} the higher order moments \rosso{of the bubble contribution $\bollatest_h$}
		\begin{equation*}
			\frac1{|\elementvem|}\int_\elementvem \bollatest_h:\m\,\dx \qquad \forall\m\in\rvf{\monomialsdue_{k-1}(\elementvem)}\setminus\monomialsdue_{k-2}(\elementvem).
		\end{equation*}
	\end{itemize}
	
	\rosso{We emphasize that the moments \textbf{dofs(v)\#3,} identify the contribution $\vtilde$ since the bubble part is orthogonal to polynomials up to degree $k-2$.} As $\V_k(\elementvem)\subset[\rvf{\enhanced{k+1}(\elementvem)}]^2$, the following proposition holds.
	
	\begin{prop}
		The degrees of freedom \textbf{dofs(v)} are unisolvent for the velocity space~$\V_k(\elementvem)$.
	\end{prop}
	
	For the pressure, we consider the following space
	\begin{equation}
		Q_k(\elementvem) = \enhanced{k}(\elementvem).
	\end{equation}
	The degrees of freedom for $q_h\in Q_k(\elementvem)$ are (see Section~\ref{sec:vem})
	\begin{itemize}
		\item the values of $q_h$ at the vertices of $\elementvem$;
		\item the values of $q_h$ in $k-1$ points on each edge $\edge\in\partial\elementvem$;
		\item the internal moments
		\begin{equation}
			\frac1{|\elementvem|} \int_\elementvem q_h\,m\,\dx\qquad\forall m\in\monomials_{k-2}(\elementvem).
		\end{equation}
	\end{itemize}
	
\begin{figure}
		\centering
		\includegraphics[width=0.7\linewidth]{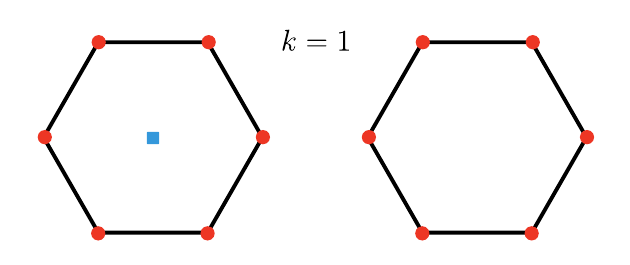}\\
		\includegraphics[width=0.7\linewidth]{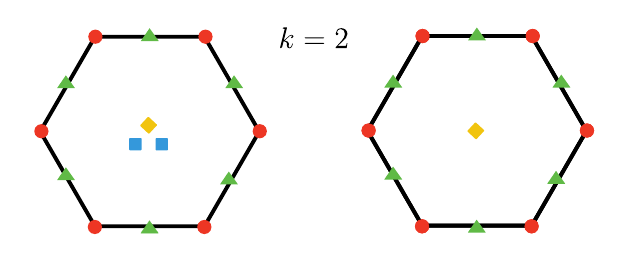}
		
		\caption{Graphical representation of the MINI mixed virtual element for $k=1,2$. The degrees of freedom of the degree $k$ contribution are denoted by red circles (vertices), green triangles (internal points on edges) and yellow diamonds (internal moments). The degrees of freedom of the additional bubbles are denoted by blue squares.}
		\label{fig:mini_vem}
	\end{figure}

	A sketch of the local degrees of freedom of the MINI--VEM of order $k=1,2$ is depicted in Figure~\ref{fig:mini_vem}.
	
	After having introduced the local spaces and discussed their main properties, we now define the global discrete spaces as
	\begin{equation}
		\begin{aligned}
			&\V_k = \{ \v_h\in\V:\, \v_h\in\V_k(\elementvem)\quad\forall\elementvem\in\mesh_h\},\\&Q_k = \{ q_h\in Q:\,q_h\in Q_k(\elementvem)\quad\forall\elementvem\in\mesh_h\}.
		\end{aligned}
	\end{equation}
	
	\begin{rem}
		The bubble enrichment of the velocity space naturally stabilizes the polynomial part of the pressure. In the next sections, we will see that a VEM type stabilization term is required for dealing with the nonpolynomial pressure contribution.
	\end{rem}
	
	\begin{rem}
		\rv{It is clear that the MINI--VEM is not a divergence free element for the Stokes equations since $\div(\V_k(\elementvem))$ is not included in $Q_k(\elementvem)$. This feature is shared with the stabilized equal-order formulations introduced in~\cite{guo2020new,li2024stabilized} and with the conforming method discussed in~\cite{MANZINI-MAZZIA}. A divergence free discretization is proposed, for instance, in~\cite{da2017divergence}.}
	\end{rem}
	
	\subsection{Construction of a computable $a_h$}\label{sec:constr_a}
	
	Let us notice that the continuous form $a$ can be written in terms of local contributions as
	$$
	a(\v,\ww) = \sum_{\elementvem\in\mesh_h} a^\elementvem(\v,\ww), \qquad a^\elementvem(\v,\ww) = \int_\elementvem \Grad\v:\Grad\ww\,\dx \qquad \forall\v,\ww\in\V.
	$$
	
	By decomposing the velocity variables as described in \eqref{eq:decomposition}, we can easily write
	\begin{equation}\label{eq:local_a}
		a^\elementvem(\u_h,\v_h) =
		\underbrace{a^\elementvem(\utilde,\vtilde) + a^\elementvem(\bolla_h,\bollatest_h)}_{\text{(I)}} + \underbrace{a^\elementvem(\utilde,\bollatest_h) + a^\elementvem(\bolla_h,\vtilde)}_{\text{(II)}}.
	\end{equation}
	
	The terms $a^\elementvem(\utilde,\vtilde)$ and $a^\elementvem(\bolla_h,\bollatest_h)$ are computable only if at least one of the two entries is a polynomial of degree $k$ and $\rvf{k+1}$, respectively, therefore we employ projections onto polynomial to construct their computable version. By applying the projectors $\PinablaG{k}$ and $\Pinablakpu$ defined in Section~\ref{sec:projectors}, we have
	\begin{equation}\label{eq:pezzi_ah}
		\begin{aligned}
			&a^\elementvem(\utilde,\vtilde) = \int_\elementvem \Grad\Pinabla{k}\utilde:\Grad\Pinabla{k}\vtilde\,\dx + \int_\elementvem \Grad\Qnabla{k}\utilde:\Grad\Qnabla{k}\vtilde\,\dx, \\
			&a^\elementvem(\bolla_h,\bollatest_h) = \int_\elementvem \Grad\Pinablakpu\bolla_h:\Grad\Pinablakpu\bollatest_h\,\dx + \int_\elementvem \Grad\Qnablakpu\bolla_h:\Grad\Qnablakpu\bollatest_h\,\dx,
		\end{aligned}
	\end{equation}
	where we adopted the compact notation $$\Qnabla{k}=\identity-\Pinabla{k},$$ which will be used from now on. It is well known that in both cases the second integral at the right hand side is not computable. Therefore, in order to define $a_h^\elementvem$, we replace the two purely virtual terms with suitable \rv{semi-positive definite} stabilization terms. To this aim, we consider two symmetric bilinear forms $S_k^\elementvem$ and $\Skpu$ satisfying
	\begin{equation}\label{eq:stab_ah}
		\begin{aligned}
			&\alpha_\star\,a^\elementvem(\vtilde,\vtilde) \le S_k^\elementvem(\vtilde,\vtilde) \le \alpha^\star\,a^\elementvem(\vtilde,\vtilde)&&\quad \forall \vtilde\in\enhancedunoS\cap\ker(\PinablaG{k}), \\
			&\gamma_\star\,a^\elementvem(\bollatest_h,\bollatest_h) \le \Skpu(\bollatest_h,\bollatest_h) \le \gamma^\star\,a^\elementvem(\bollatest_h,\bollatest_h)&&\quad \forall \bollatest_h\in\bubbleS\cap\ker(\Pinablakpu)
		\end{aligned}
	\end{equation}
	for real positive constants $\alpha_\star,\alpha^\star,\gamma_\star,\gamma^\star$. Then we replace (I) with
	\begin{equation*}
		a^\elementvem(\Pinabla{k}\utilde,\Pinabla{k}\vtilde)
		+ a^\elementvem(\Pinablakpu\bolla_h,\Pinablakpu\bollatest_h)
		+ S_k^\elementvem(\Qnabla{k}\utilde,\Qnabla{k}\vtilde)
		+ \beta_\sharp\,\Skpu(\Qnablakpu\bolla_h,\Qnablakpu\bollatest_h),
	\end{equation*}
	where $\beta_\sharp\ge0$ is a constant, whose value will be chosen later on.
	
	Let us now consider the interaction between $\vtilde$ and $\bollatest_h$. We see that the computable part is zero. More precisely, we can write
	\begin{equation}\label{eq:cross-zero0}
		a^\elementvem(\vtilde,\bollatest_h) = a^\elementvem(\Pinabla{k}\vtilde,\bollatest_h) + a^\elementvem(\Qnabla{k}\vtilde,\bollatest_h),
	\end{equation}
	where the first term is computable, while the second one is not. If $k=1$, the first term is zero thanks to~\eqref{eq:horto} since $\Pinabla{1}\vtilde$ is an harmonic polynomial. For $k\ge2$, integration by parts gives
	\begin{equation}\label{eq:cross-zero}
		a^\elementvem(\Pinabla{k}\vtilde,\bollatest_h) = -\int_\elementvem \Delta\Pinabla{k}\vtilde\cdot\bollatest_h\,\dx + \int_{\partial\elementvem} (\Grad\Pinabla{k}\vtilde\cdot\n_\elementvem)\,\bollatest_h\,\ds.
	\end{equation}
	The boundary integral vanishes because bubbles are zero at the boundary of each element. On the other hand, since $\Delta\Pinabla{k}\vtilde\in\Polydue{k-2}$, the associated integral is zero because bubbles in $\bubbleS$ have, by construction, null moments up to order $k-2$.
	As (II) is not computable, we choose to neglect such term. We will see that this will not affect the stability and convergence of the method.
	
	Finally, we define the computable discrete bilinear form as
	\begin{equation}\label{eq:ah}
		\begin{aligned}
			a_h^\elementvem(\u_h,\v_h) & =  a^\elementvem(\Pinabla{k}\utilde,\Pinabla{k}\vtilde) + a^\elementvem(\Pinablakpu\bolla_h,\Pinablakpu\bollatest_h)\\
			&\quad+ S_k^\elementvem(\Qnabla{k}\utilde,\Qnabla{k}\vtilde) + \beta_\sharp\,\Skpu(\Qnablakpu\bolla_h,\Qnablakpu\bollatest_h).
		\end{aligned}
	\end{equation}
	
	It is not difficult to see that the discrete bilinear form $a_h^\elementvem$ is $k-$consistent and stable as stated by the following propositions.
	
	\begin{prop}[$k$--consistency]
		For every $\v_h\in\V_k(\elementvem)$ it holds
		\begin{equation}
			a_h^\elementvem(\v_h,\q_k) = a^\elementvem(\v_h,\q_k)\qquad\forall\q_k\in\Polydue{k}.
		\end{equation}
	\end{prop}
	
	\begin{proof} 
		For $\q_k\in\Polydue{k}$, the term $S_k^\elementvem(\Qnabla{k}\vtilde,\Qnabla{k}\q_k)$ vanishes. Notice also that, since ${\q_k\in\V_k(\elementvem)}$ is a polynomial, the associated bubble $\q_\bolla$ is zero, therefore the stabilization term $\Skpu(\Qnablakpu\bollatest_h,\Qnablakpu\q_\bolla)$ vanishes too. This implies that, by exploiting the definition of $\Pinabla{k}$ and $\Pinablakpu$, we have
		\begin{equation*}
			a_h^\elementvem(\v_h,\q_k) = a^\elementvem(\Pinabla{k}\vtilde,\q_k) + a^\elementvem(\Pinablakpu\bollatest_h,\q_k) = a^\elementvem(\v_h,\q_k).
		\end{equation*}
	\end{proof}
	
	\begin{prop}[Stability]\label{prop:stab}
		There exist two real constants $C_\star,C^\star>0$ independent of $h$ such that
		\begin{equation}\label{eq:stab}
			C_\star\,a^\elementvem(\v_h,\v_h) \le a_h^\elementvem(\v_h,\v_h) \le C^\star\,a^\elementvem(\v_h,\v_h)\qquad \forall \v_h\in\V_k(\elementvem).
		\end{equation}
	\end{prop}
	
	\begin{proof}
		We start proving the inequality on the right. \rosso{From the definition~\eqref{eq:ah} and the stabilization properties~\eqref{eq:stab_ah}, we have
			\begin{equation*}
				\begin{aligned}
					a_h^\elementvem(\v_h,\v_h) &= a^\elementvem(\Pinabla{k}\vtilde,\Pinabla{k}\vtilde) + a^\elementvem(\Pinablakpu\bollatest_h,\Pinablakpu\bollatest_h)\\
					&\quad+ S_k^\elementvem(\Qnabla{k}\vtilde,\Qnabla{k}\vtilde) + \beta_\sharp\,\Skpu(\Qnablakpu\bollatest_h,\Qnablakpu\bollatest_h)\\
					&\le a^\elementvem(\Pinabla{k}\vtilde,\Pinabla{k}\vtilde) + a^\elementvem(\Pinablakpu\bollatest_h,\Pinablakpu\bollatest_h)\\
					&\quad+ \alpha^\star\,a^\elementvem(\Qnabla{k}\vtilde,\Qnabla{k}\vtilde)
					+ \beta_\sharp\,\gamma^\star\,a^\elementvem(\Qnablakpu\bollatest_h,\Qnablakpu\bollatest_h).
				\end{aligned}
			\end{equation*}
			By combining the terms at the right hand side as in~\eqref{eq:pezzi_ah} and adding and subtracting $2a^\elementvem(\vtilde,\bollatest_h)$, we get
			\begin{equation}
				\begin{aligned}
					a_h^\elementvem(\v_h,\v_h)
					&\le
					\rv{C} [ a^\elementvem(\vtilde,\vtilde)+a^\elementvem(\bollatest_h,\bollatest_h)+ 2\,a^\elementvem(\vtilde,\bollatest_h)-\rv{2}\,{a^\elementvem(\Qnabla{k}\vtilde,\bollatest_h)}]\\
					&\rv{\le C}\,[a^\elementvem(\v_h,\v_h)-{a^\elementvem(\Qnabla{k}\vtilde,\bollatest_h)}],%
				\end{aligned}
			\end{equation}
			where we took into account that $a^\elementvem(\vtilde,\bollatest_h)=a^\elementvem(\Qnabla{k}\vtilde,\bollatest_h)$ thanks to~\eqref{eq:cross-zero0}--\eqref{eq:cross-zero}.}
		\rev{By applying the Cauchy--Schwarz and Young inequalities to the term $a^\elementvem(\Qnabla{k}\vtilde,\bollatest_h)$, we obtain
			\begin{equation}\label{eq:young}
				|a^\elementvem(\Qnabla{k}\vtilde,\bollatest_h)|
				\le |\Qnabla{k}\vtilde|_{1,\elementvem}|\bollatest_h|_{1,\elementvem}
				\le \frac{1}{2\epsilon}|\Qnabla{k}\vtilde|_{1,\elementvem}^2+\frac{\epsilon}{2}|\bollatest_h|_{1,\elementvem}^2.
			\end{equation}
			Hence, by exploiting that $|\Qnabla{k}\vtilde|_{1,\elementvem}\leq C\,|\vtilde|_{1,\elementvem}$ and Proposition~\ref{prop:splitting}, we obtain
			\begin{equation}
				\begin{aligned}
					a^\elementvem(\v_h,\v_h)-{a^\elementvem(\Qnabla{k}\vtilde,\bollatest_h)}
					&\le |a^\elementvem(\v_h,\v_h)-{a^\elementvem(\Qnabla{k}\vtilde,\bollatest_h)}|\\
					&\le a^\elementvem(\v_h,\v_h) + |a^\elementvem(\Qnabla{k}\vtilde,\bollatest_h)|\\
					&\le a^\elementvem(\v_h,\v_h)+\frac{1}{2\epsilon}|\Qnabla{k}\vtilde|_{1,\elementvem}^2+\frac{\epsilon}{2}|\bollatest_h|_{1,\elementvem}^2\\
					&\le C\,[a^\elementvem(\v_h,\v_h)+|\vtilde|_{1,\elementvem}^2+|\bollatest_h|_{1,\elementvem}^2]\\
					&\le C\,[a^\elementvem(\v_h,\v_h)+|\v_h|_{1,\elementvem}^2]\le C^\star\,a^\elementvem(\v_h,\v_h).
				\end{aligned}
			\end{equation}
		}
		\rv{In order to prove the inequality on the left, we first observe that by exploiting again the equality $a^\elementvem(\vtilde,\bollatest_h)=a^\elementvem(\Qnabla{k}\vtilde,\bollatest_h)$, together with~\eqref{eq:young} for $\epsilon=1$, we obtain
			$$
			\begin{aligned}
				a^\elementvem(\v_h,\v_h)
				&=a^\elementvem(\Pinabla{k}\vtilde,\Pinabla{k}\vtilde) +a^\elementvem(\Qnabla{k}\vtilde,\Qnabla{k}\vtilde)+2\,a^\elementvem(\Qnabla{k}\vtilde,\bollatest_h)\\
				&\quad+ a^\elementvem(\Pinablakpu\bollatest_h,\Pinablakpu\bollatest_h)+a^\elementvem(\Qnablakpu\bollatest_h,\Qnablakpu\bollatest_h)\\
				&\le a^\elementvem(\Pinabla{k}\vtilde,\Pinabla{k}\vtilde) +2\,a^\elementvem(\Qnabla{k}\vtilde,\Qnabla{k}\vtilde)\\
				&\quad+ 2\,a^\elementvem(\Pinablakpu\bollatest_h,\Pinablakpu\bollatest_h)+2\,a^\elementvem(\Qnablakpu\bollatest_h,\Qnablakpu\bollatest_h),
			\end{aligned}
			$$
			so that~\eqref{eq:stab_ah} gives
			\begin{equation*}
				\begin{aligned}
					a^\elementvem(\v_h,\v_h) & \le a^\elementvem(\Pinabla{k}\vtilde,\Pinabla{k}\vtilde) +\frac2{\alpha_\star}S_k^\elementvem(\Qnabla{k}\vtilde,\Qnabla{k}\vtilde)  \\
					&\quad+ 2 a^\elementvem(\Pinablakpu\bollatest_h,\Pinablakpu\bollatest_h) + 2 a^\elementvem(\Qnablakpu\bollatest_h,\Qnablakpu\bollatest_h).
				\end{aligned}
			\end{equation*}
		}
		We conclude by observing that, thanks to Proposition \ref{prop:norm_bub} and \eqref{eq:stab_ah}, we have that, provided $\beta_\sharp \geq 0$, 
		\[
		a^\elementvem(\Qnablakpu\bollatest_h,\Qnablakpu\bollatest_h) \leq C \left(
		a^\elementvem(
		\Pinablakpu \bollatest_h,  \Pinablakpu \bollatest_h
		) + \beta_\sharp \Skpu(\Qnablakpu\bollatest_h,\Qnablakpu\bollatest_h)
		\right),
		\]
		\rv{where $C=1$ for $\beta_\sharp \geq 1/\gamma_\star$, while $C = \gamma_\sharp - 1$ for $\beta_\sharp= 0$}. \rv{Indeed, for $\beta_\sharp= 0$, it holds
			$$
			\begin{aligned}
				a^\elementvem(\Qnablakpu\bollatest_h,\Qnablakpu\bollatest_h)
				&= a^\elementvem(\Qnablakpu\bollatest_h,\Qnablakpu\bollatest_h)+a^\elementvem(\Pinablakpu \bollatest_h,\Pinablakpu \bollatest_h)-a^\elementvem(\Pinablakpu \bollatest_h,\Pinablakpu \bollatest_h)\\
				&\le a^\elementvem(\bollatest_h,\bollatest_h)-a^\elementvem(\Pinablakpu \bollatest_h,\Pinablakpu \bollatest_h)\\
				&\le (\gamma_\sharp-1)\,a^\elementvem(\Pinablakpu \bollatest_h,\Pinablakpu \bollatest_h).
			\end{aligned}
			$$}
	\end{proof}

	\rosso{By summing over all local contributions, we obtain the global discrete bilinear form $a_h$, that is
		\begin{equation}
			a_h(\v_h,\ww_h) = \sum_{\elementvem\in\mesh_h} a_h^\elementvem(\v_h,\ww_h).
		\end{equation}
		Continuity and coercivity derive from Proposition~\ref{prop:stab}.} The following result holds.
	\begin{prop}\label{prop:cont-cor}
		The bilinear form $a_h$ is continuous and coercive, i.e.
		\begin{equation}
			a_h(\v_h,\ww_h)\le C^\star \|\v_h\|_{1,\Omega} \|\ww_h\|_{1,\Omega},\qquad a_h(\v_h,\v_h)\ge C_\star \|\v_h\|_{1,\Omega}^2\qquad\forall\v_h,\ww_h\in\V_h.
		\end{equation}
	\end{prop}
	
	\subsection{Construction of a computable $b_h$}\label{sec:comp_bh}
	
	Also in this case, we split the continuous bilinear form $b$ over the elements of the mesh $\mesh_h$
	\begin{equation}
		b(\v,q)=\sum_{\elementvem\in\mesh_h} b^\elementvem(\v,q), \qquad  b^\elementvem(\v,q)=\int_\elementvem q\,\div\v\,\dx \qquad \forall\v\in\V,q\in Q.
	\end{equation}
	We define the local discrete bilinear form $b_h$ as
	\begin{equation}
		b_h^\elementvem(\v_h,q_h) = b^\elementvem(\v_h,\Pizero{k}q_h)\qquad\forall\v_h\in\V_k(\elementvem),\,q_h\in Q_k(\elementvem).
	\end{equation}
	As done before, by writing $\v_h\in\V_h(\elementvem)$ as $\v_h = \vtilde + \bollatest_h$, we find
	\begin{equation}\label{eq:termini}
		b_h^\elementvem(\v_h,\Pizero{k}q_h) = b_h^\elementvem(\vtilde,\Pizero{k}q_h) + b_h^\elementvem(\bollatest_h,\Pizero{k}q_h).
	\end{equation}
	Integrating by parts the second term at the right hand side, we obtain
	\begin{equation}
		b_h^\elementvem(\bollatest_h,q_h) = \int_\elementvem \Pizero{k}q_h\,\div(\bollatest_h)\,\dx = -\int_\elementvem \grad \Pizero{k}q_h\cdot\bollatest_h\,\dx
	\end{equation}
	since $\bollatest_h$ vanishes on $\partial\elementvem$. Hence, this term is computable as the degrees of freedom of the bubble function $\bollatest_h$ are the moments \rvf{of order $k-1$.}
	
	On the other hand, if we integrate by parts the first term at the right hand side of \eqref{eq:termini}, we get
	\begin{equation}
		b_h^\elementvem(\vtilde,q_h) = \int_\elementvem \Pizero{k}q_h\,\div(\vtilde)\,\dx = -\int_\elementvem \grad \Pizero{k}q_h\cdot\vtilde\,\dx + \int_{\partial\elementvem} (\vtilde\cdot\n_\elementvem)\,\Pizero{k}q_h\,\ds.
	\end{equation}
	Also in this case we have no issues with the computability of this objects: the first integral at the right hand side can be computed thanks to the definition of the enhanced space $\enhancedunoS$, while the boundary integral is computable because $\vtilde$ is a piecewise polynomial on $\partial\elementvem$ and $\Pizero{k}q_h$ is a polynomial by construction.
	
	Notice also that the use of enhanced VEM for $\V_k$ is not required by the construction of $a_h$, but it is necessary to ensure the computability of $b_h$. Moreover, it is immediate to see that, whenever $q_h\in\Poly{k}(\elementvem)$, the consistency property is satisfied and ${b_h^\elementvem(\v_h,q_h) = b^\elementvem(\v_h,q_h)}$.
	
	We then define the global form $b_h$ as
	\begin{equation}
		b_h(\v_h,q_h)=\sum_{\elementvem\in\mesh_h} b_h^\elementvem(\v_h,q_h).
	\end{equation}
	
	\subsection{The discrete problem}\label{sec:discrete_problem}
	
	Before presenting the discrete version of Problem~\ref{pro:stokes} constructed on the spaces $\VV_k=(\V_k,Q_k)$, we introduce the stabilization term $c_h$ dealing with the nonpolynomial contribution of the pressure, which is not playing any role in the definition of $b_h$. The stabilization $c_h$ is defined as
	\begin{equation}
		c_h(p_h,q_h) = \sum_{\elementvem\in\mesh_h} S_p^\elementvem(\Qzero{k}p_h,\Qzero{k}q_h),
	\end{equation}
	where $\Qzero{k} = \identity-\Pizero{k}$. More precisely, $S_p^\elementvem$ is any symmetric bilinear form satisfying the stability estimate
	\begin{equation}\label{eq:stab_sp}
		\delta_\star \|q_h\|_{0,\elementvem}^2 \le S_p^\elementvem(q_h,q_h) \le \delta^\star \|q_h\|_{0,\elementvem}^2 \qquad \forall q_h\in Q_k(\elementvem) \cap \ker(\Pizero{k}).
	\end{equation}
	Thus, $c_h$ satisfies
	\begin{equation}\label{eq:7.35}
		\delta_\star \|\rv{\Qzero{k}}q_h\|_{0,\Omega}^2 \le c_h(q_h,q_h) \le \delta^\star \|\rv{\Qzero{k}}q_h\|_{0,\Omega}^2 \qquad \forall q_h\in Q_k.
	\end{equation}
	
	The discrete version of Problem~\ref{pro:stokes}, \rv{depending on a stabilization parameter $\alpha>0$,} reads as follows.
	\begin{pb}\label{pro:stokes_discreteC}
		Find $(\u_h,p_h)\in\VV_k$ such that
		\begin{equation}
			\begin{aligned}
				a_h(\u_h,\v_h) - b_h(\v_h,p_h) &= (\f_h,\v_h)_\Omega &&\forall\v_h\in\V_k,\\
				b_h(\u_h,q_h) + \alpha\,c_h(p_h,q_h)&= 0 &&\forall q_h\in Q_k.
			\end{aligned}
		\end{equation}
	\end{pb}
	
	We also introduce the discrete counterpart of the global bilinear form $B$ defined as
	\begin{equation}
		{B}_h[(\u_h,p_h),(\v_h,q_h)] = a_h(\u_h,\v_h) - b_h(\v_h,p_h) +b_h(\u_h,q_h) + \alpha\,c_h(p_h,q_h).
	\end{equation}
	Problem~\ref{pro:stokes_discreteC} rewrites as
	\begin{pb}\label{pro:stokes_discrete_stab}
		Find $(\u_h,p_h)\in\VV_k$ such that
		\begin{equation*}
			{B}_h[(\u_h,p_h),(\v_h,q_h)] = (\f_h,\v_h)_\Omega \qquad \forall(\v_h,q_h)\in\VV_k.
		\end{equation*}
	\end{pb}
	
	The vector $\f$ at the right hand side of equation \eqref{pro:stokes} is approximated by its polynomial projection $\f_h=\rvf{\PizeroG{k-1}}\f$ (see e.g.~\cite{beirao2013basic}):
	\rvf{\begin{equation}
			( \f_h, \v_h )_\Omega = \sum_{\elementvem\in\mesh_h} \int_\elementvem \PizeroG{k-1}\f \cdot \v_h\,\dx = \sum_{\elementvem\in\mesh_h} \int_\elementvem \f \cdot \PizeroG{k-1}\v_h\,\dx.
	\end{equation}}
	
	The right hand side is computable: given again $\v_h\in \V_k(\elementvem)$ such that $\v_h=\vtilde+\bollatest_h$, we are allowed to compute $\rvf{\PizeroG{k-1}}\vtilde$ because of the properties of the enhanced VEM space $\enhanced{k}(\elementvem)$, whereas  $\rvf{\PizeroG{k-1}}\bollatest_h$ \rosso{is computable as the moments up to order \rvf{$k-1$} of the bubble space $\bubbleS$ are known, indeed those up to order $k-2$ are zero and the remaining are the dofs.}
	
	Assuming that $\f\in\bH^s(\Omega)$, for $0 \le s \le k$, the following standard error estimate holds
	\begin{equation}\label{eq:rhs_estimate}
		\left|( \f_h, \v_h )_\Omega - (\f,\v_h)_\Omega \right| \le C h^{s+1} \|\f\|_{s,\Omega}|\v_h|_{1,\Omega}.
	\end{equation}
	
	\begin{rem}
		We may consider adding a stabilization term to the discrete divergence form $b_h$ with the aim of removing the need for the pressure stabilization $c_h$. The design of this alternative formulation requires additional effort and will be investigated in our future works.
	\end{rem}
	
	\rv{\begin{rem}
			We chose $\f_h=\rvf{\PizeroG{k-1}}\f$ as it is the most natural choice since the definition of the bubble space descends from plain VEM and the operator $\PizeroG{k}$ is not computable from the degrees of freedom.  On the other hand, the choice $\f_h=\rvf{\PizeroG{k-2}}\f$ kills the bubble contribution due to orthogonality with respect to polynomials of degree $k-2$. More complex quadrature might also by considered for computing the right hand side.
	\end{rem}}
	
	\section{Well-posedness}\label{sec:well_posedness}
	
	In this section, we prove the well-posedness of the discrete problem. To this aim, we first introduce the following proposition, stating a weak inf-sup condition for $b_h$. 
	
	\begin{prop}
		\rosso{For all $q_h\in Q_k$, there exists $\v_h \rv{= \vtilde+\bollatest_h} \in\V_k$ with $\|\v_h\|_{1,\Omega}=\|q_h\|_{0,\Omega}$ such that
			\begin{equation}
				b_h(\v_h,q_h) \ge \zeta_1\,\|q_h\|_{0,\Omega}^2-\zeta_2\, c_h(q_h,q_h),
			\end{equation}
			where $\zeta_1$ and $\zeta_2$ are two positive constants, independent of $h$.
		}
	\end{prop}
	
	\begin{proof}
		We first observe that the discrete space $Q_k$ is a subspace of $\Ldo$. \rv{\rev{Then}, thanks to \eqref{continfsup}, for all $q_h\in Q_k$, there exists $\v(q_h)\in \V$ such that
			\begin{equation}\label{eq:norm_equal_2}
				b(\v(q_h),q_h)\ge C\|\v(q_h)\|_{0,\Omega} \|q_h\|_{0,\Omega}.
			\end{equation}
			Without loss of generality, after possibly rescaling $\v(q_h)$, we can always assume that \begin{equation}\label{eq:normal_equal1}\|\v(q_h)\|_{1,\Omega}=\|q_h\|_{0,\Omega},\end{equation} so that \eqref{eq:norm_equal_2} becomes
			\begin{equation}\label{eq:norm_equal}
				b(\v(q_h),q_h)\ge C \|q_h\|^2_{0,\Omega}.
		\end{equation}}
		Let now $\Pclem{\v}\in\enhanced{k}$ denote \rv{the Cl\'ement interpolant} of $\v(q_h)$ \rv{given by Proposition \ref{clement}}. From the inequality above, we can write
		\begin{equation}
			\begin{aligned}
				C\,\|q_h\|_{0,\Omega}^2 &\le b(\Pclem{\v},q_h) + b(\v(q_h)-\Pclem{\v},q_h) \\
				& \le b_h(\Pclem{\v},q_h)+ \sum_{\elementvem\in\mesh_h}b^\elementvem(\Pclem{\v},q_h - \Pizero{k}q_h) + b(\v(q_h)-\Pclem{\v},q_h),
			\end{aligned}		
		\end{equation}
		where we used the definition of $b_h$. Integration by parts gives
		\begin{equation}
			b(\v(q_h)-\Pclem{\v},q_h) = -\int_\Omega (\v(q_h)-\Pclem{\v})\cdot\grad q_h\,\dx,
		\end{equation}
		where the boundary term vanishes because $\v(q_h)$ and $\Pclem{\v}$ have zero trace on $\partial\Omega$. Hence, we have
		\begin{equation}
			\begin{aligned}
				C\|q_h\|_{0,\Omega}^2 &\le b_h(\Pclem{\v},q_h) + \sum_{\elementvem\in\mesh_h}b^\elementvem(\Pclem{\v},\Qzero{k}q_h) -\int_\Omega (\v(q_h)-\Pclem{\v})\cdot\grad q_h\,\dx\\
				&\rv{= }\,b_h(\Pclem{\v},q_h) + \sum_{\elementvem\in\mesh_h}b^\elementvem(\Pclem{\v},\Qzero{k}q_h)\\
				&\quad - \sum_{\elementvem\in\mesh_h} \int_\elementvem (\v(q_h)-\Pclem{\v})\cdot\grad\Pizero{k}q_h\,\dx- \sum_{\elementvem\in\mesh_h} \int_\elementvem (\v(q_h)-\Pclem{\v})\cdot\grad\Qzero{k}q_h\,\dx,
			\end{aligned}
		\end{equation}
		where, we recall, $\Qzero{k} q_h$ stands for $q_h-\Pizero{k} q_h$. We observe that, by construction of $\Pclem \v$, the moments up to order $k-2$ of the difference $\v(q_h)-\Pclem \v$ are zero. Moreover, $\grad\Pizero{k}q_h\in\Polydue{k-1}$. \rv{Then, letting $\bollatest_h\in\bubbleS$ denote the interpolant of $\v(q_h)-\Pclem \v$ in the bubble space, defined by the conditions
			\[
			\int_{\elementvem} \bollatest_h :\m\,\dx  = 	\int_{\elementvem} (\v(q_h)-\Pclem \v ) :\m\,\dx \qquad \forall\m\in\rvf{\monomialsdue_{k-1}}(\elementvem)\setminus\monomialsdue_{k-2}(\elementvem)\rev{,}
			\]}
		the following equality holds
		\begin{equation}
			\int_\elementvem (\v(q_h)-\Pclem \v)\cdot\grad\Pizero{k}q_h\,\dx = \int_\elementvem \bollatest_h\cdot\grad\Pizero{k}q_h\,\dx.
		\end{equation}
		Integrating again by parts, we find
		\begin{equation}
			-\int_\elementvem \bollatest_h\cdot\grad\Pizero{k}q_h\,\dx = \int_\elementvem \div\bollatest_h\,\Pizero{k}q_h\,\dx
		\end{equation}
		so that
		\begin{equation}
			\begin{aligned}
				{C}\|q_h\|_{0,\Omega}^2 &\le b_h(\Pclem \v,q_h) + b_h(\bollatest_h,q_h)\\
				&\quad + \sum_{\elementvem\in\mesh_h} b^\elementvem(\Pclem \v,\Qzero{k}q_h)- \sum_{\elementvem\in\mesh_h}\int_\elementvem (\v(q_h)-\Pclem \v)\cdot\grad\Qzero{k}q_h\,\dx.
			\end{aligned}
		\end{equation}
		\rv{Applying Young's inequality yields, for $\varepsilon>0$,
			\begin{equation}
				\begin{aligned}
					{C}\|q_h\|_{0,\Omega}^2 &\le b_h(\Pclem \v+\bollatest_h,q_h) + \varepsilon\|\Pclem \v\|_{1,\Omega}^2 + \frac{c}\varepsilon\|\Qzero{k}q_h\|_{0,\Omega}^2\\
					&\quad+\frac{\varepsilon}{h^2}\|\v(q_h)-\Pclem \v\|_{0,\Omega}^2 + \rv{ \frac{c'}\varepsilon} h^2 \|\Qzero{k}q_h\|_{1,h}^2,
				\end{aligned}
			\end{equation}
			where $c$ and $c'$ are two positive constants independent of $h$ and $\varepsilon$.
			Setting $\v_h = \Pclem \v + \bollatest_h$, together with \eqref{clementapprox},\eqref{eq:normal_equal1} and \eqref{eq:norm_equal},  this eventually gives a bound of the form
			\begin{equation}
				c_0\|q_h\|_{0,\Omega}^2 \le b_h(\v_h,q_h) + \frac{c_1}\varepsilon	\|\Qzero{k}q_h\|_{0,\Omega}^2 + c_2 \varepsilon\|q_h\|_{0,\Omega}^2,
			\end{equation}
			$\varepsilon > 0$ arbitrary, and with $c_0$, $c_1$ and $c_2$ positive constants also independent of $h$ and $\varepsilon$. Subtracting $c_2\varepsilon \| \rev{q_h} \|_{0,\Omega}^2$ from both sides and using \eqref{eq:7.35} we obtain
			\begin{equation}
				(c_0 - c_2 \varepsilon)\|q_h\|_{0,\Omega}^2 \le b_h(\v_h,q_h)+ \frac{c_1}\varepsilon \delta_\star^{-1}\ c_h(q_h,q_h).
			\end{equation}
			Finally, choosing $\varepsilon = c_0/ (2 c_2)$ and subtracting from both sides the last term on the right hand side we obtain the  desired bound with $\zeta_1=c_0/2$ and $\zeta_2=2 c_1 c_2 /(c_0 \delta_\star)$.}
	\end{proof}

	The well-posedness result can then be proved.
	
	\begin{prop}
		For any $(\u_h,p_h)\in\VV_k$ there exists a positive constant $\omega$ independent of $h$ such that
		\begin{equation}\label{eq:Btilde_stab}
			\sup_{(\v_h,q_h)\in\VV_k} \frac{{B}_h[(\u_h,p_h),(\v_h,q_h)]}{\|(\v_h,q_h)\|_{\VV}} \ge \omega\,\|(\u_h,p_h)\|_{\VV}.
		\end{equation}
		Moreover, ${B}_h$ is continuous, i.e. it satisfies
		\begin{equation}\label{eq:Btilde_cont}
			{B}_h[(\u_h,p_h),(\v_h,q_h)] \le C \, \|(\u_h,p_h)\|_{\VV} \|(\v_h,q_h)\|_{\VV}.
		\end{equation}
	\end{prop}
	
	\begin{proof}
		
		The result is proved by adapting~\cite[Theorem 4.5]{guo2020new}.
		
	\end{proof}
	
	{As a consequence, the discrete problem admits a unique solution.}
	
	\begin{thm}
		The discrete Problem~\ref{pro:stokes_discreteC}/\ref{pro:stokes_discrete_stab} has a unique solution ${(\u_h,p_h)\in(\V_k,Q_k)}$ satisfying
		\begin{equation}
			\|\u_h\|_{1,\Omega}+\|p_h\|_{0,\Omega} \le C \|\f\|_{0,\Omega}.
		\end{equation}
	\end{thm}

	\section{Error analysis}\label{sec:error_analysis}
	
	In this section, we present the convergence analysis for the proposed MINI--VEM for Stokes. We derive estimates in both the $H^1$ and $L^2$ norm for the velocity and in $L^2$ norm for the pressure. 
	
	\subsection{Main error estimates for velocity and pressure}
	
	\begin{thm}\label{thm:est}
		Let $(\u,p)\in\VV$ be the solution of the continuous Stokes Problem~\ref{pro:stokes}/\ref{pro:stokes2} and $(\u_h,p_h)\in\VV_k$ the solution of the discrete Problem~\ref{pro:stokes_discreteC}/\ref{pro:stokes_discrete_stab}. The following estimate holds true
		\begin{equation}\label{eq:est}
			\begin{aligned}
				\rv{|\u-\u_h|_{1,\Omega}} + \|p-p_h\|_{0,\Omega} &\le C\,( \rv{|\u-\u_I|_{1,\Omega}} + | \u-\Pinabla{k}\u |_{1,h} + \|\div\u-\Pizero{k}\div\u\|_{0,\Omega}\\
				&\qquad+ \|p-p_I\|_{0,\Omega}  + \|p-\Pizero{k}p\|_{0,\Omega} + \|\f-\f_h\|_{0,\Omega} ),
			\end{aligned}
		\end{equation}
		where $\u_I$ is the interpolant of $\u$ in $\enhancedunoG\subset\V_k$ and $p_I$ is the interpolant of $p$ in $Q_k$.
	\end{thm}
	
	\begin{proof}
		We start by applying the triangle inequality:
		\begin{equation}\label{eq:tri_in_1}
			\rv{| \u - \u_h |_{1,\Omega}} + \| p - p_h \|_{0,\Omega} \le \rv{| \u - \u_I |_{1,\Omega}} + \rv{| \u_I - \u_h |_{1,\Omega}} + \| p - p_I \|_{0,\Omega} + \| p_I - p_h \|_{0,\Omega}.
		\end{equation}
		In the following we estimate the terms $| \u_I - \u_h |_{1,\Omega}$, $\| p_I - p_h \|_{0,\Omega}$.
		
		For the discrete bilinear form $a_h$, we obtain the following equality by exploiting the $k$--consistency property
		\begin{equation}
			\begin{aligned}
				a_h(\u_I-\u_h,\v_h) &= \sum_{\elementvem\in\mesh_h} [a_h^\elementvem(\u_I,\v_h)] - a_h(\u_h,\v_h)\\
				&= \sum_{\elementvem\in\mesh_h} [a_h^\elementvem(\u_I-\Pinabla{k}\u,\v_h) \,\rv{+}\, a_h^\elementvem(\Pinabla{k}\u,\v_h)] - a_h(\u_h,\v_h)\\
				&= \sum_{\elementvem\in\mesh_h} [a_h^\elementvem(\u_I-\Pinabla{k}\u,\v_h)\,\rv{+}\, a^\elementvem(\Pinabla{k}\u,\v_h)] - a_h(\u_h,\v_h).
			\end{aligned}
		\end{equation}
		Then, by simple manipulations, we find
		\begin{equation}\label{eq:a_1}
			\begin{aligned}
				a_h(\u_I-\u_h,\v_h) &
				=\sum_{\elementvem\in\mesh_h} [a_h^\elementvem(\u_I-\u,\v_h) + a_h^\elementvem(\u-\Pinabla{k}\u,\v_h)-a^\elementvem(\u-\Pinabla{k}\u,\v_h)]\\
				&\qquad+ a(\u,\v_h)- a_h(\u_h,\v_h).
			\end{aligned}
		\end{equation}
		
		For $b_h$, it holds instead
		\begin{equation}\label{eq:b_21}
			\begin{aligned}
				b_h(\v_h,p_I-p_h) &= \sum_{\elementvem\in\mesh_h}[(\div\v_h,\Pizero{k}p_I)_{\elementvem}] - b_h(\v_h,p_h)\\
				&= \sum_{\elementvem\in\mesh_h}[(\div\v_h,\Pizero{k}(p_I-p))_{\elementvem}+(\div\v_h,\Pizero{k} p-p)_\elementvem]\\
				&\quad+ b(\v_h,p) - b_h(\v_h,p_h),
			\end{aligned}
		\end{equation}
		and, by exploiting the following relation stemming from the definition of $\Pizero{k}$
		\begin{equation}
			(\Pizero{k}\div\u,q_h)_\elementvem = (\Pizero{k}\div\u,\Pizero{k}q_h)_\elementvem = (\div\u,\Pizero{k}q_h)_\elementvem,
		\end{equation}
		we can easily write
		\begin{equation}
			\begin{aligned}\label{eq:b_22}
				b_h(\u_I-\u_h,q_h) &= \sum_{\elementvem\in\mesh_h}[(\div \u_I,\Pizero{k}q_h)_\elementvem] - b_h(\u_h,q_h)\\
				&= \sum_{\elementvem\in\mesh_h}[(\div\u_I-\div\u,\Pizero{k}q_h)_\elementvem] + b(\u,q_h) - b_h(\u_h,q_h)\\
				&\quad \rv{-}\sum_{\elementvem\in\mesh_h}[(\div\u-\Pizero{k}\div\u,q_h)_\elementvem].\\
			\end{aligned}
		\end{equation}
		Moreover, for the stabilization term $c_h$, we have
		\begin{equation}\label{eq:c_2}
			c_h(p_I-p_h,q_h) = c_h(p_I,q_h) - c_h(p_h,q_h).
		\end{equation} 
		
		By summing up~\eqref{eq:a_1}, \eqref{eq:b_21}, \eqref{eq:b_22}, \eqref{eq:c_2} and taking into account continuous and discrete problems, it holds
		\begin{equation}
			\begin{aligned}
				a(\u,\v_h)- a_h(\u_h,\v_h) &- b(\v_h,p) + b_h(\v_h,p_h)\\
				&+ b(\u,q_h) - b_h(\u_h,q_h) - \alpha\,c_h(p_h,q_h)= (\f-\f_h,\v_h)_\Omega,
			\end{aligned}
		\end{equation}
		and
		\begin{equation*}
			\begin{aligned}
				{B}_h[(\u_I-\u_h,p_I-p_h),(\v_h,q_h)] = &\sum_{\elementvem\in\mesh_h} [a_h^\elementvem(\u_I-\u,\v_h) + a_h^\elementvem(\u-\Pinabla{k}\u,\v_h)-a^\elementvem(\u-\Pinabla{k}\u,\v_h)]\\
				&-\sum_{\elementvem\in\mesh_h}[(\div\v_h,\Pizero{k}(p_I-p))_{\elementvem}]-\sum_{\elementvem\in\mesh_h}[(\div\v_h,\Pizero{k} p-p)_\elementvem]\\
				&+ \sum_{\elementvem\in\mesh_h}[(\div\u_I-\div\u,\Pizero{k}q_h)_\elementvem]
				\,\rv{-} \sum_{\elementvem\in\mesh_h}[(\div\u-\Pizero{k}\div\u,q_h)_\elementvem]\\
				&+ \alpha\,c_h(p_I,q_h) + (\f-\f_h,\v_h)_\Omega.
			\end{aligned}
		\end{equation*}
		Notice that, by projection estimate (see Lemma~\ref{lem:proj_err}) and the triangle inequality, we have
		\begin{equation}
			\begin{aligned}
				c_h(p_I,q_h) &\le \delta^\star\,\|p_I-\Pizero{k}p_I\|_{0,\Omega}\|q_h-\Pizero{k}q_h\|_{0,\Omega}\le C \, \|p_I-\Pizero{k}p_I\|_{0,\Omega}\|q_h\|_{0,\Omega}\\
				&\le C\,(\|p_I-p\|_{0,\Omega}+\|p-\Pizero{k}p\|_{0,\Omega})\|q_h\|_{0,\Omega}.
			\end{aligned}
		\end{equation}
		At this point, we obtain
		\begin{equation}
			\begin{aligned}
				&{B}_h[(\u_I-\u_h,p_I-p_h),(\v_h,q_h)]\\
				&\qquad\le C\,( \rv{|\u-\u_I |_{1,\Omega}} + | \u-\Pinabla{k}\u |_{1,h} + \|\div\u-\Pizero{k}\div\u\|_{0,\Omega}\\
				&\qquad\qquad+ \|p-p_I\|_{0,\Omega}  + \|p-\Pizero{k}p\|_{0,\Omega} + \|\f-\f_h\|_{0,\Omega} )\,(|\v_h |_{1,\Omega}+\|q_h\|_{0,\Omega}),
			\end{aligned}
		\end{equation}
		and, as a consequence of the well-posedness,
		\begin{equation}\label{eq:last}
			\begin{aligned}
				\rv{|\u_I-\u_h|_{1,\Omega}} + \|p_I-p_h\|_{0,\Omega} &\le C\,( |\u-\u_I|_{1,\Omega} + | \u-\Pinabla{k}\u |_{1,h} + \|\div\u-\Pizero{k}\div\u\|_{0,\Omega}\\
				&\qquad+ \|p-p_I\|_{0,\Omega}  + \|p-\Pizero{k}p\|_{0,\Omega} + \|\f-\f_h\|_{0,\Omega} ).
			\end{aligned}
		\end{equation}
		Now, if we combine~\eqref{eq:last} with~\rv{\eqref{eq:tri_in_1}}, we find~\eqref{eq:est}. 
	\end{proof}
	
	\rosso{The following Corollary is a direct consequence of Theorem~\ref{thm:est} combined \rv{with the Poincar\'e inequality for  $H^1_0(\Omega)$}, the estimate for the right hand side (see~\eqref{eq:rhs_estimate}), the interpolation estimates in Propositions~\ref{prop:interp_enh} and \ref{clement}, and the projection estimates in Lemma~\ref{lem:proj_err}.}
	
	\begin{cor}
		In the same framework of Theorem~\ref{thm:est}, if $\u\in\bH^{s+1}(\Omega)$, $p\in\mathrm{H}^s(\Omega)$ and $\f\in \bH^s(\Omega)$, $0 \leq s \leq k$, we also have
		\begin{equation}\label{eq:est2}
			\| \u - \u_h \|_{1,\Omega} + \| p - p_h \|_{0,\Omega} \le C\,h^s\, ( |\u|_{s+1,\Omega} + |p|_{s,\Omega} + |\f|_{s,\Omega}).
		\end{equation}
	\end{cor}
	
	\subsection{Error estimates for the $L^2$ norm of the velocity}
	
	The following theorem asserts estimates for the error $\| \u-\u_h \|_{0,\Omega}$.
	
	\begin{thm}
		\rosso{Assume that the domain $\Omega$ is convex.} Let $(\u,p)\in\VV$ be the solution of the continuous Stokes Problem~\ref{pro:stokes}/\ref{pro:stokes2} and $(\u_h,p_h)\in\VV_k$ the solution to Problem~\ref{pro:stokes_discreteC}/\ref{pro:stokes_discrete_stab}. If $\u\in\bH^{s+1}(\Omega)$, $p\in\mathrm{H}^s(\Omega)$ and $\f\in \bH^s(\Omega)$, $0 \leq s \leq k$, the following estimate holds true
		\begin{equation}\label{eq:l2_vel}
			\|\u-\u_h\|_{0,\Omega} \le Ch^{s+1}\,(|\u|_{s+1,\Omega}+|p|_{s,\Omega}+|\f|_{s,\Omega}).
		\end{equation}
	\end{thm}
	
	\begin{proof}
		
		In order to derive the $L^2$ error estimate for the velocity, we resort to the usual duality argument. We denote by $(\ppsi,\rho)\in [\bH^2(\Omega)\cap \dueHunozeroO]\times[\mathrm{H}^1(\Omega)\cap\Ldo]$ the solution of the following problem.
		
		\begin{equation}\label{eq:stokes_dual}
			\begin{aligned}
				-\Delta\ppsi - \grad \rho &= \u-\u_h && \text{in }\Omega, \\
				\div\ppsi &=0&& \text{in }\Omega, \\
				\ppsi&=\mathbf{0}&& \text{on }\partial\Omega.
			\end{aligned}
		\end{equation}
		
		Notice that, thanks to the convexity of the domain $\Omega$, the pair $(\ppsi,\rho)$ satisfies
		\begin{equation}
			\| \ppsi \|_{2,\Omega} + \|\rho\|_{1,\Omega} \le C\,\|\u-\u_h\|_{0,\Omega}.
		\end{equation}
		
		The dual problem in variational formulation reads
		\begin{equation}\label{dual_var}
			B[(\v,q),(\ppsi,\rho)] = (\u-\u_h,\v)_\Omega\qquad\forall(\v,q)\in\VV.
		\end{equation}
		From standard theory, this admits a unique solution $(\u,p)\in\VV$.
		
		We consider the interpolant $\ppsi_I$ of $\ppsi$ in $\enhancedunoG$ and the interpolant $\rho_I$ of $\rho$ in $Q_k$. They satisfy the following properties
		\begin{equation}\label{eq:interp_dual}
			\begin{aligned}
				\| \ppsi-\ppsi_I \|_{1,\Omega} &\le Ch\,\|\ppsi\|_{2,\Omega} \le Ch\,\|\u-\u_h\|_{0,\Omega},\\
				\| \rho-\rho_I \|_{0,\Omega} &\le Ch\,\|\rho\|_{1,\Omega} \le Ch\,\|\u-\u_h\|_{0,\Omega}.
			\end{aligned}
		\end{equation}
		
		From~\eqref{dual_var}, it is easy to see that, by taking $(\v,q)=(\u-\u_h,0)$, we obtain
		\begin{equation}\label{eq:starting_pooint}
			\begin{aligned}
				\|\u-\u_h\|_{0,\Omega}^2 &= a(\u-\u_h,\ppsi) + b(\u-\u_h,\rho)\\
				&= a(\u-\u_h,\ppsi-\ppsi_I) + b(\u-\u_h,\rho-\rho_I) + a(\u-\u_h,\ppsi_I)+ b(\u-\u_h,\rho_I).
			\end{aligned}
		\end{equation}
		
		We treat each term separately. First, by continuity and~\eqref{eq:interp_dual}, we have
		\begin{equation}\label{eq:I}
			\begin{aligned}
				a(\u-\u_h,\ppsi-\ppsi_I) &\le C\,\|\u-\u_h\|_{1,\Omega}\|\ppsi-\ppsi_I\|_{1,\Omega}\\
				&\le Ch\,\|\u-\u_h\|_{1,\Omega}\|\ppsi\|_{2,\Omega}\\
				&\le Ch\,\|\u-\u_h\|_{1,\Omega}\|\u-\u_h\|_{0,\Omega},
			\end{aligned}
		\end{equation}
		and, by applying the same strategy, we also find
		\begin{equation}\label{eq:II}
			\begin{aligned}
				b(\u-\u_h,\rho-\rho_I) &\le C\, \|\u-\u_h\|_{1,\Omega}\|\rho-\rho_I\|_{0,\Omega}\\
				&\le Ch\,\|\u-\u_h\|_{1,\Omega}\|\rho\|_{1,\Omega}\\
				&\le Ch\,\|\u-\u_h\|_{1,\Omega}\|\u-\u_h\|_{0,\Omega}.
			\end{aligned}
		\end{equation}
		
		In addition, by exploiting the definition of $B$ and ${B}_h$, it easy to see that the last two terms in~\eqref{eq:starting_pooint} yield
		\begin{equation}
			\begin{aligned}
				a(\u-\u_h,\ppsi_I)+ b(\u-\u_h,\rho_I) = &a_h(\u_h,\ppsi_I) - a(\u_h,\ppsi_I)\\
				&+ b(\ppsi_I,p) - b_h(\ppsi_I,p_h)\\
				&+ b_h(\u_h,\rho_I) - b(\u_h,\rho_I)\\
				&+ \rv{\alpha}\,c_h(p_h,\rho_I) + (\f-\f_h,\ppsi_I)_\Omega.
			\end{aligned}
		\end{equation}
		
		Before looking for an estimate for the term $a_h(\u_h,\ppsi_I) - a(\u_h,\ppsi_I)$, we introduce terms containing the projectors $\Pinabla{k}\u$, $\Pinabla{k}\ppsi$. By applying $k$-consistency of the discrete bilinear form, we can write
		\begin{equation*}
			\begin{aligned}
				a_h(\u_h,\ppsi_I) - a(\u_h,\ppsi_I) &= \sum_{\elementvem\in\mesh_h} [a_h^\elementvem(\u_h-\Pinabla{k}\u,\ppsi_I) - a^\elementvem(\u_h-\Pinabla{k}\u,\ppsi_I)]\\
				&= \sum_{\elementvem\in\mesh_h} [a_h^\elementvem(\u_h-\Pinabla{k}\u,\ppsi_I-\Pinabla{k}\ppsi) - a^\elementvem(\u_h-\Pinabla{k}\u,\ppsi_I-\Pinabla{k}\ppsi)],\\
			\end{aligned}
		\end{equation*}
		which implies
		\begin{equation*}
			\rosso{
				|a_h(\u_h,\ppsi_I) - a(\u_h,\ppsi_I)| \le C\, |\u_h-\Pinabla{k}\u|_{1,h}|\ppsi_I-\Pinabla{k}\ppsi|_{1,h}.}
		\end{equation*}
		Now, by triangular inequality, the projection estimates and~\eqref{eq:interp_dual}, it holds
		\begin{equation*}
			\rosso{|\ppsi_I-\Pinabla{k}\ppsi|_{1,h}} \le \|\ppsi_I-\ppsi \|_{1,\Omega} + \|\ppsi-\Pinabla{k}\ppsi\|_{1,h} \le Ch\,\|\ppsi\|_{2,\Omega}\le Ch\,\|\u-\u_h\|_{0,\Omega},
		\end{equation*}
		\rosso{and
			\begin{equation}
				|\u_h-\Pinabla{k}\u|_{1,h}
				\le
				|\u-\u_h|_{1,\Omega}
				+|\u-\Pinabla{k}\u|_{1,h},
			\end{equation}	
		}
		therefore
		\begin{equation}\label{eq:III}
			|a_h(\u_h,\ppsi_I) - a(\u_h,\ppsi_I)| \le Ch\,\rosso{( |\u-\u_h|_{1,\Omega}
				+|\u-\Pinabla{k}\u|_{1,h})}\|\u-\u_h\|_{0,\Omega}.
		\end{equation}
		
		For the  next term, we make use of $\div\ppsi=0$ and we apply again~\eqref{eq:interp_dual}, hence
		\begin{equation}
			\begin{aligned}
				b(\ppsi_I,p) - b_h(\ppsi_I,p_h) &= b_h(\ppsi_I,p-p_h) + b(\ppsi_I,p) - b_h(\ppsi_I,p)\\
				& = b_h(\ppsi_I,p-p_h) + \sum_{\elementvem\in\mesh_h}b^\elementvem(\ppsi_I,p-\Pizero{k}p)\\
				& = b_h(\ppsi_I-\ppsi,p-p_h) + \sum_{\elementvem\in\mesh_h}b^\elementvem(\ppsi_I-\ppsi,p-\Pizero{k}p)\\
			\end{aligned}
		\end{equation}
		so that
		\begin{equation}
			\begin{aligned}\label{eq:IV}
				|b(\ppsi_I,p) - b_h(\ppsi_I,p_h)| &\le C\,(\|p-p_h\|_{0,\Omega}+\|p-\Pizero{k}p\|_{0,\Omega})\|\ppsi-\ppsi_I\|_{1,\Omega}\\
				&\le Ch\,(\|p-p_h\|_{0,\Omega}+\|p-\Pizero{k}p\|_{0,\Omega})\|\u-\u_h\|_{0,\Omega}.
			\end{aligned}
		\end{equation}
		With similar computations and taking into account that $\div\u=0$, we also have
		\begin{equation*}
			\begin{aligned}
				b_h(\u_h,\rho_I) - b(\u_h,\rho_I) &= b_h(\u_h,\rho_I-\rho) + \sum_{\elementvem\in\mesh_h} b^\elementvem(\u_h,\Pizero{k}\rho-\rho_I)\\
				&= b_h(\u_h,\rho_I-\rho) + \sum_{\elementvem\in\mesh_h} b^\elementvem(\u_h-\u,\Pizero{k}\rho-\rho_I)\\
				&= b_h(\u_h,\rho_I-\rho) + \sum_{\elementvem\in\mesh_h} b^\elementvem(\u_h-\u,\Pizero{k}\rho-\rho) + \sum_{\elementvem\in\mesh_h} b^\elementvem(\u_h-\u,\rho-\rho_I).\\
			\end{aligned}
		\end{equation*}
		Consequently,
		\begin{equation}\label{eq:V}
			\begin{aligned}
				|b_h(\u_h,\rho_I) - b(\u_h,\rho_I)| &\le C\,(2\|\rho_I-\rho\|_{0,\Omega}+\|\Pizero{k}\rho-\rho\|_{0,\Omega})\|\u_h-\u\|_{1,\Omega}\\
				&\le Ch\,\|\rho\|_{1,\Omega}\|\u_h-\u\|_{1,\Omega}\\
				&\le Ch\,\|\u_h-\u\|_{1,\Omega}\|\u_h-\u\|_{0,\Omega}.
			\end{aligned}
		\end{equation}
		
		At this point, it remains to estimate the stabilization term and the right hand side. For the stabilization term we have
		\begin{equation}\label{eq:VI}
			\begin{aligned}
				|\rv{\alpha}\,c_h(p_h,\rho_I)| &= |-\rv{\alpha}\,c_h(p-p_h,\rho_I)+\rv{\alpha}\,c_h(p,\rho_I)|\\
				&\le \rv{\alpha}\,\delta^\star \sum_{\elementvem\in\mesh_h} \big(\|\Qzero{k}(p-p_h)\|_{0,\elementvem}  + \|\Qzero{k}p\|_{0,\elementvem}\big)\|\Qzero{k}\rho_I\|_{0,\elementvem}\\
				&\le Ch\,\big(\|p-p_h\|_{0,\Omega}+\|p-\Pizero{k}p\|_{0,\Omega}\big)\|\u-\u_h\|_{0,\Omega}.
			\end{aligned}
		\end{equation}
		On the other hand, for the right hand side, we obtain
		\begin{equation}\label{eq:VII}
			\begin{aligned}
				|(\f-\f_h,\ppsi_I)_\Omega| &= |(\f-\f_h,\ppsi_I-\Pinabla{0}\ppsi_I)_\Omega|\\
				&\le Ch\,\|\f-\f_h\|_{0,\Omega} \|\ppsi_I\|_{1,\Omega}\\
				& \le Ch\,\|\f-\f_h\|_{0,\Omega} \|\u-\u_h\|_{0,\Omega}.
			\end{aligned}
		\end{equation}
		
		Finally, by combining the bounds \eqref{eq:I}, \eqref{eq:II}, \eqref{eq:III}, \eqref{eq:IV}, \eqref{eq:V}, \eqref{eq:VI}, \eqref{eq:VII}, we end up with
		\begin{equation}
			\begin{aligned}
				\| \u-\u_h \|_{0,\Omega} \le Ch\,(&\rosso{\|\u-\u_h\|_{1,\Omega} + |\u-\Pinabla{k}\u|_{1,h}}\\
				&+\| p-p_h\|_{0,\Omega} + \|p-\Pizero{k}p\|_{0,\Omega} + \|\f-\f_h\|_{0,\Omega}).
			\end{aligned}
		\end{equation}	
		The proof is completed taking into account the error estimate given by Theorem~\ref{thm:est} and the usual projection estimates.		
	\end{proof}
	
	\begin{rem}
		\rosso{We observe that the bubble contribution to the velocity does not play any role in the best approximation estimate.}
		Indeed, bubbles are just a mathematical tool aimed at stabilizing the formulation, and they do not improve the convergence property to the method. 
	\end{rem}
	
	\section{Numerical tests}\label{sec:num_tests}
	
	We analyze the numerical performance of the proposed method. This section is divided into two parts. First, we study the condition number of the matrix arising from the~$\text{MINI--VEM}$~discretization in dependence of the choice of polynomial basis and varying the value of the pressure stabilizing parameter $\alpha$ in the formulation of Problem~\ref{pro:stokes_discreteC}. In the second part, we show some convergence results confirming the theoretical estimates presented in the previous section.
	
	The numerical tests are performed by considering four different geometric discretizations of the domain. We construct the discrete bilinear form $a_h$ without considering the bubble stabilization term $\Skpu$, i.e. we set $\beta_\sharp=0$ in~\eqref{eq:ah}, which is an admissible choice as explained and proved in Section~\ref{sec:conforming_vem}. The stabilization term $S_k^\elementvem$ for contribution of order $k$ to the discrete velocity and the pressure stabilization $S_p^\elementvem$ are both chosen to be the well-known \textit{dofi--dofi stabilization}, which is generically defined as
	\begin{equation*}
		S_{\mathrm{dofi}}^\elementvem(u,v) = \sum_{i = 1}^{N_{dof}} \mathrm{dof}_i(u)\mathrm{dof}_i(v),
	\end{equation*}
	where $N_{dof}$ is the number of VEM degrees of freedom in the element $\elementvem$.
	
	As usual in VEM literature, the error computation is performed by means of polynomial projections. More precisely, we consider the errors
	\begin{equation}\label{eq:errors}
		\begin{gathered}
			\err{0}{\u_h} = \frac{\|\u-\Pizero{k}\u_h\|_{0,\Omega}}{\|\u\|_{0,\Omega}},\qquad
			\err{1}{\u_h} = \frac{|\u-\Pizero{k}\u_h|_{1,\Omega}}{|\u|_{1,\Omega}},\\
			\err{0}{p_h} = \frac{\|p-\Pizero{k}p_h\|_{0,\Omega}}{\|p\|_{0,\Omega}}.
		\end{gathered}
	\end{equation}
	We remark again that we are going to compute the error just projecting the discrete solution in~$\Poly{k}$: we are dropping out the contribution given by the bubbles, which are used just as a mathematical tool for stabilizing the discrete formulation.
	
	\subsection{Condition number and pressure stabilization}
	
	In this subsection, we analyze two features of the proposed virtual element method.
	
	We study how the choice of the weight $\alpha$ in front of the pressure stabilization term affects the invertibility of the system matrix. To this aim, we measure the condition number of the system varying the value of $\alpha$ from $10^{-15}$ to $10^3$.
	
	At the same time, we also study how the choice of basis $\monomials_{k}$ for the polynomial space $\Poly{k}$ affects the condition number: it is well known (see e.g.~\cite{mascotto2018ill}) that high order approximations require a suitable choice to avoid ill conditioning of the system. Our investigation is motivated by the fact that the MINI--VEM of degree $k$ may behave like a plain VEM of degree $\rvf{k+1}$, since it is constructed by means bubble functions of degree~$\rvf{k+1}$. More precisely, we compare the condition number of the system for $\monomials_{k}(\elementvem)$ being either the standard choice of scaled monomials introduced in~\cite{beirao2013basic}, i.e.
	\begin{equation*}
		\mathcal{M}_{k}(\elementvem) = \rv{\left\{ m=\left(\frac{\x-\x_\elementvem}{h_\elementvem}\right)^\mdeg \,\text{with }\mdeg\in\{0,1,\dots,k\}^2\text{ and }\mdeg_1+\mdeg_2\le k\right\},}
	\end{equation*}
	\rv{where for $\x=(x_1,x_2)$, $\x^\mdeg=x_1^{\mdeg_1}x_2^{\mdeg_2}$}, or the $L^2$ orthonormal basis $\mathcal{Q}_{k}(\elementvem)$ presented in~\cite{mascotto2018ill}, which is obtained by applying a Gram--Schmidt orthonormalization process to $\mathcal{M}_{k}(\elementvem)$. \rv{Due to the arbitrary shape of the meshes we are going to consider, we clarify that the centroid $\x_\elementvem=(x_{\elementvem,1},x_{\elementvem,2})$ of $\elementvem$ is computed as
		$$
		x_{\elementvem,i} = \frac{1}{|\elementvem|} \int_\elementvem x_i\,\dx\qquad\qquad\text{for}\quad i=1,2.
		$$}
	
\begin{figure}
		\subfloat[Hexagons]{\includegraphics[trim = 40 0 40 0, width=0.23\linewidth]{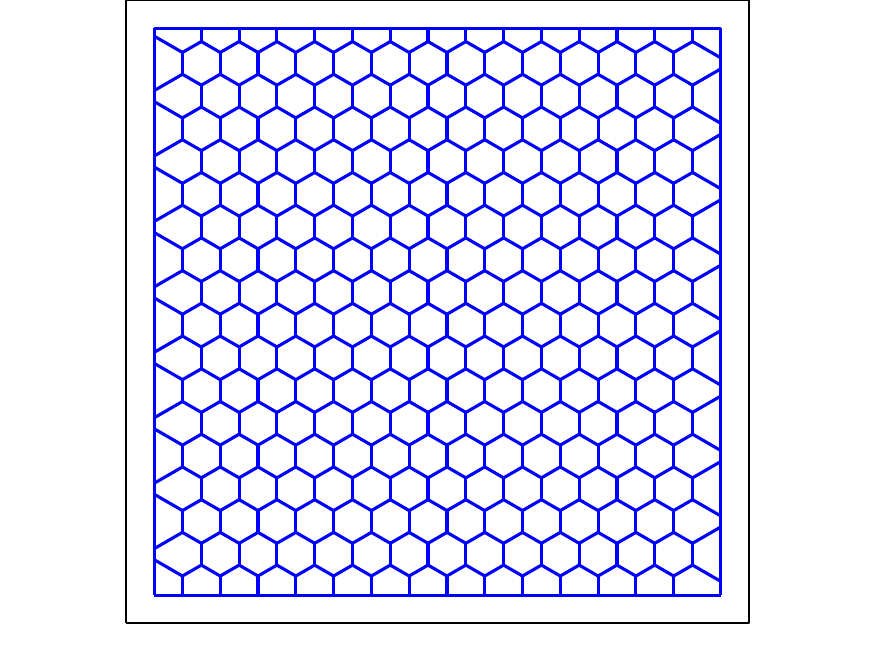}}\quad
		\subfloat[Voronoi]{\includegraphics[trim = 40 0 40 0, width=0.23\linewidth]{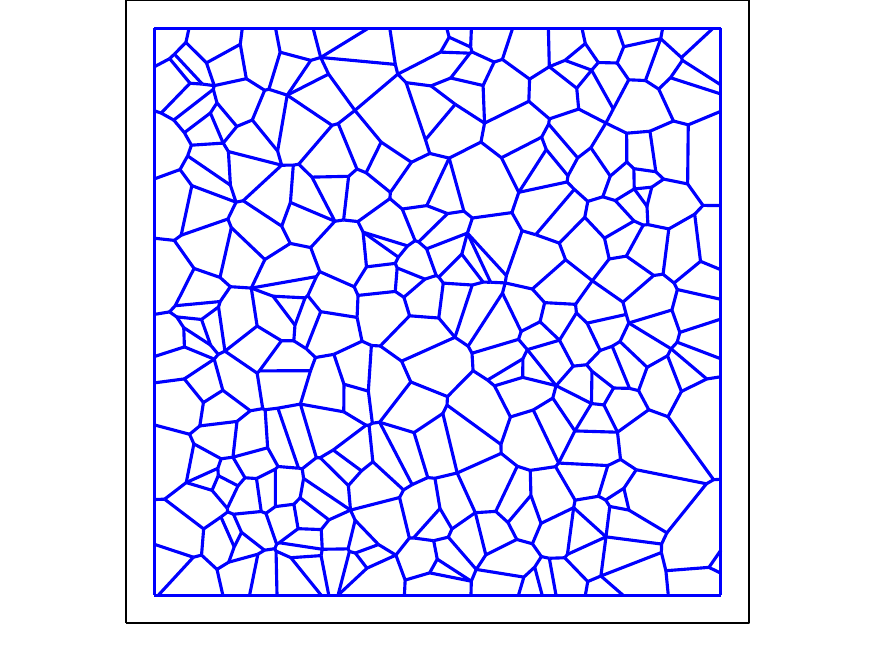}}\quad
		\subfloat[Random]{\includegraphics[trim = 40 0 40 0, width=0.23\linewidth]{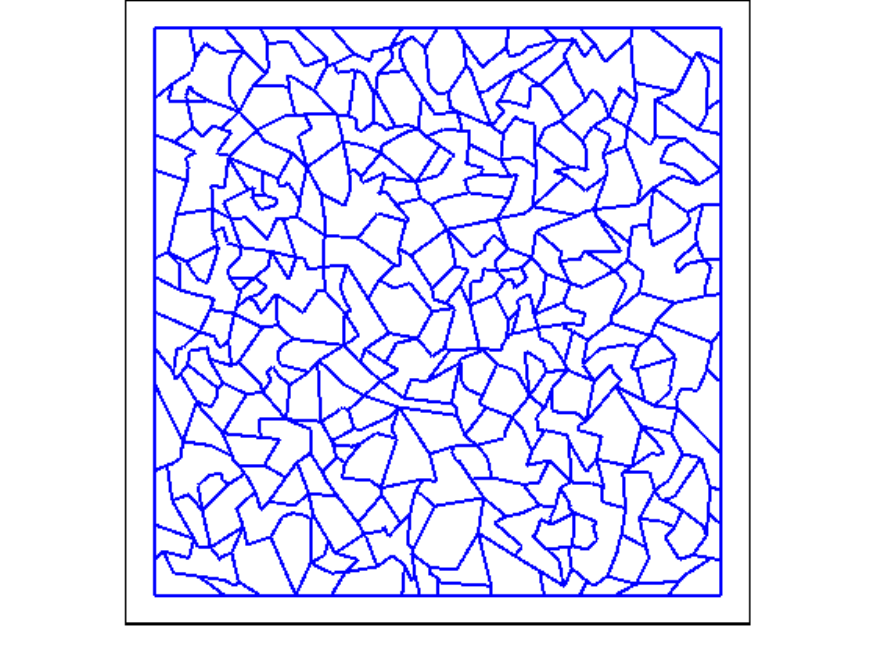}}\quad
		\subfloat[Diamond]{\includegraphics[trim = 40 0 40 0, width=0.23\linewidth]{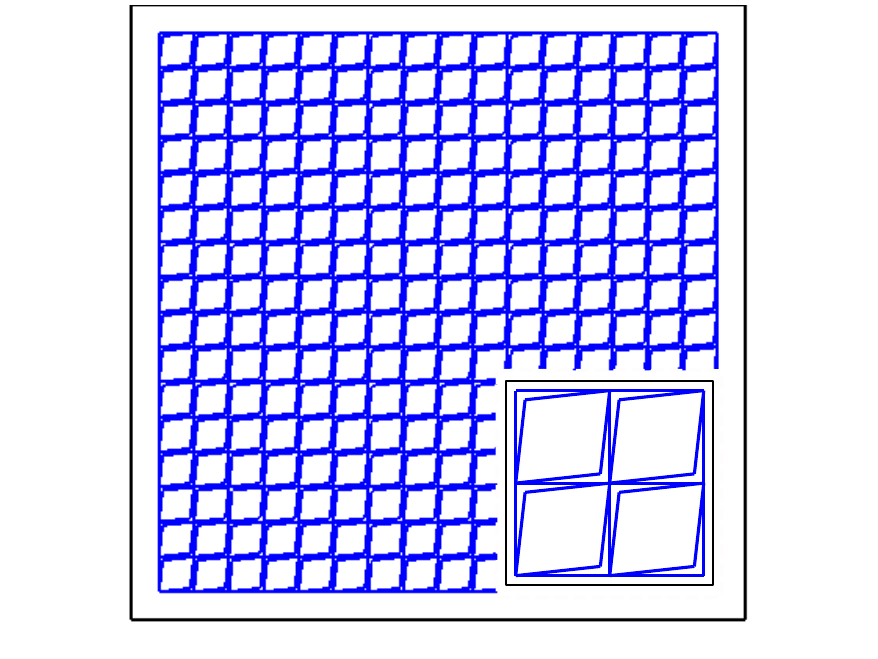}}
		\caption{Meshes used for studying the condition number. Mesh features are collected in Table~\ref{tab:mesh_info}. From left to right: hexagons (Level 3), Voronoi (Level 2), random polygons (Level 1), diamond with zoom (Level 5).}
		\label{fig:meshes}
	\end{figure}

\begin{figure}
		
		\centering
		\textbf{Condition number\\Scaled monomials $\mathcal{M}_{k}$ \textit{vs} $L^2$ orthonormal basis $\mathcal{Q}_{k}$}
		
		\
		
		\includegraphics[width=0.45\linewidth]{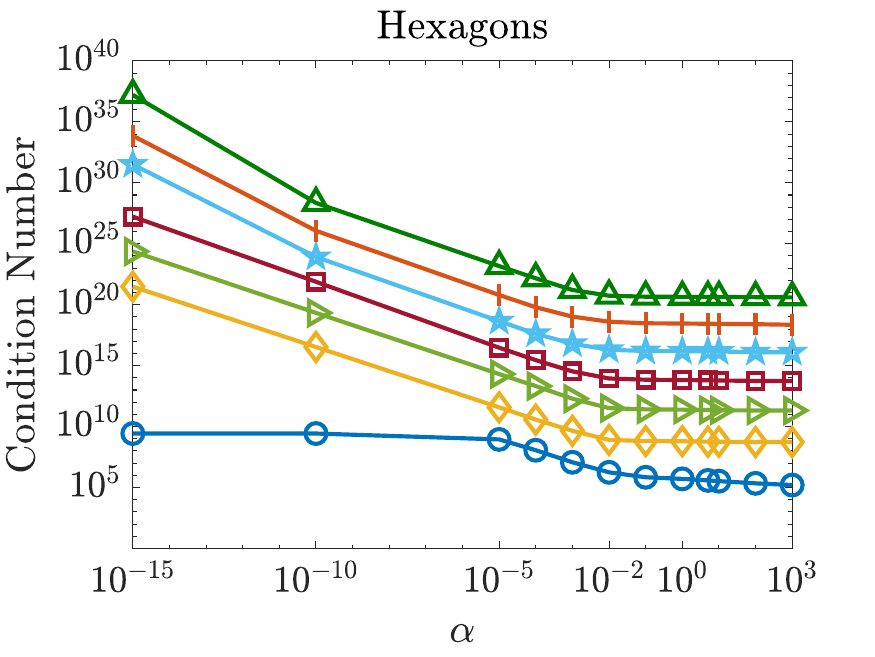}\quad
		\includegraphics[width=0.45\linewidth]{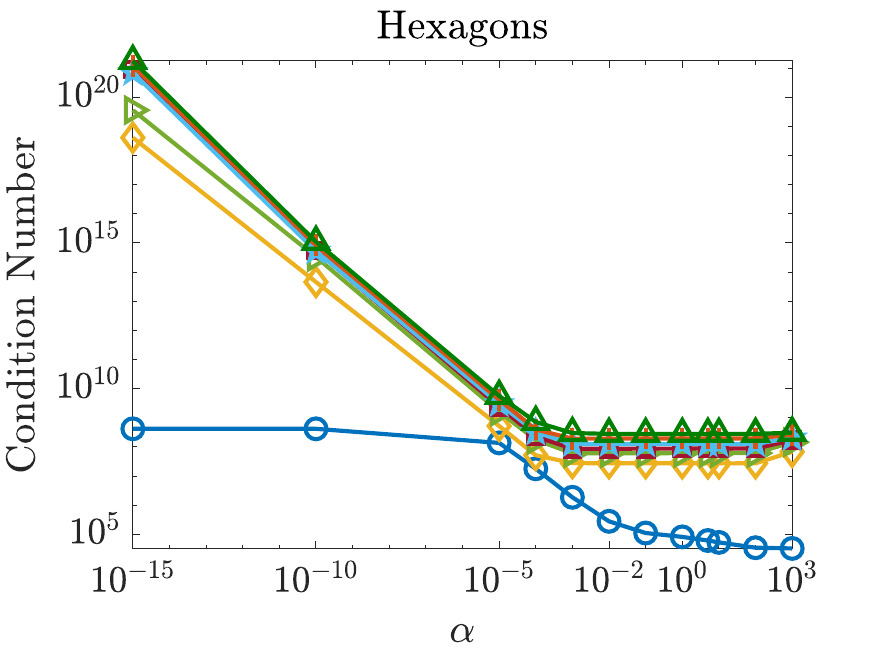}\\\medskip
		\includegraphics[width=0.45\linewidth]{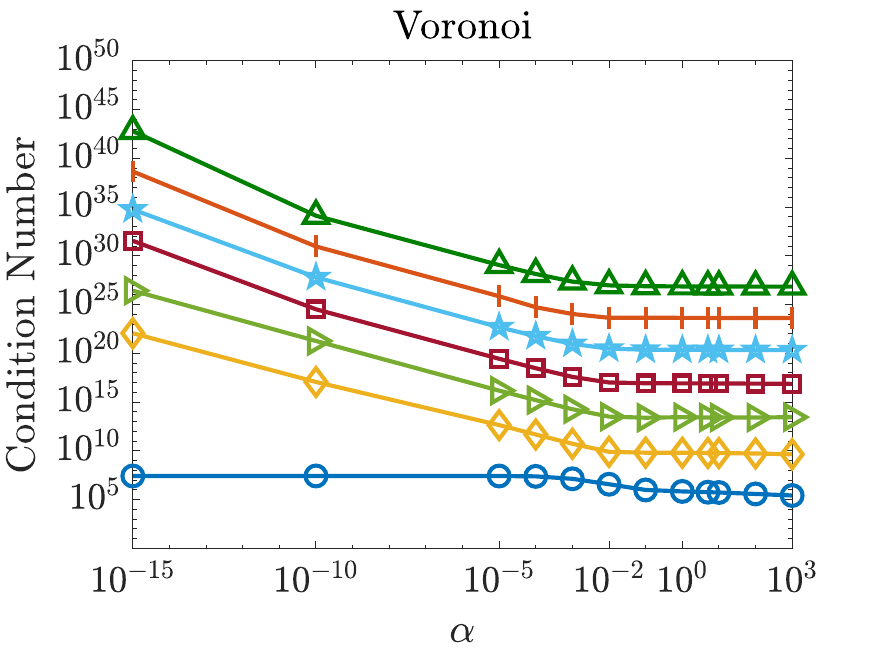}\quad
		\includegraphics[width=0.45\linewidth]{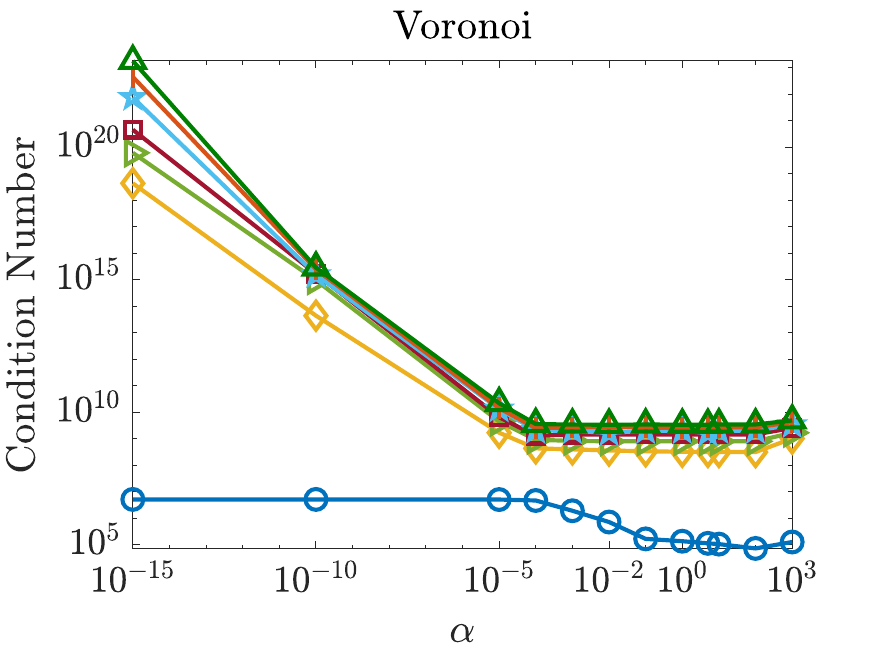}\\\medskip
		\includegraphics[width=0.45\linewidth]{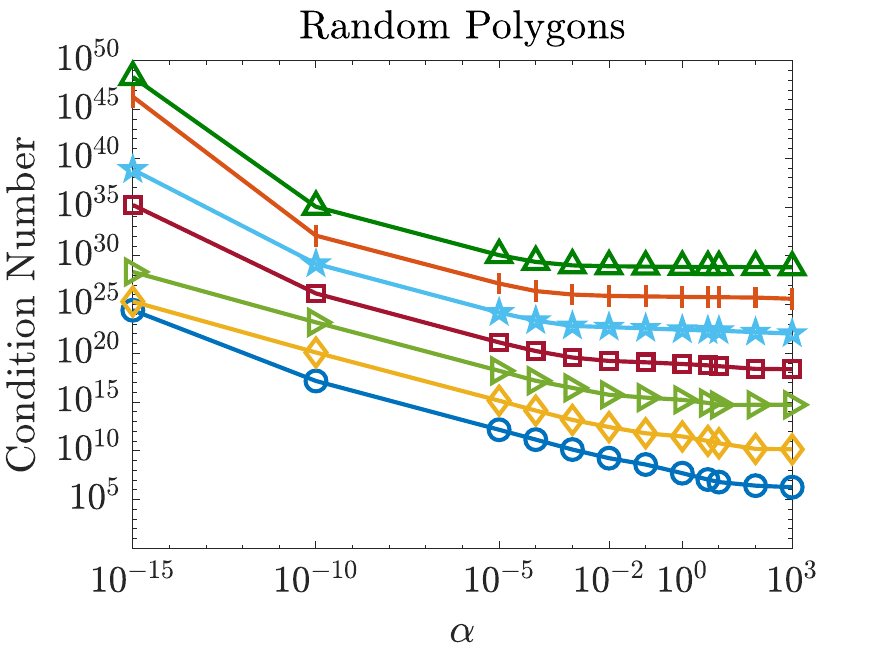}\quad
		\includegraphics[width=0.45\linewidth]{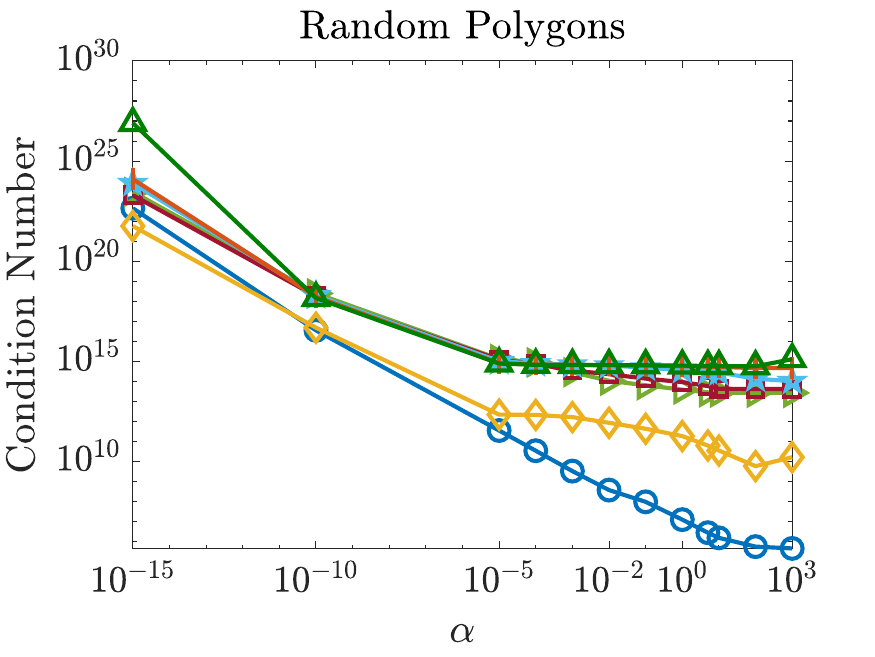}\\\medskip
		\includegraphics[width=0.45\linewidth]{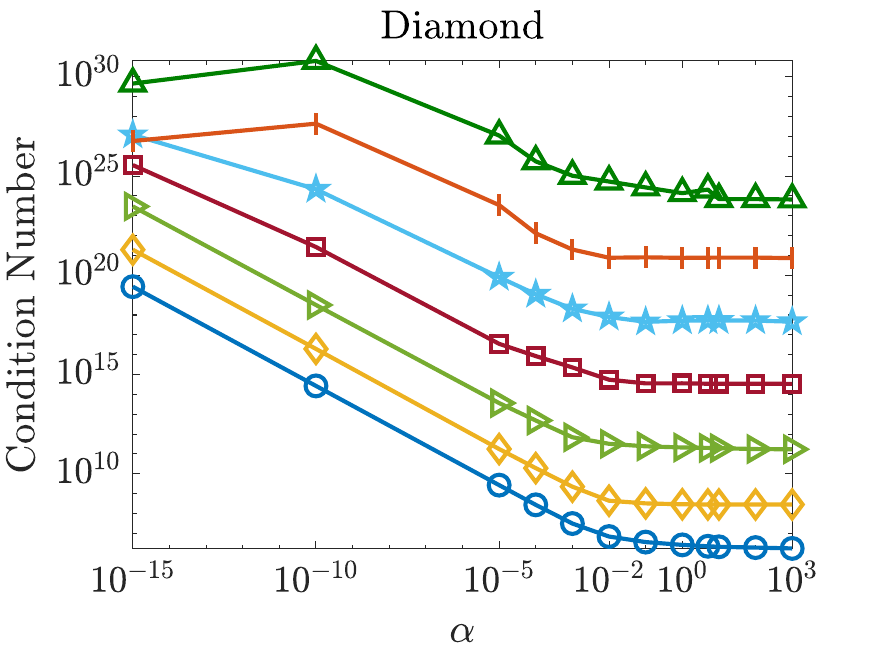}\quad
		\includegraphics[width=0.45\linewidth]{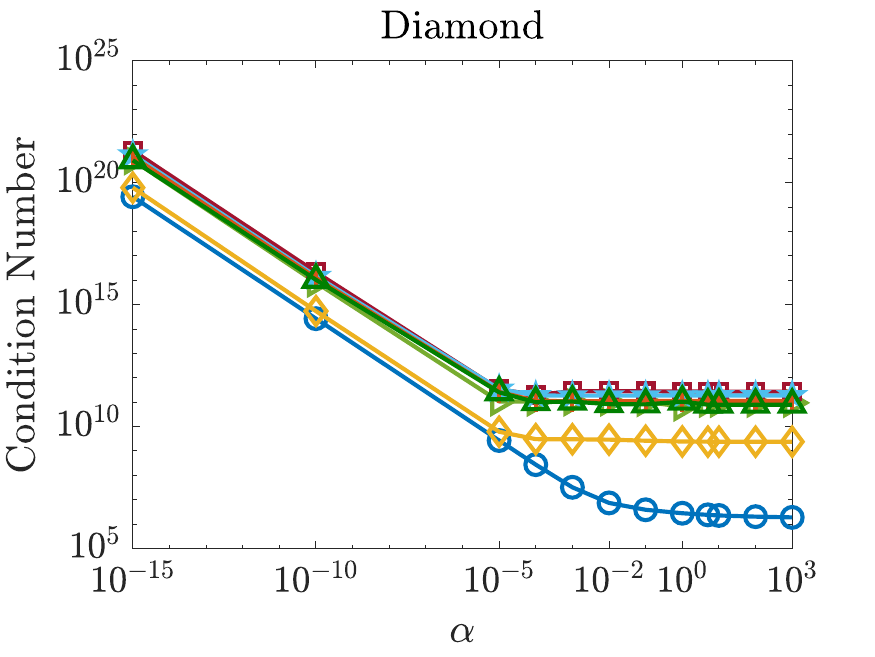}\\
		\includegraphics[width=0.6\linewidth]{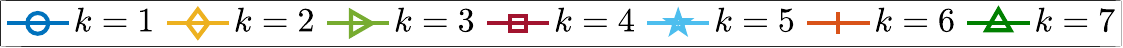}\medskip
		
		\caption{Condition number of MINI--VEM with respect to pressure stabilization parameter $\alpha$. Effect of polynomial basis: scaled monomials $\mathcal{M}_{k}$ (left column) compared with $L^2$ orthonormal basis $\mathcal{Q}_{k}$ (right column).}
		\label{fig:conditioning}
	\end{figure}

	We set $\Omega=(0,1)^2$ and we assemble the linear system up to order \rv{$k=7$} for the meshes depicted in Figure~\ref{fig:meshes}: a mesh of hexagons, a Voronoi mesh, a mesh of random nonconvex polygons, a mesh with diamond-shaped elements. \rv{Notice that the elements of the ``random'' mesh are not star-shaped in general: we consider this example to analyze the robustness of the proposed method when the standard geometrical requirements of virtual elements are relaxed.} The behavior of the condition number for $\mathcal{M}_{k}$ is reported in the left column of Figure~\ref{fig:conditioning}, while the results for $\mathcal{Q}_{k}$ are in the right column of the same figure. Condition number is plotted with respect to the stabilization parameter $\alpha$.
	
	In agreement with the results presented in~\cite{mascotto2018ill}, the system assembled by means of the scaled monomials has a larger condition number than the system constructed by the $L^2$ orthonormal basis. Increasing the polynomial degree, the condition number of $\mathcal{M}_{k}$ quickly increases to prohibitive values. This phenomenon affects also the behavior of the conditioning with respect to the stabilizing parameter $\alpha$. In the case of assembly with scaled monomials, the condition number reaches the steady state when $\alpha\ge10^{-2}$. On the other hand, if the $L^2$ orthonormal basis $\mathcal{Q}_{k}$ is considered, the condition number stabilizes for $\alpha\ge10^{-5}$. Moreover, especially for the random polygons mesh, the condition number of the lowest order method becomes stable for larger values of $\alpha$ than the higher order approximations. 
	\rv{We finally observe that, for hexagonal and Voronoi meshes, the condition number for lowest-order method does not appear to fully blow up as $\alpha$ goes to $0$. Indeed, unlike what happens in all the other tests, the stiffness matrix for $\alpha = 0$ appears, for the specific meshes considered, to be invertible (though extremely ill conditioned: for an order one method with a mesh size of the order $10^{-1}$, a condition number of the order $10^9$ is extremely high).}
	
	\subsection{Convergence tests}
	
	In this subsection, we study the convergence of the MINI--VEM with respect to mesh refinement. We consider four sequences of meshes of the kind depicted in Figure~\ref{fig:meshes}: the information regarding number of vertices, edges and elements for each level of refinement are reported in Table~\ref{tab:mesh_info}. We solve again the Stokes problem in the unit square with polynomial accuracy up to order \rv{$k=6$}. All the tests in this subsection are performed by setting~$\alpha=1$ and by considering the polynomial basis $\mathcal{Q}_{k}$.
	
	\renewcommand{\arraystretch}{1.15}
	\begin{table}
		\begin{tabular}{ccccc|ccccc}
			\multicolumn{10}{c}{\textbf{MESH DATA}}\\
			\hline
			\multicolumn{10}{c}{$N_V$ = number vertices, $N_e$ = number edges, $N_\elementvem$ = number elements}\\
			\hline
			\multicolumn{5}{c|}{\textbf{Hexagons}}&\multicolumn{5}{c}{\textbf{Voronoi}}\\
			\hline
			Level	&	$N_V$	&	$N_e$	&	$N_\elementvem$	& $h$ &Level	&	$N_V$	&	$N_e$	&	$N_\elementvem$	& $h$\\
			1		&	$62$	&	$91$	&	$30	$			& $0.290$ &	1 &	$124$	&	$187$	&	$64$			&	$0.0375$\\
			2		&	$242$	&	$361$	&	$120$			& $0.144$ &	2	&	$458$	&	$713$	&	$256$			&	$0.189$\\
			3		&	$542$	&	$811$	&	$270$			& $0.096$ &	3	&	$1832$  &	$2855$	&	$1024$			&	$0.049$\\
			4		&	$922$	&	$1381$	&	$460$			& $0.082$ &	4	&	$7428$  &	$11523$	&	$4096$			&	$0.038$\\
			5		&	$3682$	&	$5521$	&	$1840$			& $0.043$ &	--	&	--		&	--		&	--				&	--\\
			6		&	$14882$	&	$22321$	&	$7440$			& $0.018$&	--	&	--		&	--		&	--				&	--\\
			\hline
			\multicolumn{5}{c|}{\textbf{Random polygons}}&\multicolumn{5}{c}{\textbf{Diamond}}\\
			\hline
			Level	&	$N_V$	&	$N_e$	&	$N_\elementvem$	& $h$ 		&Level	&	$N_V$	&	$N_e$	&	$N_\elementvem$	& $h$\\
			1		&	$256$	&	$331$	&	$64$			&	$0.357$	&	1	&	$57$	&	$104$	&	$48$			&	$0.707$\\
			2		&	$501$	&	$656$	&	$128$			&	$0.295$	&	2	&	$86$	&	$160$	&	$75$			&	$0.354$\\
			3		&	$908$	&	$1196$	&	$256$			&	$0.192$	&	3	&	$209$	&	$400$	&	$192$			&	$0.283$\\
			4		&	$1747$	&	$2327$	&	$512$			&	$0.147$	&	4	&	$321$	&	$620$	&	$300$			&	$0.177$\\
			5		&	$3411$	&	$4566$	&	$1024$			&	$0.104$	&	5	&	$801$	&	$1568$	&	$768$			&	$0.141$\\
			6		&	$13350$	&	$17985$	&	$4096$			&	$0.053$	&	6	&	$1241$	&	$2440$	&	$1200$			&	$0.088$\\
			--		&	--		&	--		&	--				&	--		&	7	&	$3137$	&	$6208$	&	$3072$			&	$0.044$\\
			\hline
		\end{tabular}
		\caption{Mesh data for each sequence used in the numerical investigation.}
		\label{tab:mesh_info}
	\end{table}
	\renewcommand{\arraystretch}{1}
	
	We consider two different analytical solutions (chosen from~\cite{cioncolini2019mini}) and we compute the right hand side accordingly.
	
	\
	
	\textbf{Test 1} 
	\begin{gather*}
		u_1(x,y) = \sin(2\pi y)(1-\cos(2 \pi x)),\qquad\qquad
		u_2(x,y) = \sin(2\pi x)(\cos(2 \pi y)-1),\\
		p(x,y) = 2\pi(\cos(2\pi y)-\cos(2\pi x)).\\
	\end{gather*}
	
	\
	
	\textbf{Test 2} 
	\begin{gather*}
		u_1(x,y) = (x^4-2x^3+x^2)(2y^3-y),\qquad\qquad
		u_2(x,y) = -(2x^3-3x^2+x)(y^4-y^2),\\
		p(x,y) = (4x^3-6x^2+2x)(2y^3-y)+\frac15(6x^5-15x^4+10x^3)y -\frac1{10}.\\
	\end{gather*}
	
	The velocity in Test~1 has homogeneous Dirichlet boundary conditions on the entire $\partial\Omega$, whereas the velocity in Test~2 has nonzero Dirichlet boundary conditions on the top edge of the domain. 
	
\begin{figure}
		\centering
		\subfloat[Hexagons]{
			\includegraphics[width=0.35\linewidth]{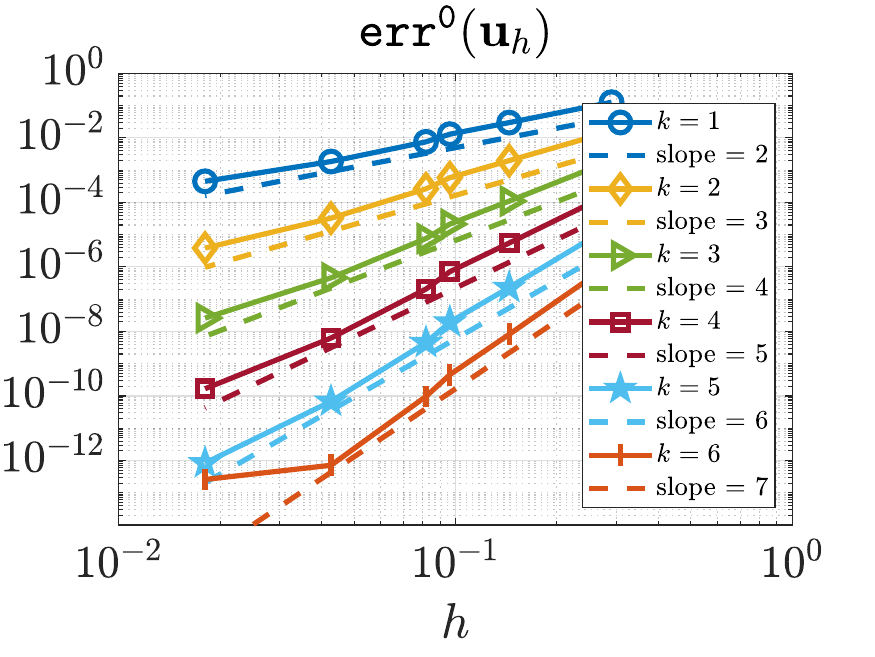}
			\includegraphics[width=0.35\linewidth]{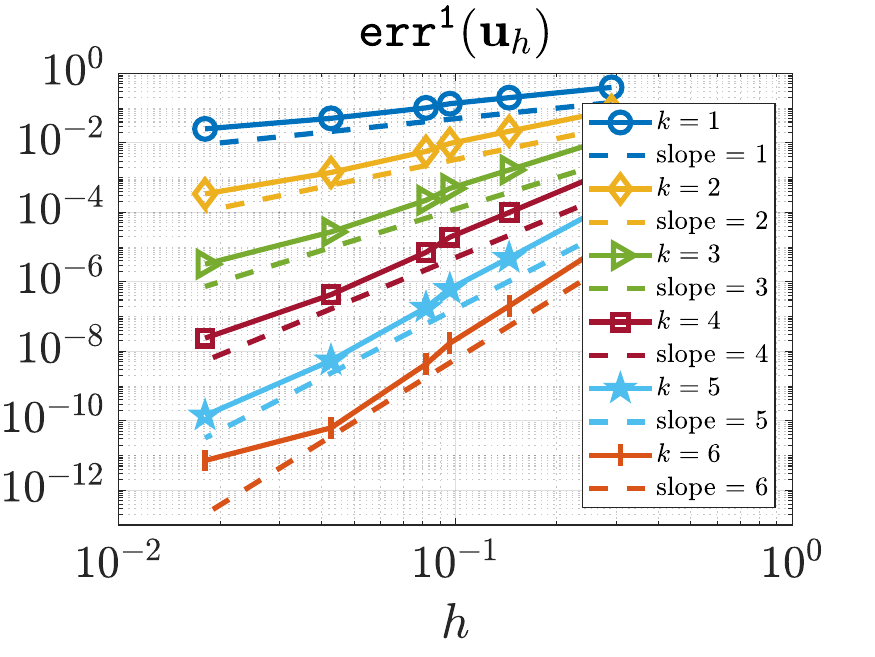}
			\includegraphics[width=0.35\linewidth]{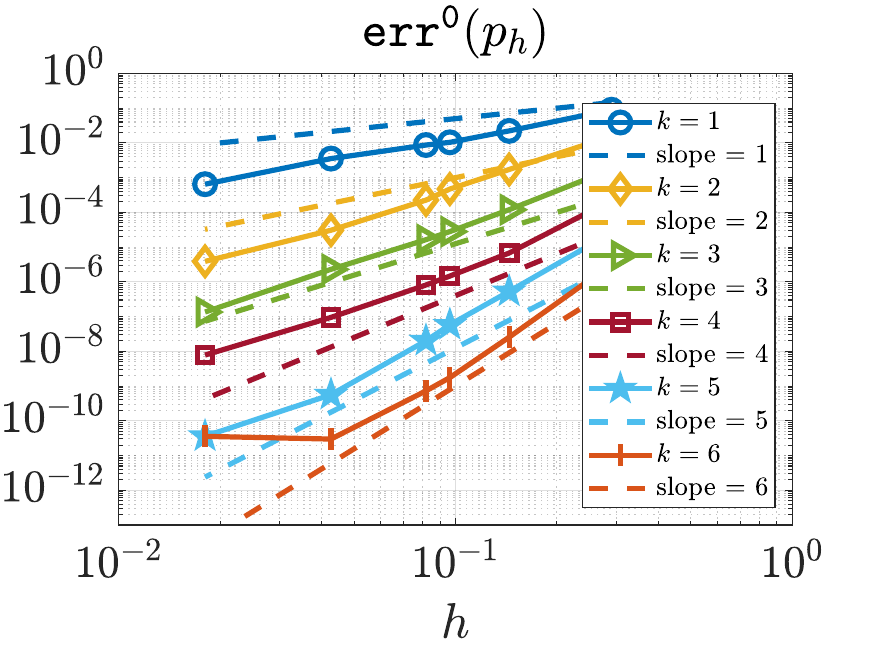}
		}
		
		\
		
		\subfloat[Voronoi]{
			\includegraphics[width=0.35\linewidth]{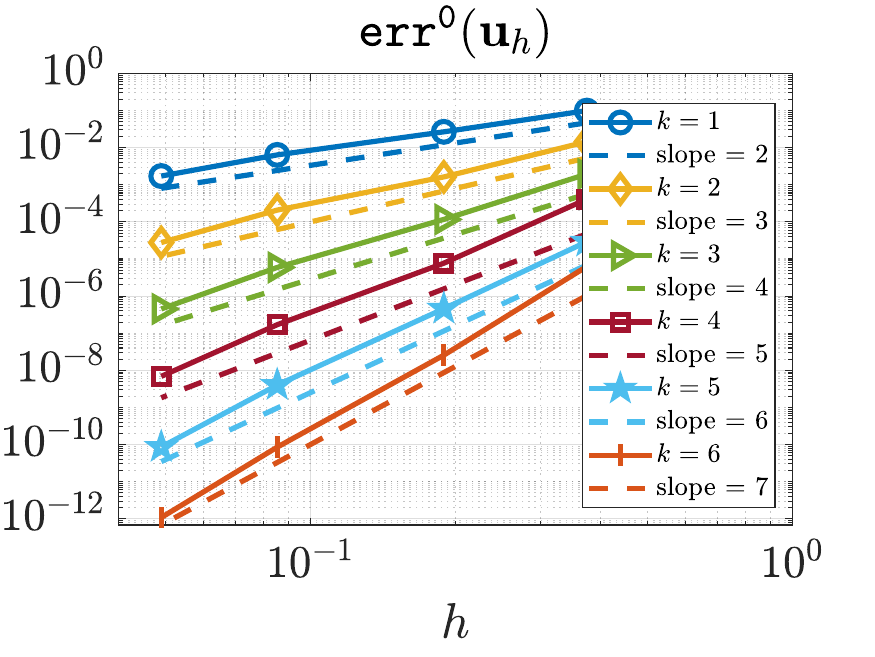}
			\includegraphics[width=0.35\linewidth]{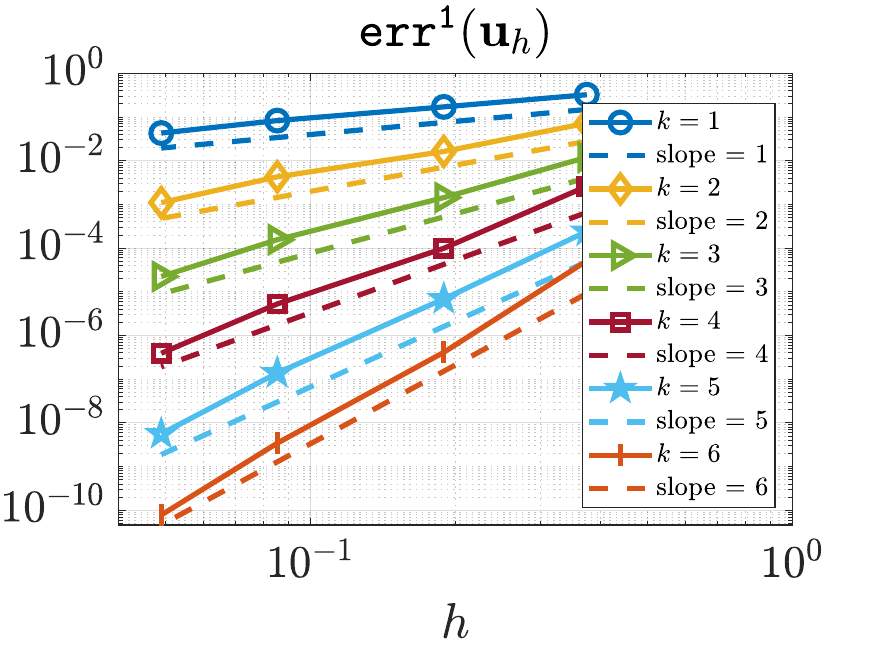}
			\includegraphics[width=0.35\linewidth]{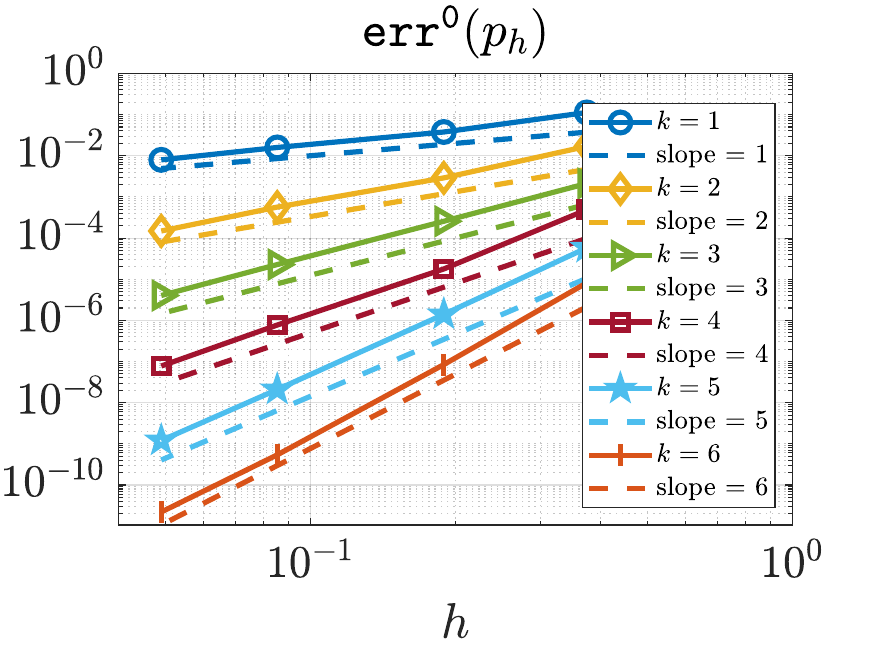}
		}
		
		\
		
		\subfloat[Random polygons]{
			\includegraphics[width=0.35\linewidth]{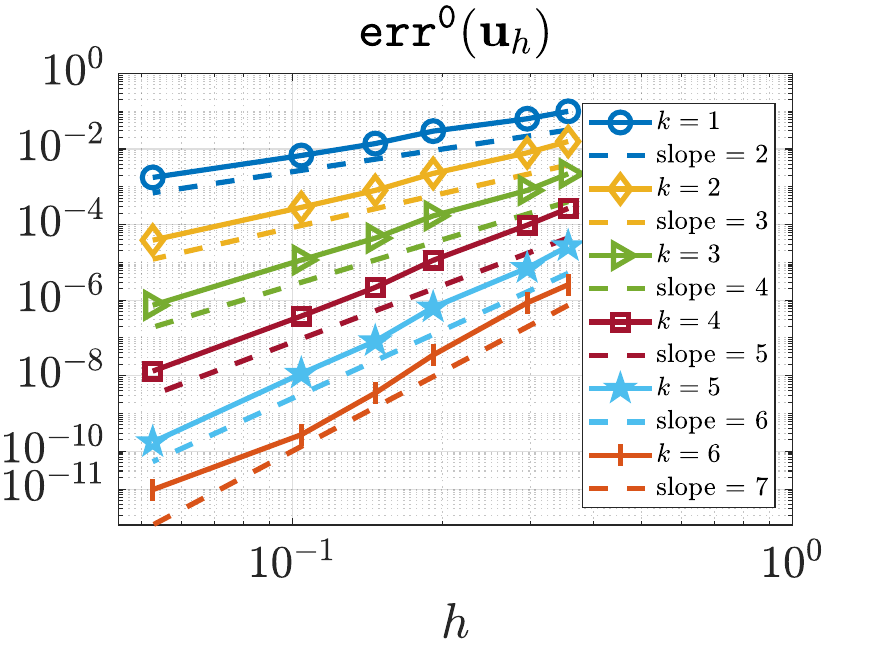}
			\includegraphics[width=0.35\linewidth]{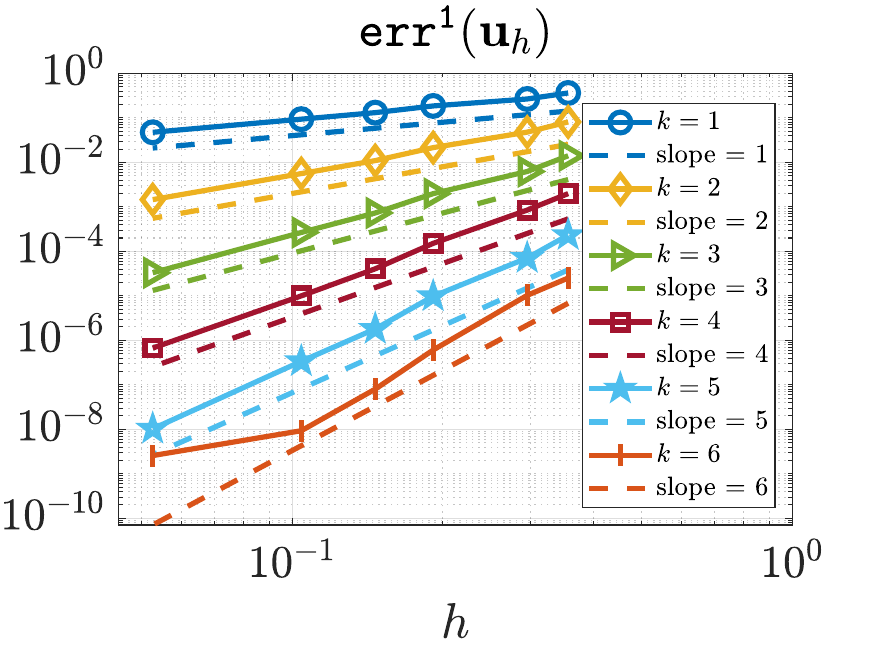}
			\includegraphics[width=0.35\linewidth]{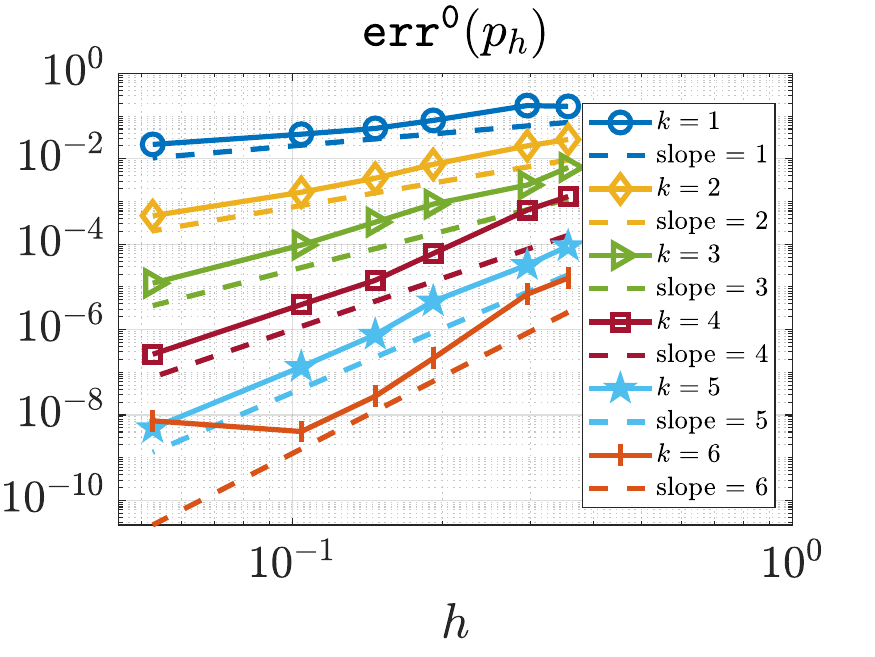}
		}
		
		\
		
		\subfloat[Diamond]{
			\includegraphics[width=0.35\linewidth]{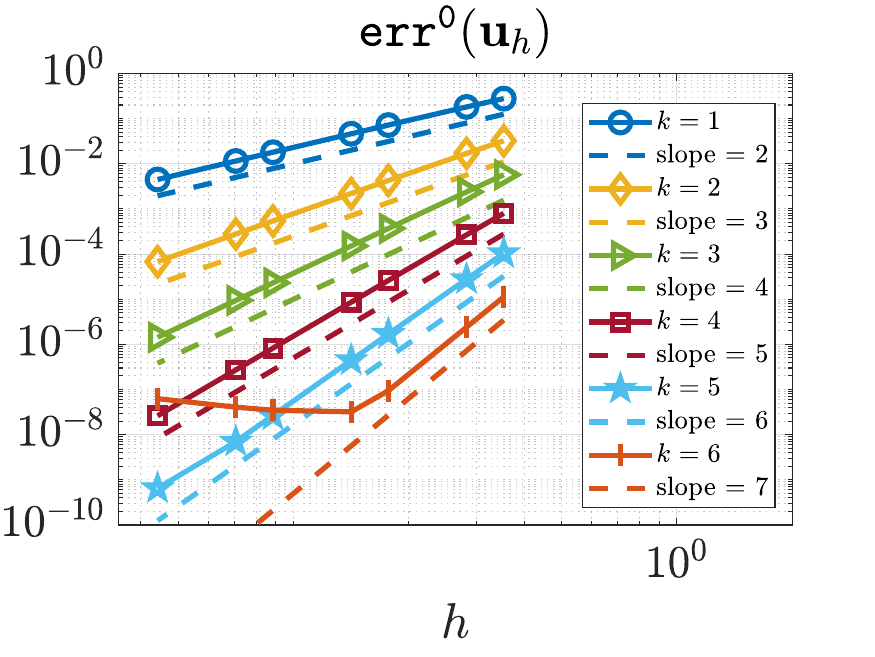}
			\includegraphics[width=0.35\linewidth]{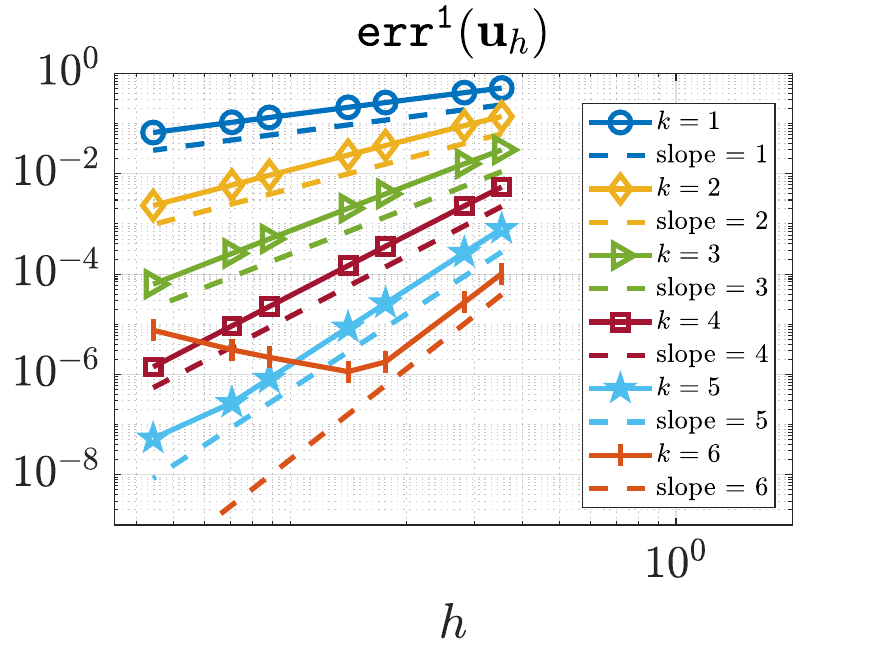}
			\includegraphics[width=0.35\linewidth]{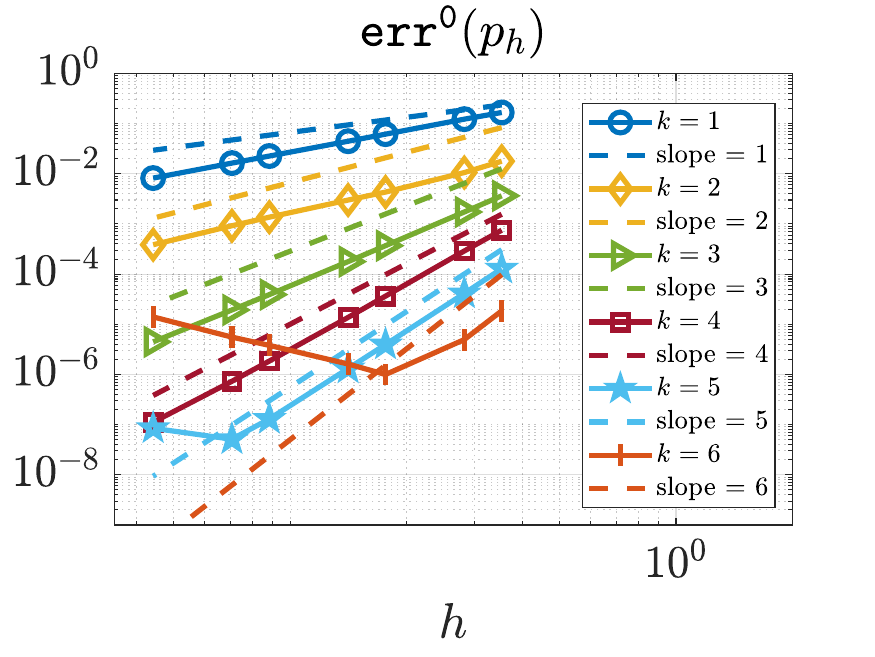}
		}
		
		\caption{Convergence analysis for Test 1. \rv{The errors, defined in~\eqref{eq:errors},} are represented with solid lines, whereas dashed lines are employed to draw the theoretical slopes.}
		\label{fig:test1}
	\end{figure}

\begin{figure}
		
		\subfloat[Hexagons]{
			\includegraphics[width=0.35\linewidth]{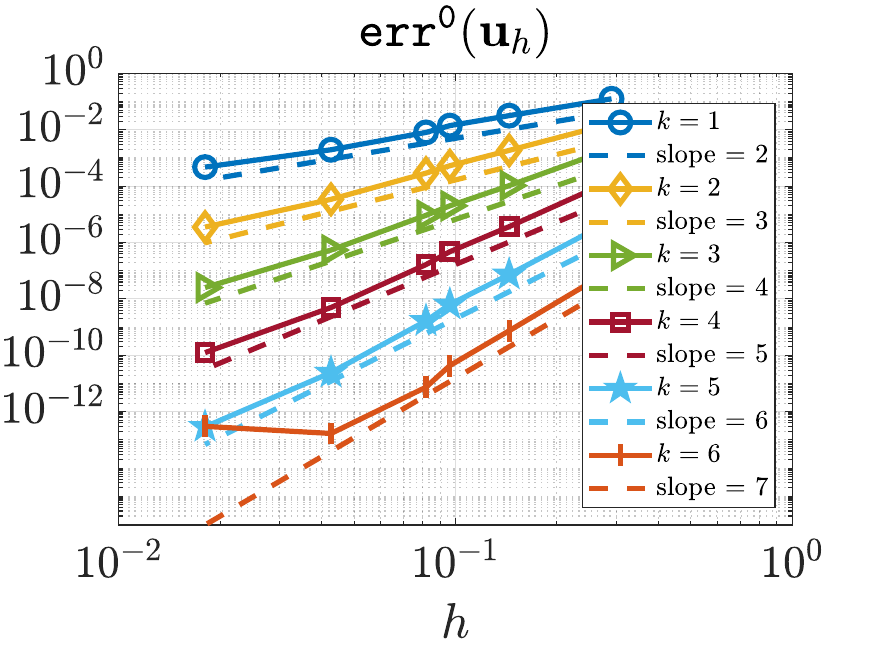}
			\includegraphics[width=0.35\linewidth]{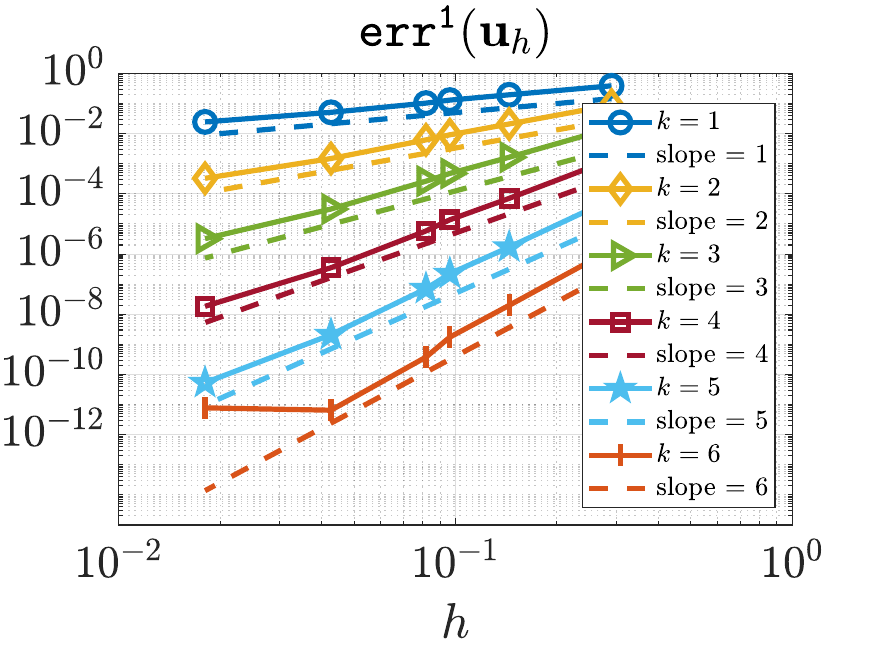}
			\includegraphics[width=0.35\linewidth]{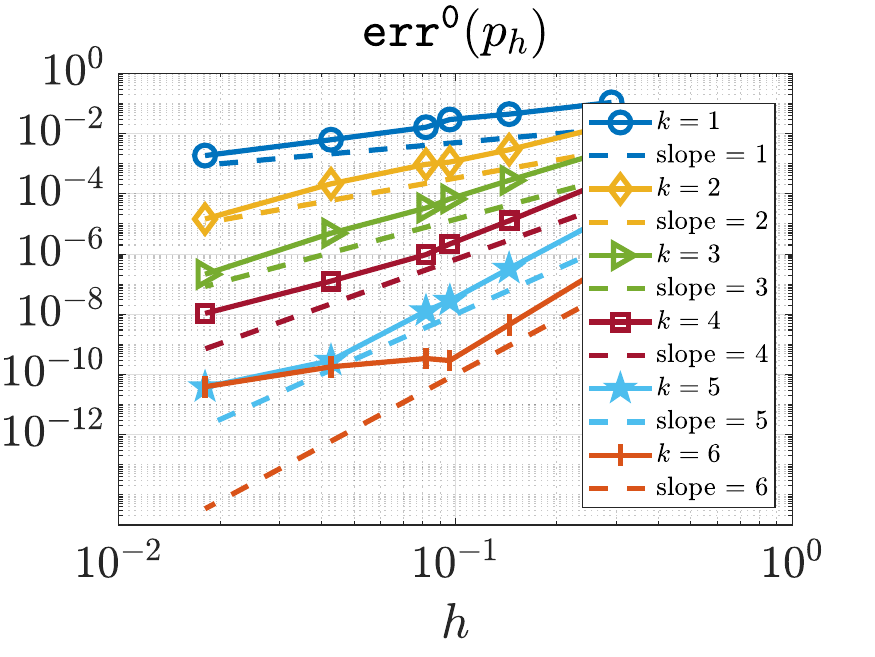}
		}
		
		\
		
		\subfloat[Voronoi]{
			\includegraphics[width=0.35\linewidth]{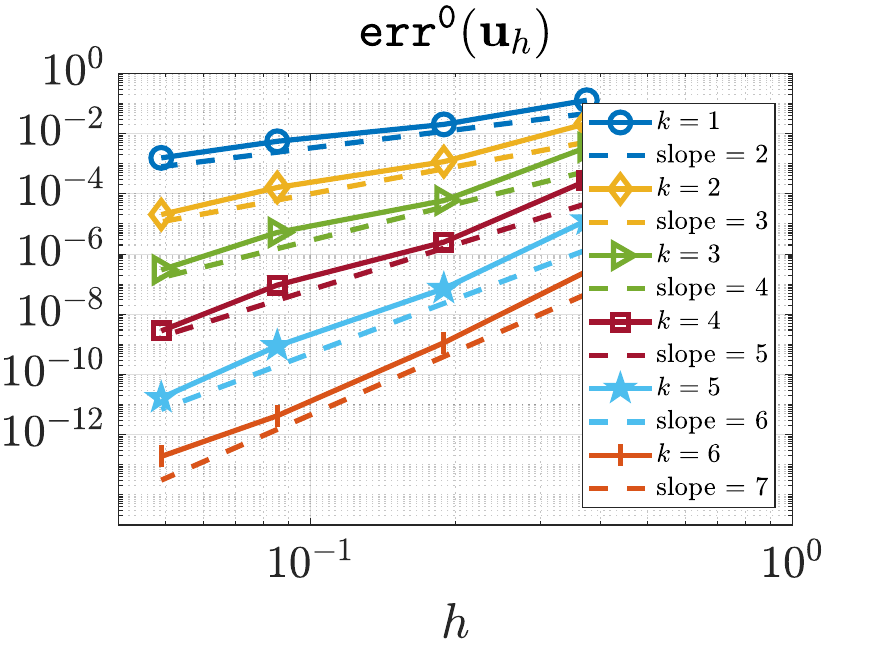}
			\includegraphics[width=0.35\linewidth]{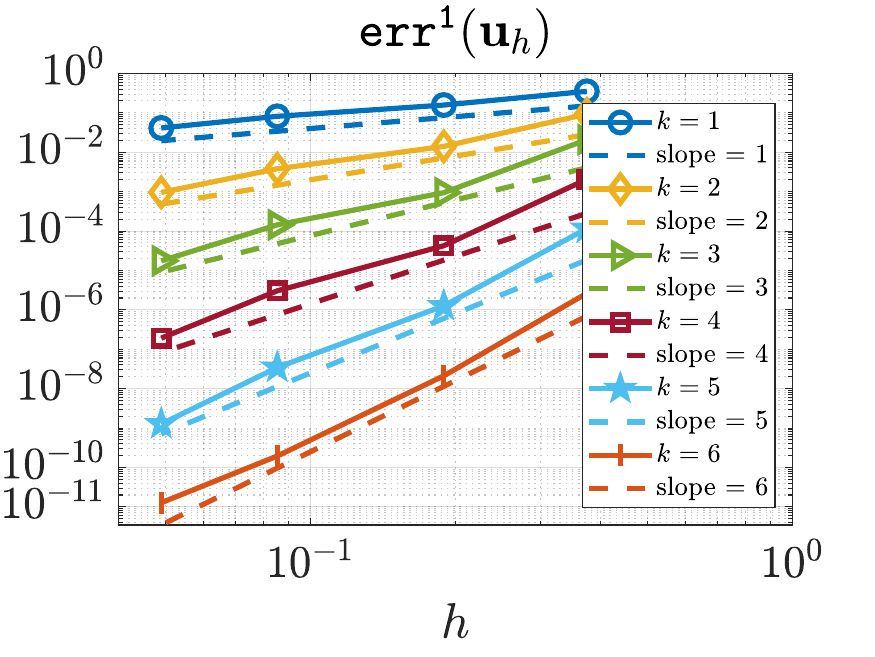}
			\includegraphics[width=0.35\linewidth]{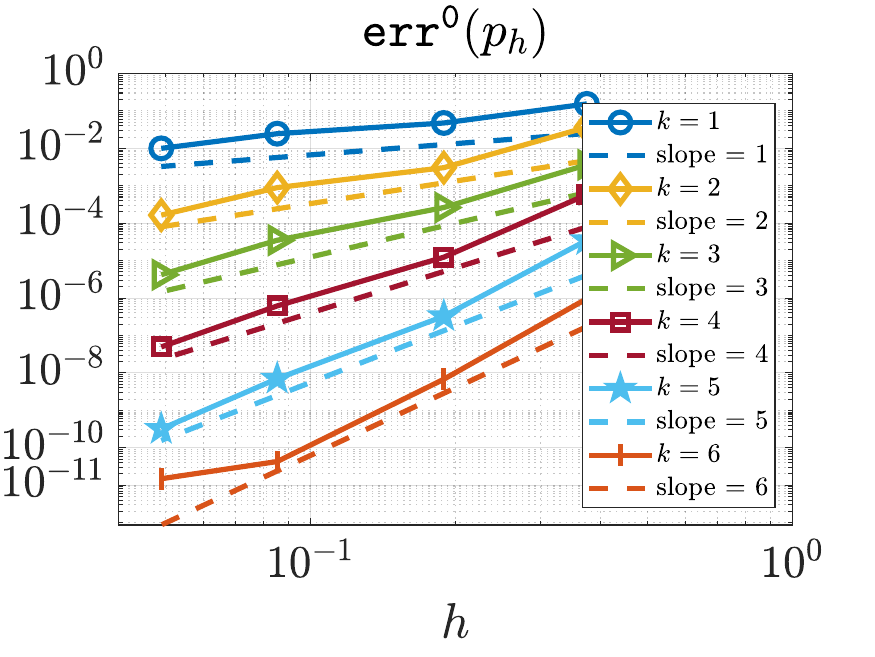}
		}
		
		\
		
		\subfloat[Random polygons]{
			\includegraphics[width=0.35\linewidth]{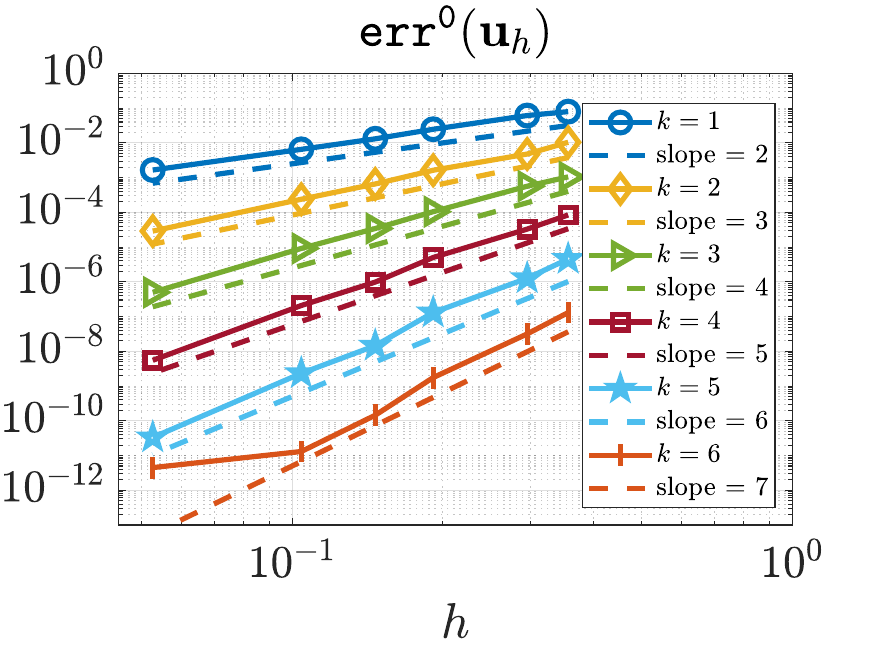}
			\includegraphics[width=0.35\linewidth]{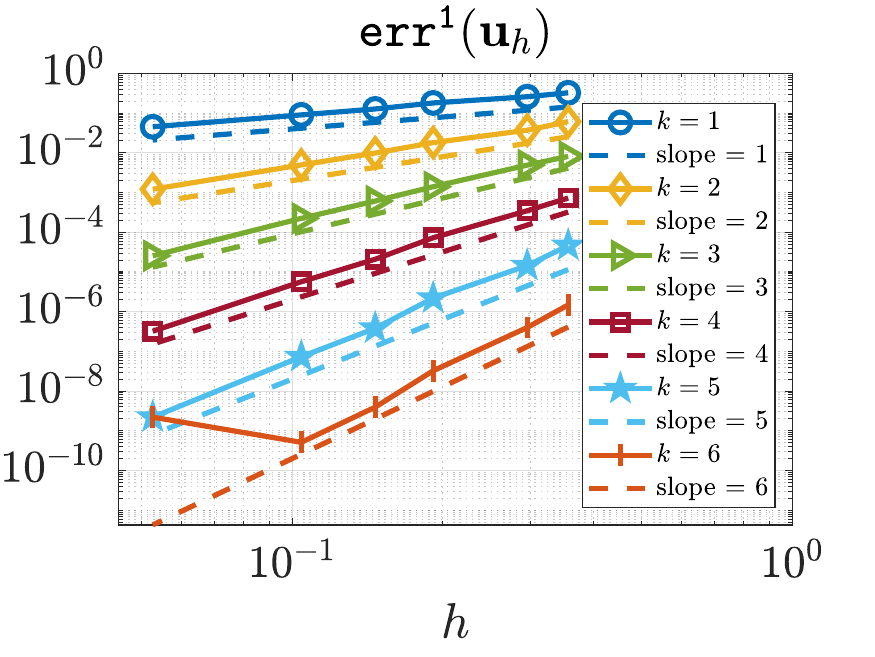}
			\includegraphics[width=0.35\linewidth]{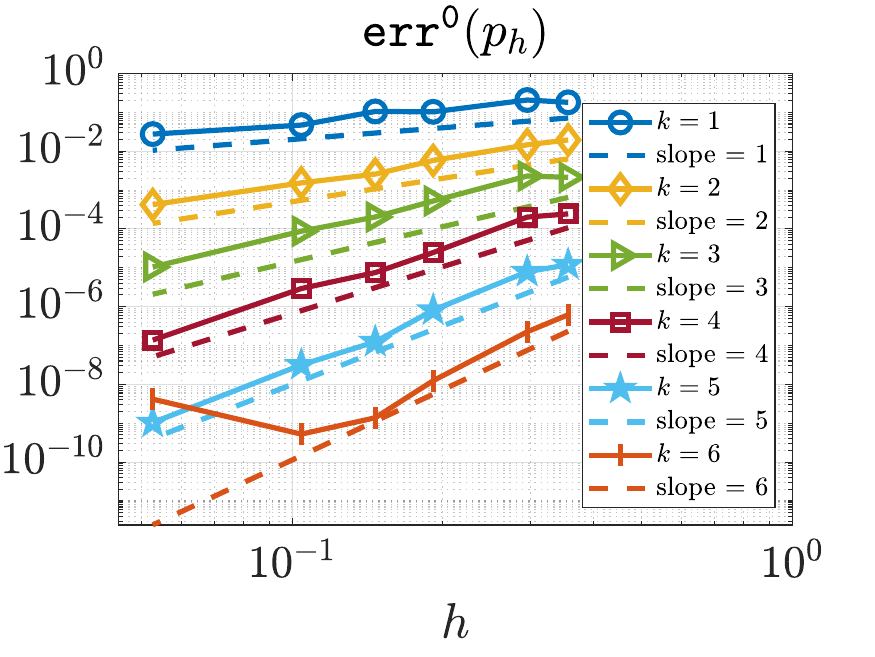}
		}
		
		\
		
		\subfloat[Diamond]{
			\includegraphics[width=0.35\linewidth]{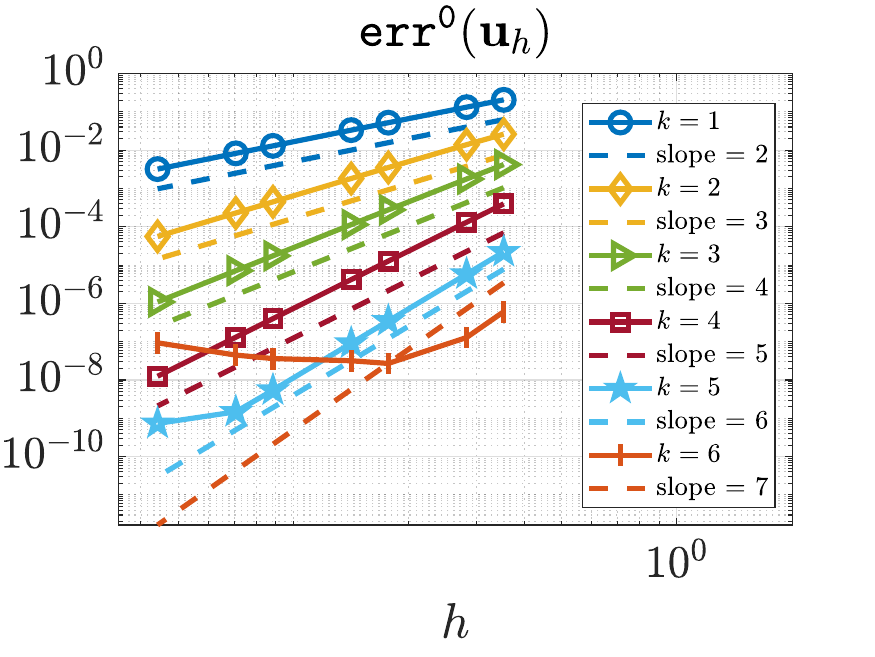}
			\includegraphics[width=0.35\linewidth]{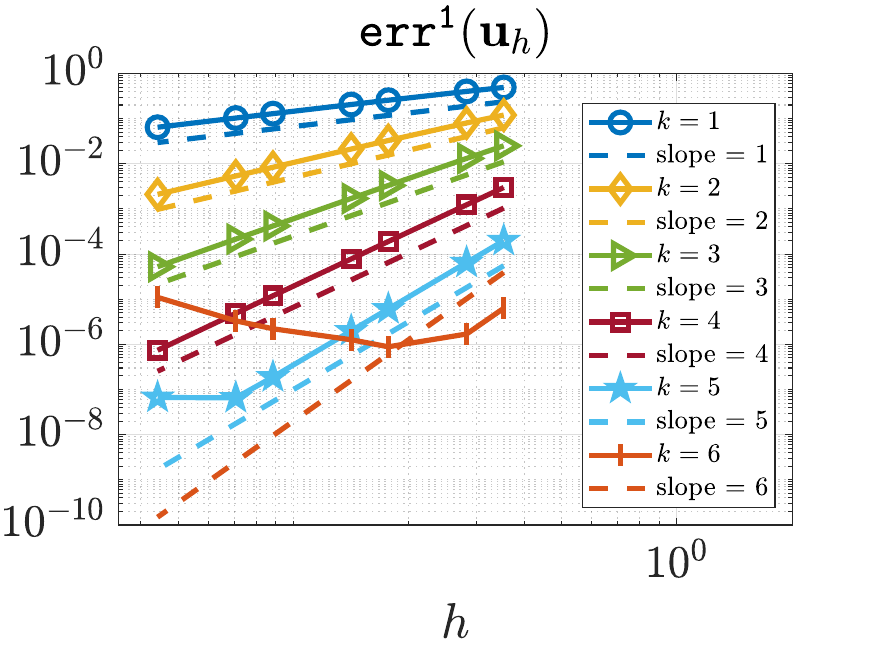}
			\includegraphics[width=0.35\linewidth]{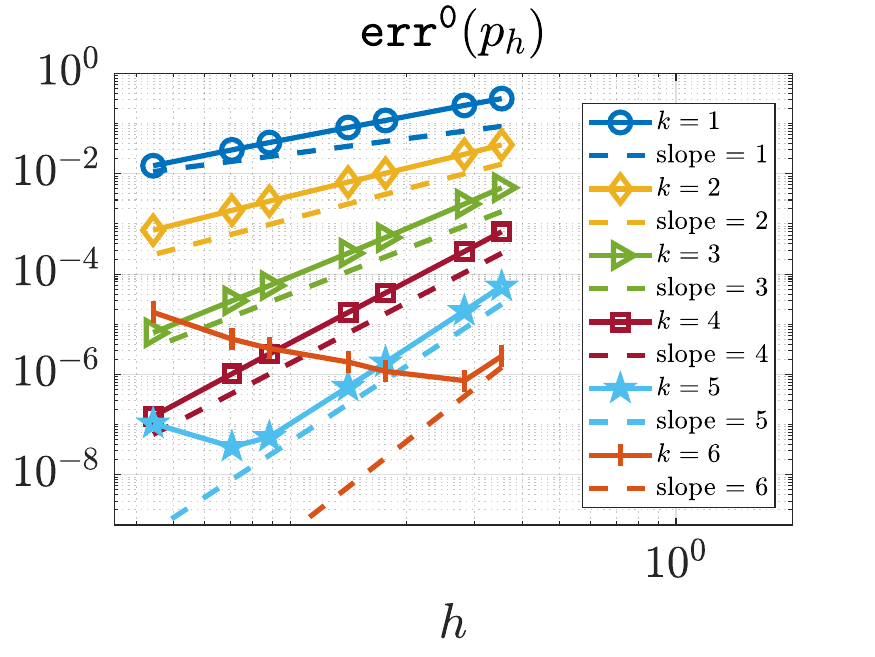}
		}
		
		\caption{Convergence analysis for Test 2. \rv{The errors, defined in~\eqref{eq:errors},} are represented with solid lines, whereas dashed lines are employed to draw the theoretical slopes.}
		\label{fig:test2}
	\end{figure}

	The results for Test~1 are collected in Figure~\ref{fig:test1}, while the results for Test~2 can be found in Figure~\ref{fig:test2}. Error curves are represented by solid line, while dashed lines are employed to represent the theoretical slopes. From the convergence plots, it is evident that the estimates presented in Section~\ref{sec:error_analysis} are confirmed by all the numerical tests and satisfied for all values of polynomial accuracy. \rv{We point out that high order approximations on the random and diamond meshes suffer from bad conditioning since elements are not shape regular.}
	
	\rv{
		
		\subsection{Comparison with the MINI finite element}\label{sec:tests_fem_vem}
		
		When $k=1$, the local discrete spaces are defined as
		$$
		\V_1(\elementvem) = [\enhanced{1}(\elementvem)]^2\oplus[\rvf{\bubble{2}}]^2\qquad\text{and}\qquad Q_1(\elementvem)=\enhanced{1}(\elementvem).
		$$
		More precisely, the degrees of freedom of $\enhanced{1}(\elementvem)$ are just the values at the vertices of $\elementvem$, while the bubble space consists of virtual \rvf{quadratic} bubbles described by their moments \rvf{of} order \rvf{zero}.
		
		In this framework, we find the virtual element version of the MINI mixed finite element. Indeed, if we consider the special case where $\elementvem$ is a triangle, the space $\enhanced{1}(\elementvem)$ coincides with the space of linear polynomials in $\elementvem$, while the bubble space $\bubble{2}$ is identified by just one internal degree of freedom. The only difference between this particular case of MINI--VEM and the original MINI--FEM is given by the polynomial order of the bubble function: in the former, we have virtual quadratic bubbles, in the latter, cubic bubbles given by the product of barycentric coordinates. Moreover, since the pressure space is purely polynomial, the bilinear form $b$ is computed exactly and $c_h$ vanishes. Similarly,  for $a_h$, we have
		$$
		\begin{aligned}
			a_h^\elementvem(\u_h,\v_h) & =  a^\elementvem(\Pinabla{1}\utilde,\Pinabla{1}\vtilde) + a^\elementvem(\Pinabla{2}\bolla_h,\Pinabla{2}\bollatest_h)+\beta_\sharp\,S_2^\elementvem(\Qnabla{2}\bolla_h,\Qnabla{2}\bollatest_h)\\
			& =a^\elementvem(\utilde,\vtilde) + a^\elementvem(\Pinabla{2}\bolla_h,\Pinabla{2}\bollatest_h)+\beta_\sharp\,S_2^\elementvem(\Qnabla{2}\bolla_h,\Qnabla{2}\bollatest_h).
		\end{aligned}
		$$
		
		We now briefly compare the performance of the lowest order MINI--VEM on triangular meshes with the original MINI--FEM. We thus consider a sequence of uniform triangulations of size $n\times n$ with ${n=5,10,20,30,40,80}$.
		
		We first analyze the condition number by varying the bubble stabilization parameter $\beta_\sharp$. The results are collected in Figure~\ref{fig:cond_fem_vem}. It is clear that conditioning increases fast when $\beta_\sharp>1$. On the other hand, for $\beta_\sharp\le1$, MINI--VEM has better conditioning than its finite element counterpart, especially when the polynomial basis $\mathcal{Q}_k$ is chosen: in this case the difference is of one order of magnitude.
		
		\begin{figure}
			\centering
			\includegraphics[width=0.45\linewidth]{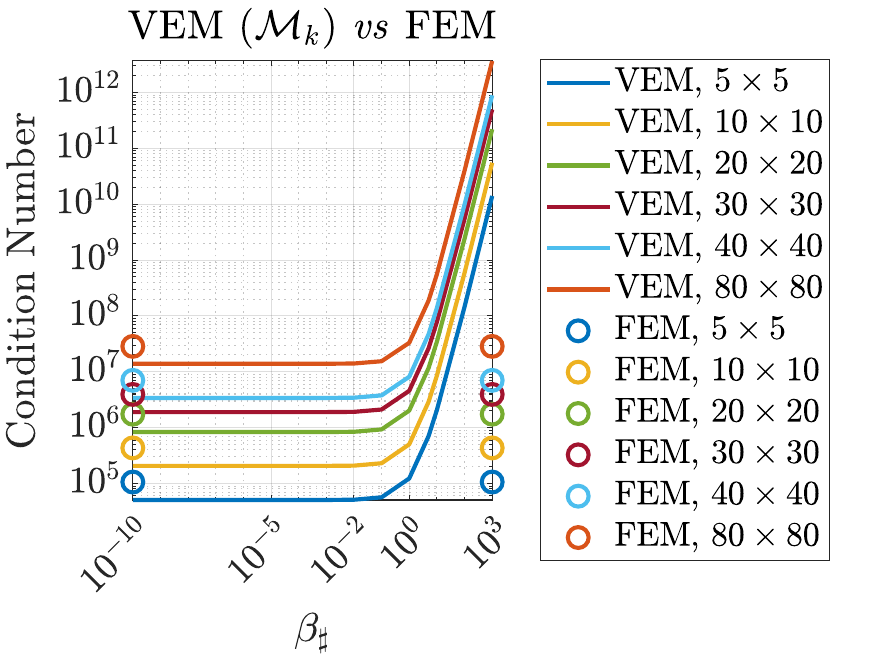}
			\quad
			\includegraphics[width=0.45\linewidth]{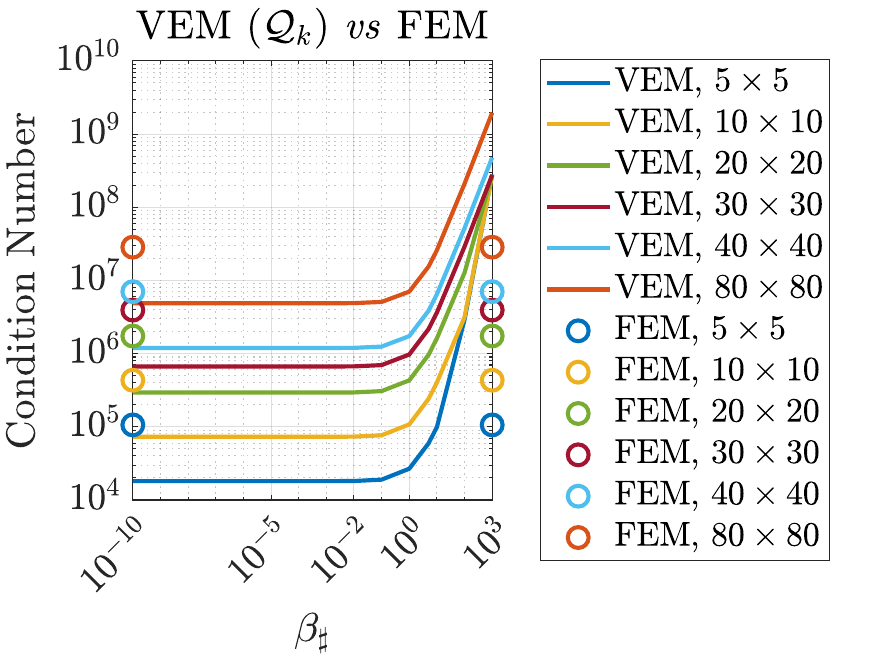}
			\caption{\rv{Comparison of condition number between MINI--VEM and MINI--FEM. Solid lines represent conditioning of MINI--VEM with respect to $\beta_\sharp$, while circles denotes the condition number of MINI--FEM (we placed them at both far left and far right of the picture to easy of readability).}}
			\label{fig:cond_fem_vem}
		\end{figure}
		
		We then compare the considered methods in terms of convergence. We consider Test~1 and we solve it for selected values of $\beta_\sharp$. Convergence plots are reported in Figure~\ref{fig:fem_err} for both polynomial basis and do not consider the bubble contribution to the error. The performance of MINI--VEM and MINI--FEM is basically equivalent. We point out that the pressure error super-converges as extensively studied in~\cite{cioncolini2019mini} in case of finite elements. 
		
		\begin{figure}
			\centering
			\includegraphics[width=0.32\linewidth]{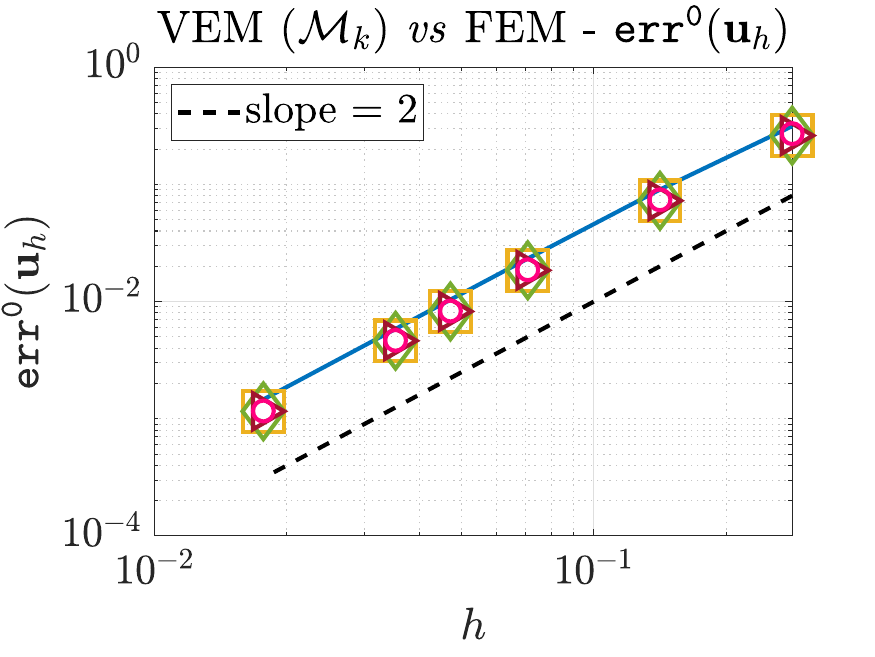}
			\includegraphics[width=0.32\linewidth]{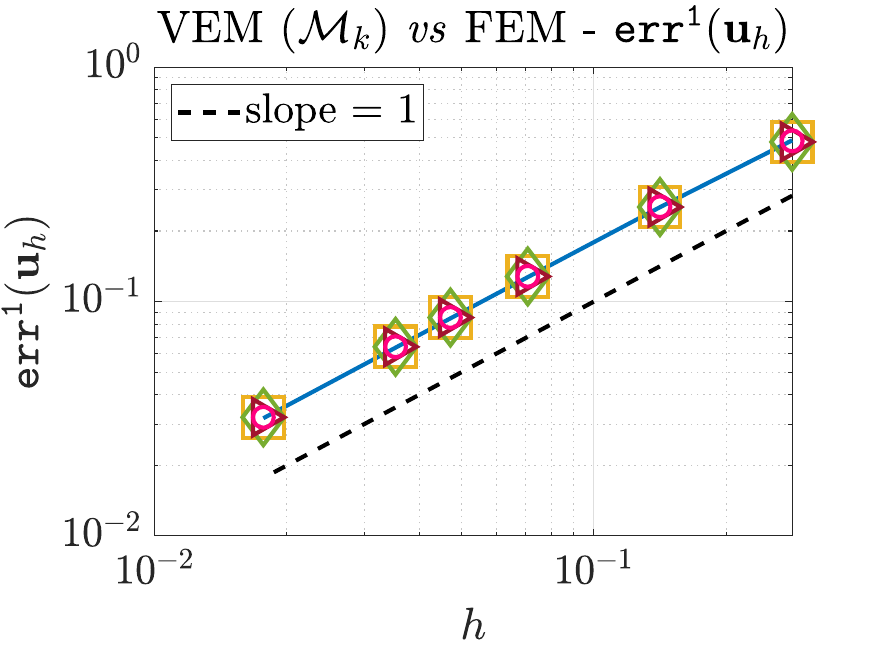}
			\includegraphics[width=0.32\linewidth]{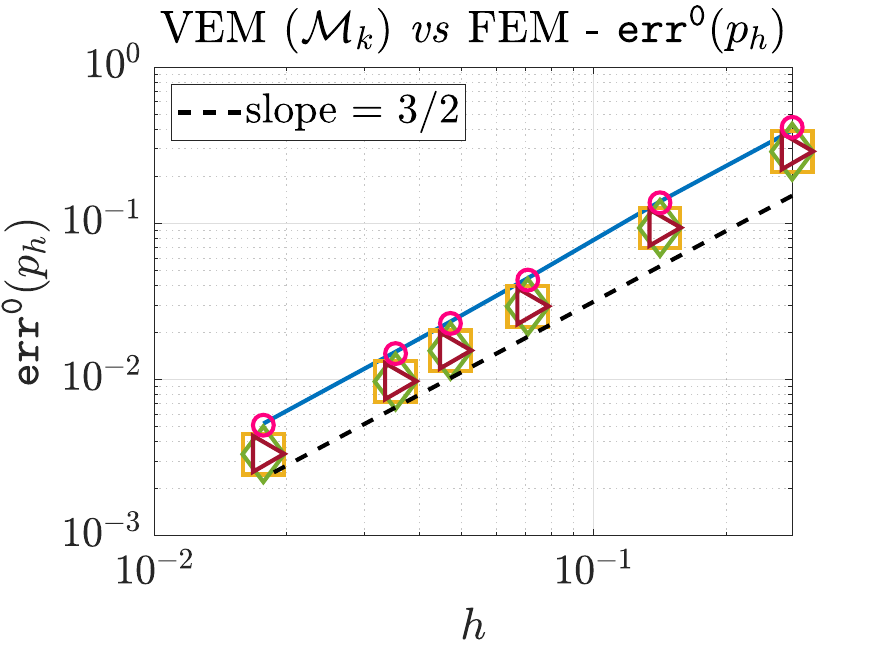}
			
			\
			
			\includegraphics[width=0.32\linewidth]{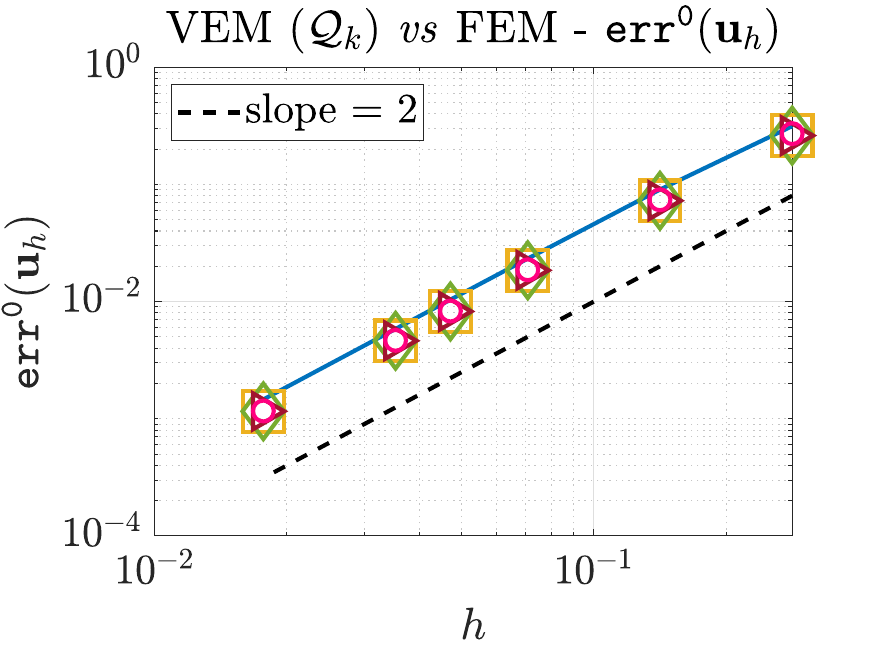}
			\includegraphics[width=0.32\linewidth]{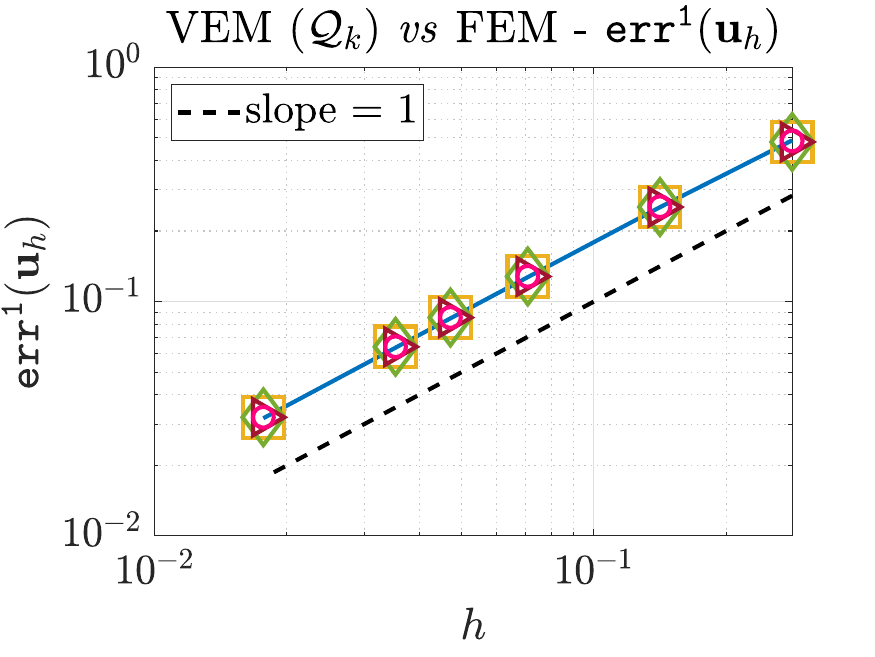}
			\includegraphics[width=0.32\linewidth]{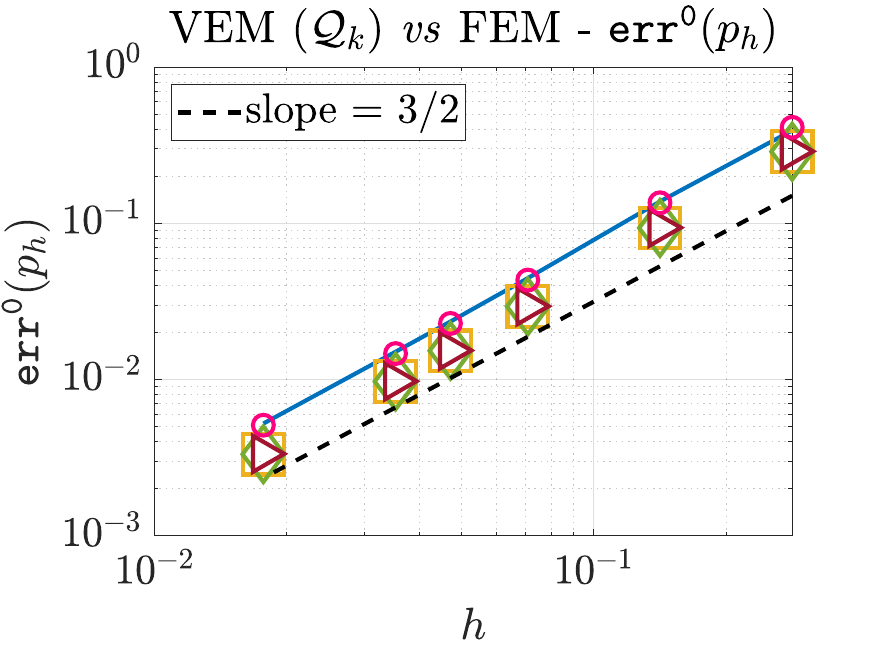}
			
			\
			
			\includegraphics[width=0.75\linewidth]{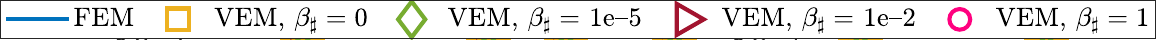}
			
			\caption{\rv{Error comparison between MINI--VEM and MINI--FEM. Top line: $\mathcal{M}_k$; bottom line: $\mathcal{Q}_k$}.}
			\label{fig:fem_err}
		\end{figure}
		
	}
	
	\section{{Static condensation of bubbles}}\label{sec:static_cond}
	
	It is possible to relate our MINI--VEM mixed method with an equal-order virtual element method similar to the one proposed in~\cite{guo2020new,li2024stabilized}. Since bubble functions vanish on each element boundary, they can be eliminated from the discrete system, lowering the global number of degrees of freedom. This process is known as \textit{static condensation}~\cite{franca1996approximation,qu2004static}. We are going to carry out this process in algebraic framework.
	
	We first rewrite Problem~\ref{pro:stokes_discreteC} by isolating the bubble contribution. As previously done at local level, a function $\v_h$ belonging to the global velocity space $\V_k$ can be expressed as the sum of $\vtilde\in\enhancedunoG$ plus a bubble function
	$$
	\bollatest_h\in\bigoplus_{\elementvem\in\mesh_h}\bubbleS.
	$$
	In other words,
	\[
	\bollatest_h = \sum_{\elementvem\in\mesh_h} \bollatest_h^\elementvem \qquad \text{ with } \qquad \bollatest_h^\elementvem\in\bubbleS.
	\]
	We reformulate the discrete Problem~\ref{pro:stokes_discreteC} as follows.
	
	\begin{pb}\label{pb:pb_split}
		Find $\utilde\in\enhancedunoG$, $\bolla_h\in\bigoplus_{\elementvem}\bubbleS$, $p_h\in Q_k$ such that
		\begin{equation*}
			\begin{aligned}
				&a_h^\vmatr(\utilde,\vtilde) - b_h(\vtilde,p_h) = (\f_h,\vtilde)_\Omega &&\forall\vtilde\in\enhancedunoG,\\
				&a_h^{\bolla}(\bolla_h,\bollatest_h) - b_h(\bollatest_h,p_h) = (\f_h,\bollatest_h)_\Omega  &&\forall\bollatest_h\in\bigoplus_{\elementvem}\bubbleS,\\
				&b_h(\utilde+\bolla_h,q_h) + \alpha\,c_h(p_h,q_h)= 0 &&\forall q_h\in Q_k, 
			\end{aligned}
		\end{equation*}
		where
		\begin{gather*}
			a_h^\vmatr(\utilde,\vtilde) = \sum_{\elementvem\in\mesh_h}a_h^{\elementvem,\vmatr}(\utilde,\vtilde), \qquad \text{} \qquad a_h^{\bolla}(\bolla_h^\elementvem,\bollatest_h^\elementvem) =\sum_{\elementvem\in\mesh_h}a_h^{\elementvem,\bolla}(\bolla_h^\elementvem,\bollatest_h^\elementvem),\\[2pt]	
			\begin{aligned}
				&a_h^{\elementvem,\vmatr}(\utilde,\vtilde) = a^\elementvem(\Pinabla{k}\utilde,\Pinabla{k}\vtilde) + S_k^\elementvem(\Qnabla{k}\utilde,\Qnabla{k}\vtilde), \\[2pt]
				&a_h^{\elementvem,\bolla}(\bolla_h^\elementvem,\bollatest_h^\elementvem)=a^\elementvem(\Pinablakpu\bolla_h,\Pinablakpu\bollatest_h)+ \beta_\sharp\,\Skpu(\Qnablakpu\bolla_h,\Qnablakpu\bollatest_h).
			\end{aligned}	
		\end{gather*}
	\end{pb}
	
	Clearly, the following property is satisfied: $a_h(\u_h,\v_h) = a_h^\vmatr(\utilde,\vtilde) + a_h^{\bolla}(\bolla_h^\elementvem,\bollatest_h^\elementvem)$. 
	
	Since bubbles functions are defined on each element $\elementvem\in\mesh_h$ and vanish at the boundary $\partial\elementvem$, the second equation in Problem~\ref{pb:pb_split} can be localized. By exploiting that $$
	b_h^\elementvem(\bollatest_h^\elementvem,p_h)=-(\bollatest_h^\elementvem,\grad\Pizero{k}p_h)_\elementvem,
	$$
	it is clear that each local bubble $\bolla_h^\elementvem$ solves a discrete Poisson equation in $\elementvem$.
	\begin{pb}[Local Bubble problem]\label{pb:local_bubble}
		Find $\bolla_h^\elementvem\in\bubbleS\subset\dueHunozeroK$ such that
		\begin{equation}
			a_h^{\elementvem,\bolla}(\bolla_h^\elementvem,\bollatest_h^\elementvem) = (\f_h - \grad\Pizero{k}p_h,\bollatest_h^\elementvem)_\elementvem\qquad\forall\bollatest_h^\elementvem\in\bubbleS.
		\end{equation}
	\end{pb}
	Thanks to this fact, we can compute explicitly the bubble contribution and plug it into the other two equations of Problem~\ref{pb:pb_split}. 
	
	We first rewrite Problem~\ref{pb:pb_split} in algebraic form.
	
	\begin{pb}\label{pb:operators}
		Find $\utilde\in\enhancedunoG$, $\bolla_h\in\bigoplus_{\elementvem}\bubbleS$, $p_h\in Q_k$ such that
		\begin{equation}
			\begin{aligned}
				&\Amatr_\vmatr\utilde - \Bmatr_\vmatr^\top p_h = \Fmatr_\vmatr,\\
				&\Amatr_\bolla\bolla_h - \Bmatr_\bolla^\top p_h = \Fmatr_\bolla,\\
				&\Bmatr_\vmatr\utilde + \Bmatr_\bolla \bolla_h + \Cmatr_\pmatr p_h = \mathbf{0}.
			\end{aligned}
		\end{equation}
	\end{pb}
	
	We consider the second equation in Problem~\ref{pb:operators} and we isolate the bubble $\bolla_h$; we obtain
	\begin{equation}\label{eq:bubble_expression}
		\bolla_h = \Amatr_\bolla^{-1}\Fmatr_\bolla + \Amatr_\bolla^{-1}\Bmatr_\bolla^\top p_h.
	\end{equation}	
	The stiffness matrix $\Amatr_\bolla$, constructed on the bubble space, shows a block structure with each block representing the local Problem~\ref{pb:local_bubble} in an element $\elementvem\in\mesh_h$. Therefore, $\Amatr_\bolla$ is nonsingular.
	
	We notice that the first equation in Problem~\ref{pb:operators} does not contain bubble contribution, therefore we keep it as it is. On the other hand, we take the third equation and we plug in the expression~\eqref{eq:bubble_expression}. We find
	\begin{equation}
		\Bmatr_\vmatr\utilde + (\Bmatr_\bolla\Amatr_\bolla^{-1}\Bmatr_\bolla^\top+ \Cmatr_\pmatr)p_h = -\Bmatr_\bolla\Amatr_\bolla^{-1}\Fmatr_\bolla.
	\end{equation}
	
	Finally, we state the condensed problem.
	\begin{pb}\label{pb:pb_condensed}
		Find $\utilde\in\enhancedunoG$, $p_h\in Q_k$ such that
		\begin{equation}
			\begin{aligned}
				&\Amatr_\vmatr\utilde - \Bmatr_\vmatr^\top p_h = \Fmatr_\vmatr,\\
				&\Bmatr_\vmatr\utilde + (\Bmatr_\bolla\Amatr_\bolla^{-1}\Bmatr_\bolla^\top+ \Cmatr_\pmatr)p_h = -\Bmatr_\bolla\Amatr_\bolla^{-1}\Fmatr_\bolla.
			\end{aligned}
		\end{equation}
	\end{pb}
	
	Notice that the pressure stabilization consists now of two terms: $\Cmatr_\pmatr$ is the ``original'' one, dealing with the nonpolynomial part, whereas the new term $\Bmatr_\bolla\Amatr_\bolla^{-1}\Bmatr_\bolla^\top$ stabilizes the polynomial contribution and naturally descends from the construction of the MINI--VEM. The condensed Problem~\ref{pb:pb_condensed} is a stabilized equal-order virtual element formulation of the Stokes equation. \rosso{Such formulation is strictly related to those in~\cite{guo2020new,li2024stabilized}: two differences are the right hand side of the second equation and the pressure stabilization term. In our case, the bubble condensation process gives a criterion for the choice of the pressure stabilization.}
	
	\begin{figure}
		\centering
		\includegraphics[width=0.4\linewidth]{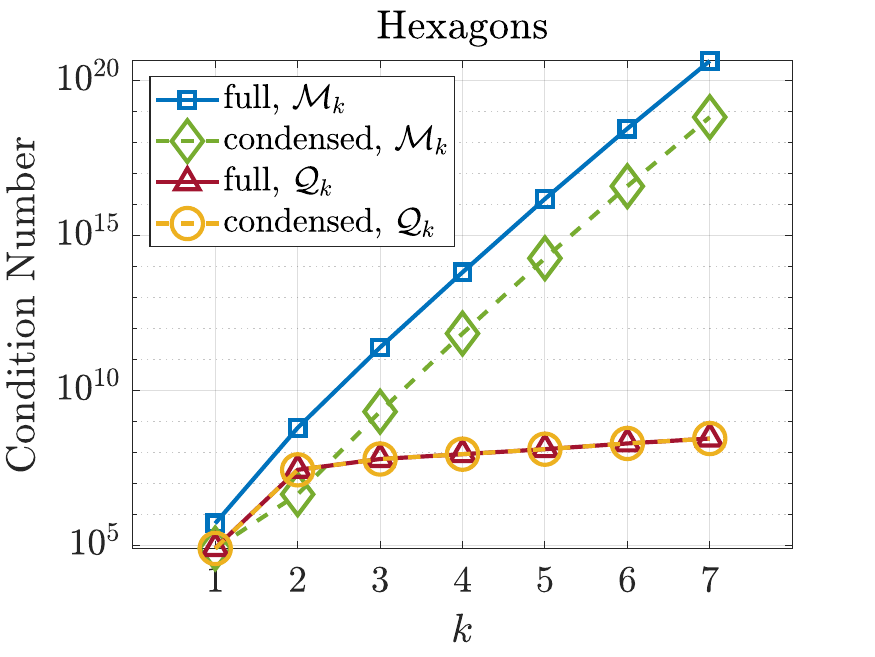}
		\includegraphics[width=0.4\linewidth]{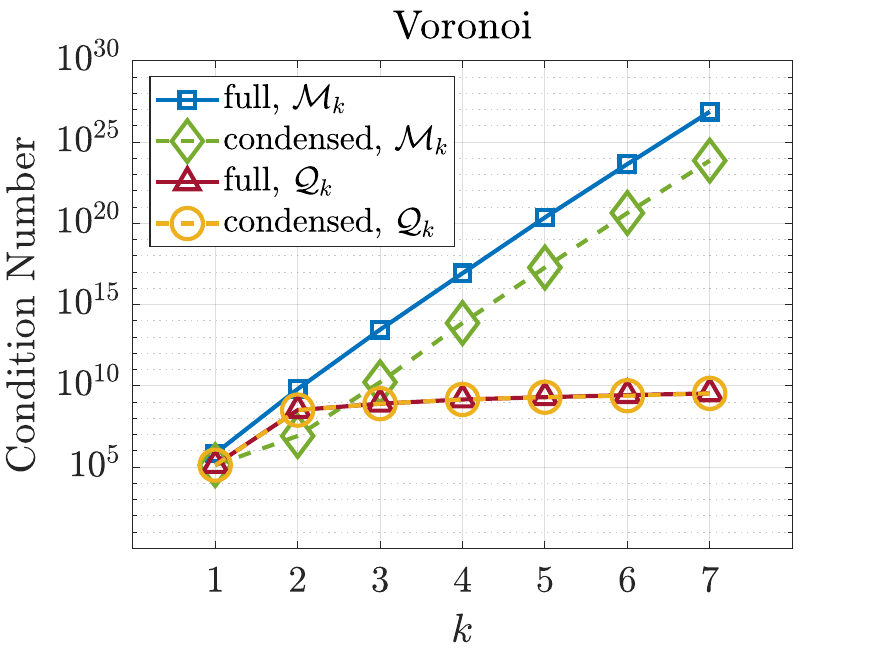}
		
		\
		
		\includegraphics[width=0.4\linewidth]{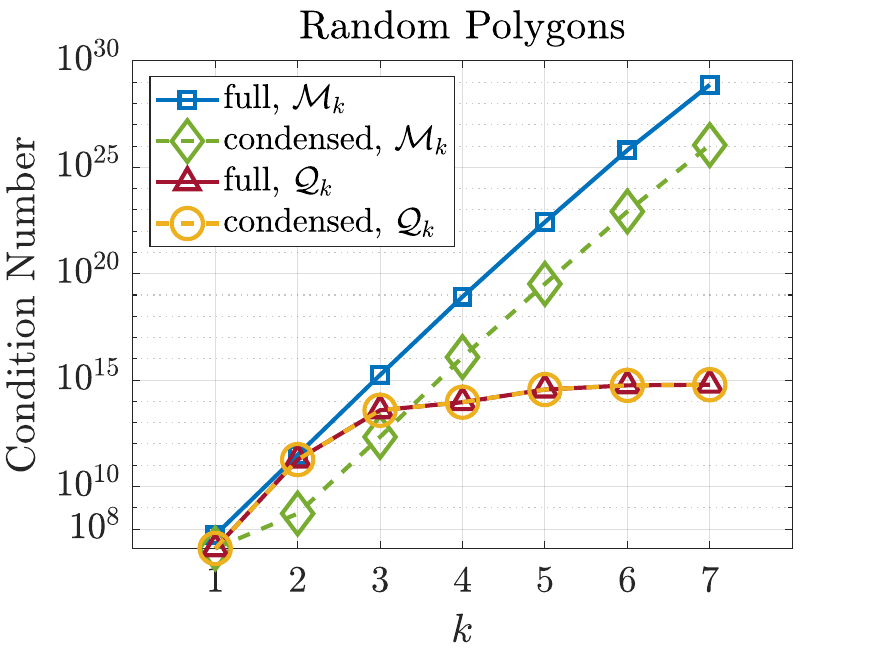}
		\includegraphics[width=0.4\linewidth]{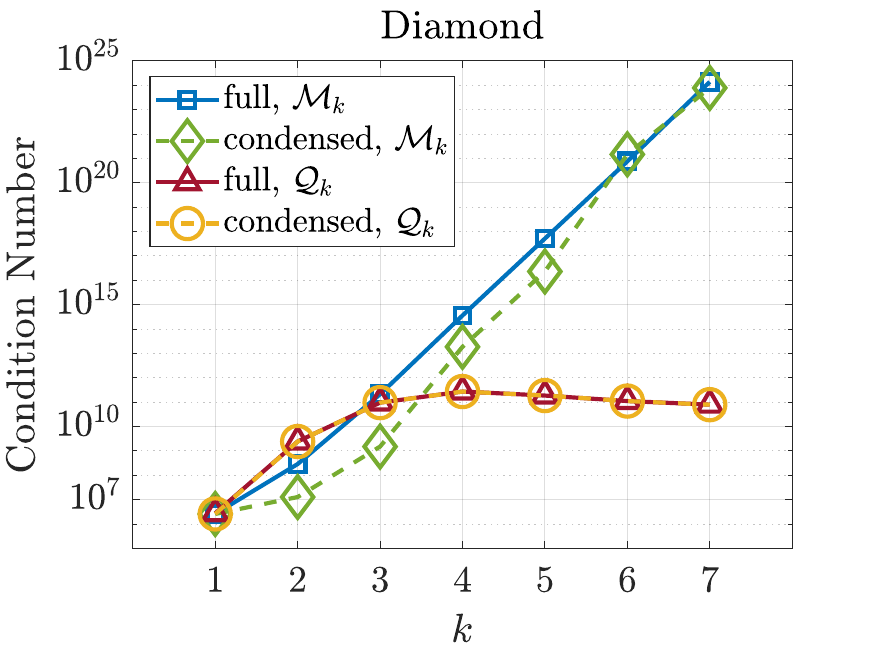}
		
		\caption{\rv{Condition number of the linear system arising from the static condensation of bubbles: comparison with the full problem for both $\mathcal{M}_k$ and $\mathcal{Q}_k$.}}
		\label{fig:conditioning_full_vs_cond}
	\end{figure}
	
\begin{figure}
		\centering
		\includegraphics[width=0.4\linewidth]{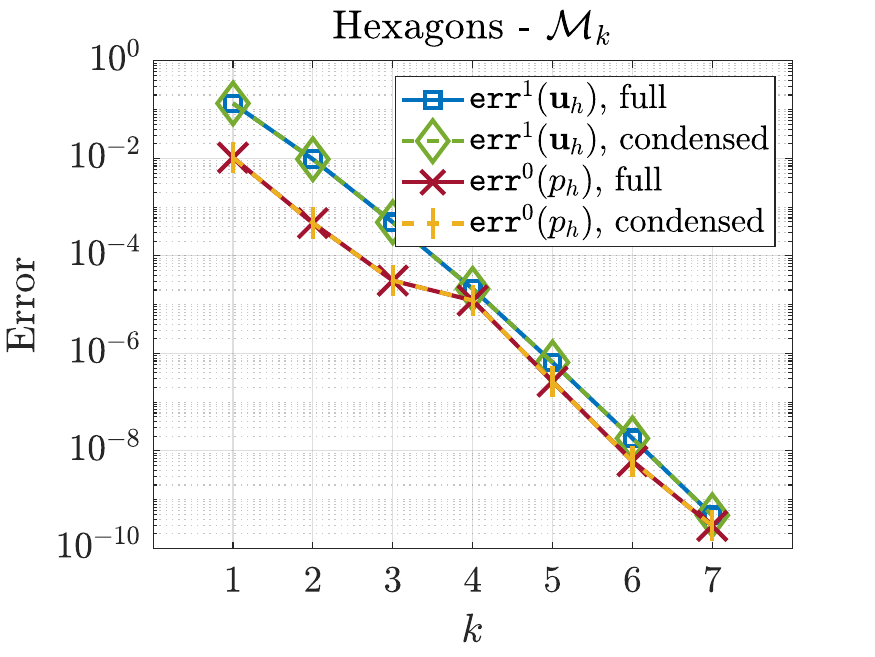}\qquad
		\includegraphics[width=0.4\linewidth]{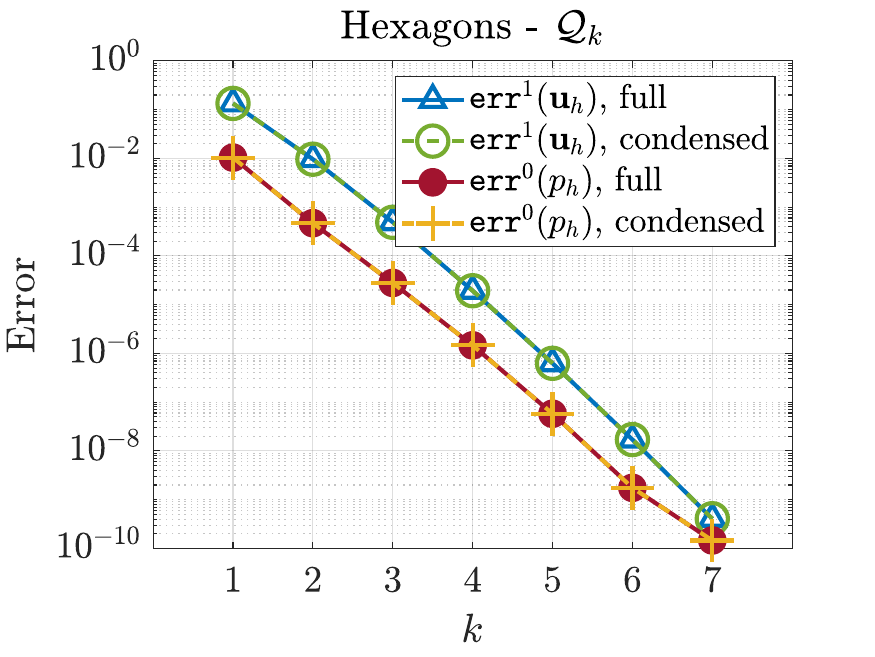}
		
		\
		
		\includegraphics[width=0.4\linewidth]{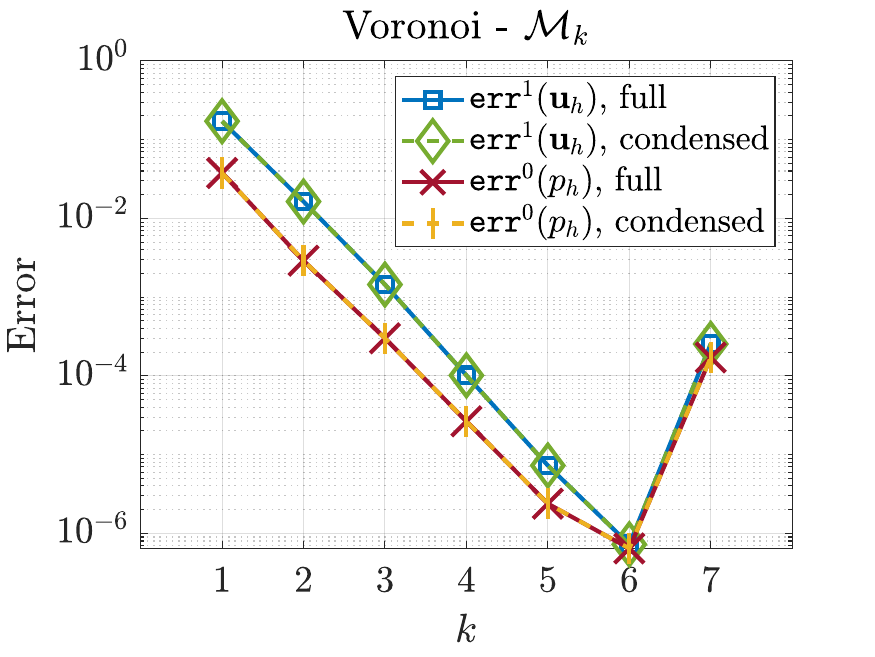}\qquad
		\includegraphics[width=0.4\linewidth]{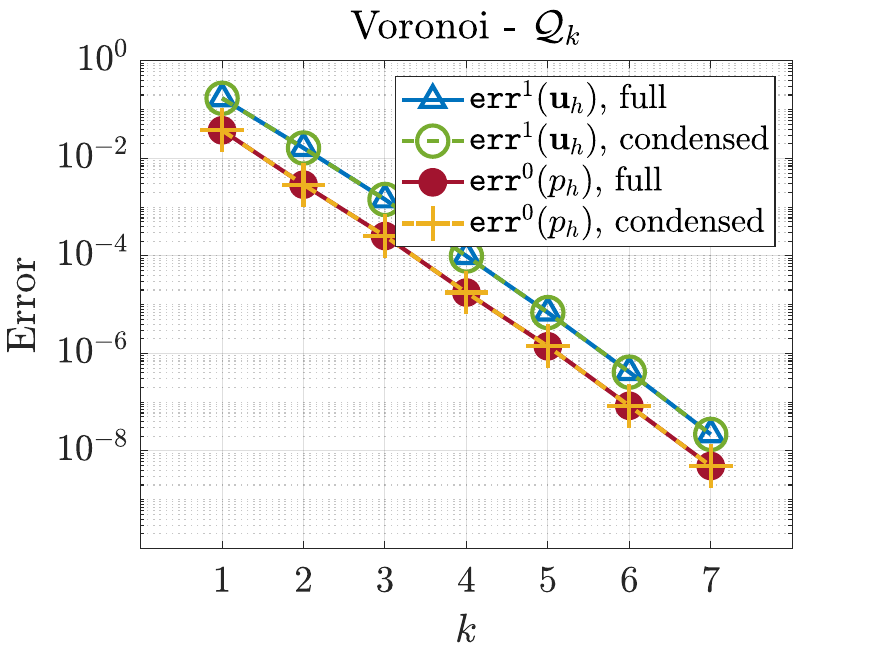}
		
		\
		
		\includegraphics[width=0.4\linewidth]{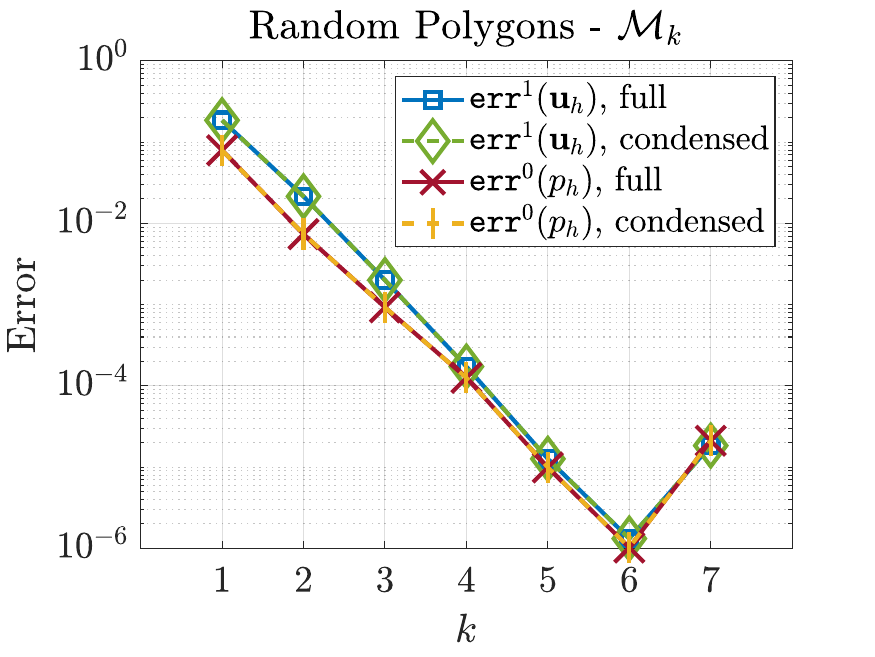}\qquad
		\includegraphics[width=0.4\linewidth]{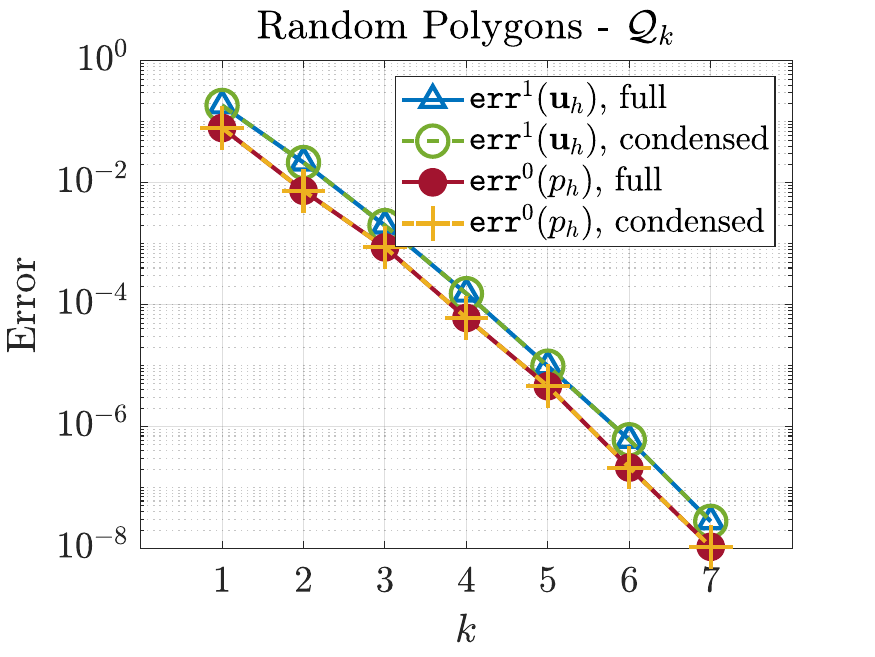}
		
		\
		
		\includegraphics[width=0.4\linewidth]{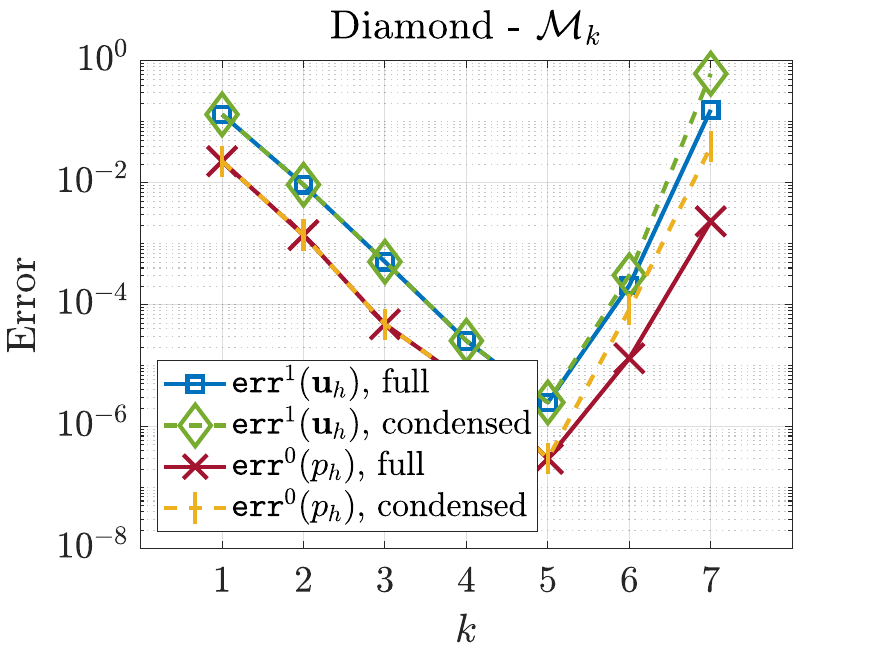}\qquad
		\includegraphics[width=0.4\linewidth]{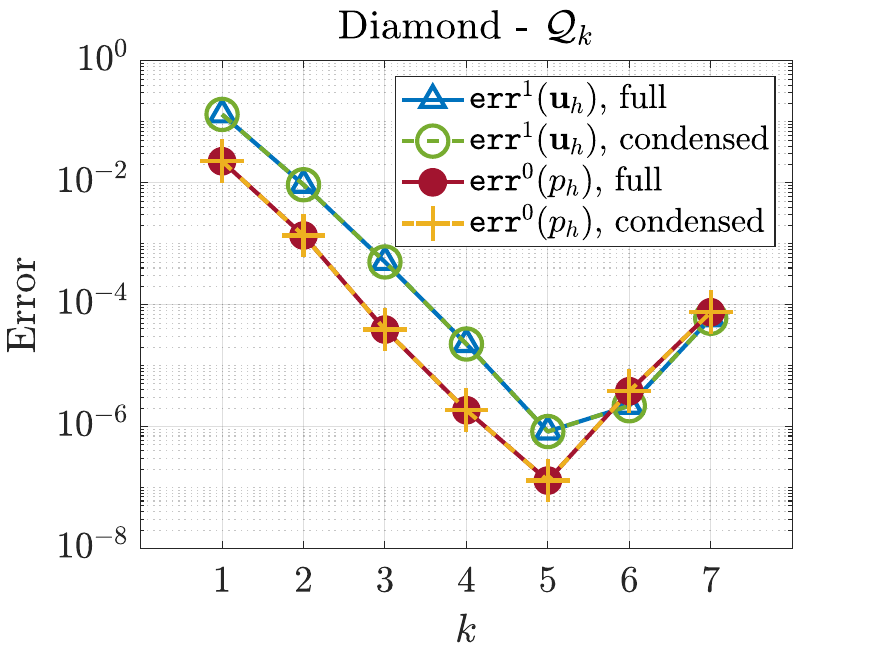}
		\caption{\rv{Error comparison, in terms of $k$, between full and condensed problem. Both $\mathcal{M}_k$ and $\mathcal{Q}_k$ are considered.}}
		\label{fig:err_full_vs_cond}
	\end{figure}

	\rv{We now carry out a comparison between full and condensed problem in terms of condition number and approximation. We consider again the four meshes depicted in Figure~\ref{fig:meshes}, and we vary the degree $k$ from 1 to 7. We plot the evolution of  the condition number in Figure~\ref{fig:conditioning_full_vs_cond} by considering both $\mathcal{M}_k$ and $\mathcal{Q}_k$: there are no significant differences between full and condensed problem, but, for $\mathcal{M}_k$, conditioning of the condensed system is slightly better than that of the full system. The related error analysis is performed on Test 1 by comparing the error evolution for the velocity (in $H^1$ norm) and pressure, as shown in Figure~\ref{fig:err_full_vs_cond}. Also in this case, the results provided by full and condensed problem are equivalent. As already previously observed, conditioning dominates the error when high order approximations are considered on badly shaped meshes, especially when the chosen polynomial basis is $\mathcal{M}_k$.}
	
\begin{figure}
		\centering
		\includegraphics[width=0.325\linewidth]{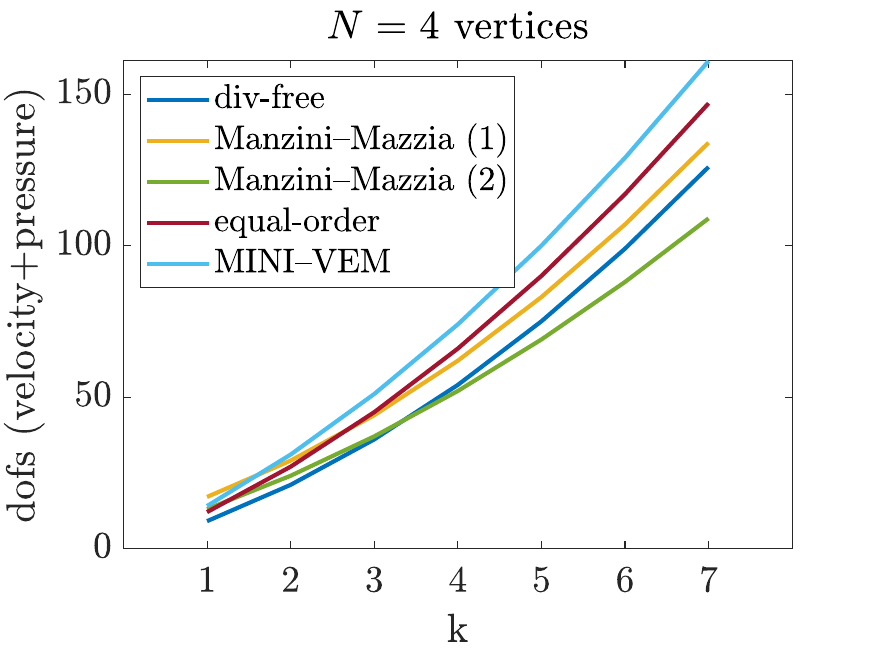}
		\includegraphics[width=0.325\linewidth]{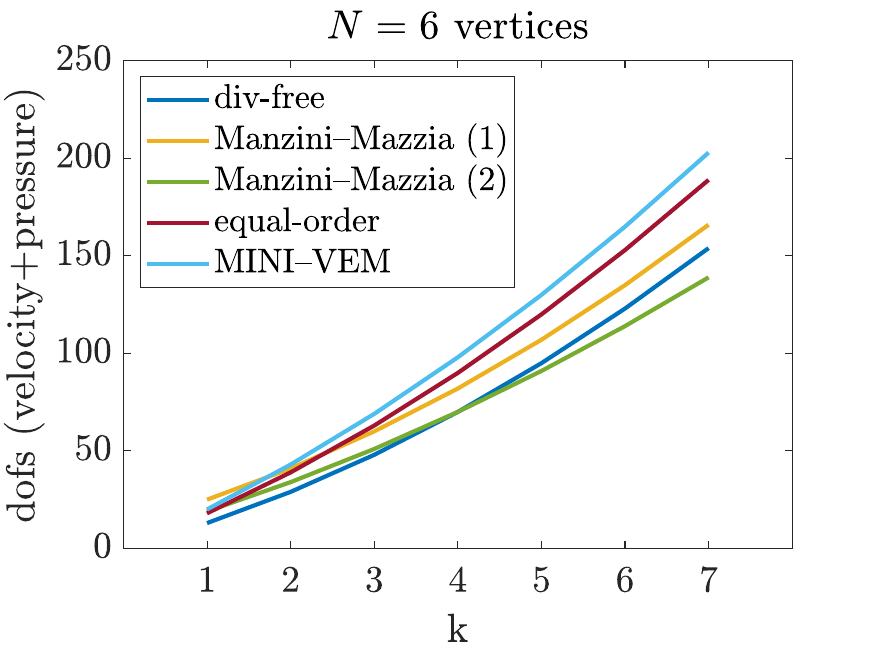}
		\includegraphics[width=0.325\linewidth]{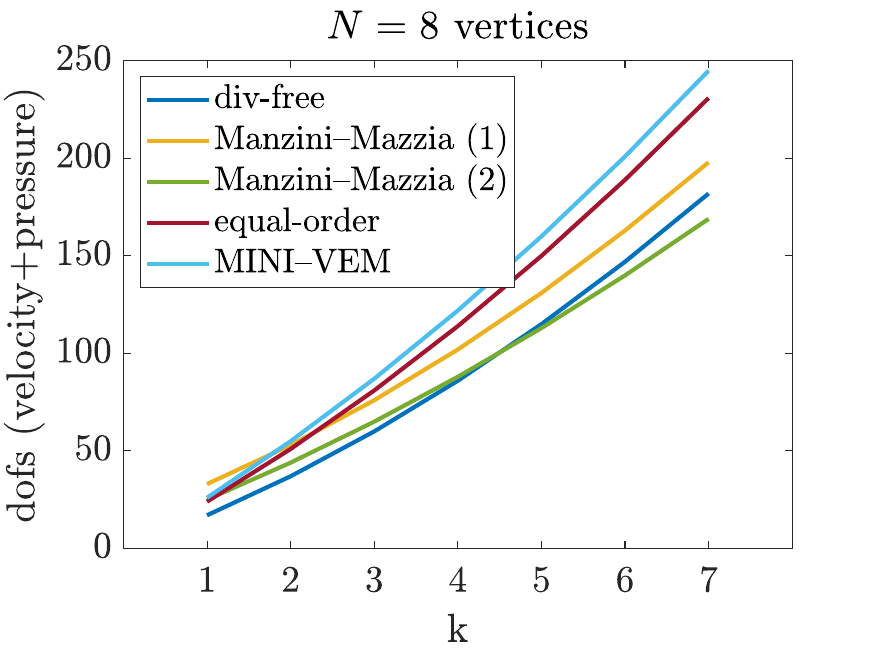}
		\caption{\rv{Comparison, in terms of degrees of freedom, between MINI--VEM and existing virtual element methods for the Stokes equation.}}
		\label{fig:dofs_comparison}
	\end{figure}

	\begin{rem}
		\rv{The static condensation of bubbles removes $2k$ internal moments per element, regardless the number of vertices. Thus, for a polygon $\elementvem$ with $N$ vertices, the total number (velocity and pressure) of local degrees of freedom of an equal-order discretization is
			$$
			3\left(kN+\frac{k(k-1)}{2}\right).
			$$
			On the other hand, the two formulations introduced by Manzini--Mazzia~\cite{MANZINI-MAZZIA} have the following total number of local dofs
			$$
			2\left( N+kN+\frac{k(k-1)}2 \right) + \frac{k(k+1)}2,
			\qquad
			2N+kN+(k-1)N+\frac{k(k-1)}2  + \frac{k(k+1)}2,
			$$
			respectively. Finally, the div-free method~\cite{da2017divergence} requires
			$$
			2Nk + \frac{(k-1)(k-2)}2 + \frac{k(k+1)}2 - 1 
			$$
			local degrees of freedom. As depicted in Figure~\ref{fig:dofs_comparison} for $N=4,6,8$, the MINI--VEM and the equal-order methods have a larger number of local degrees of freedom. This is not surprising since the div-free formulation and those proposed by Manzini and Mazzia consider the velocity as a single object, while MINI--VEM and the equal-order methods are constructed component-wise. Moreover, in~\cite{MANZINI-MAZZIA,da2017divergence} the pressure is discretized by pure polynomials. In any case, the MINI--VEM, in both its full and condensed version, can be easily implemented since it is based on the original virtual element method proposed for the Laplace equation~\cite{beirao2013basic}.}
	\end{rem}
	
	\section{Conclusions}\label{sec:conclusions}
	
	In this paper we presented the MINI--VEM mixed method for the Stokes equation as natural evolution of the popular MINI mixed finite element method for polygonal meshes and high order approximation.
	
	Both velocity and pressure spaces are defined starting from the enhanced virtual element space of degree $k$ presented in~\cite{ahmad2013equivalent}. The velocity space is then enriched with virtual bubble functions of degree $\rvf{k+1}$ which, together with a pressure stabilization, ensure the well-posedness of the discrete method. We proved optimal error estimates for velocity and pressure in energy norm and for the velocity in $L^2$ norm, which were then confirmed by several numerical tests. We also investigated how the choice of polynomial basis and the value of pressure stabilization parameter affect the condition number of the linear system arising from the MINI--VEM discretization. From the MINI--VEM formulation we derived an equal-order virtual element method by exploiting static condensation of the additional bubble functions.
	
	Within the paper, we also proved that virtual bubble functions are self-stabilized and orthogonal to harmonic polynomials with respect to the product of gradients.

	\section{Acknowledgments}
	
	The authors are member of the INdAM -- GNCS research group. 
		
	This paper has been co--funded by the MUR Progetti di Ricerca di Rilevante Interesse Nazionale (PRIN) Bando 2020 (grant 20204LN5N5) and Bando 2022/PNRR (grant P2022BH5CB, NextGenerationEU) .

\bibliographystyle{abbrv}
\bibliography{biblio}

\end{document}